\documentclass[a4paper,oneside,10pt]{article}%
\usepackage{amsmath}
\usepackage{amsfonts}
\usepackage{amssymb}
\usepackage{multirow}
\usepackage{graphicx}
\usepackage{epsfig}
\usepackage{booktabs}
\usepackage[square,numbers,sort&compress]{natbib}
\usepackage{diagbox}%
\usepackage[pagewise]{lineno}
\usepackage{caption}
\usepackage[UTF8, scheme=plain]{ctex}
\usepackage{xcolor}
\usepackage{float}
\usepackage{placeins}
\usepackage{needspace}

\usepackage{enumitem}
\usepackage{hyperref}

\usepackage[ruled,linesnumbered]{algorithm2e}
\SetAlgoLined
\SetKwInput{KwIn}{Input}
\SetKwInput{KwOut}{Output}
\SetKw{KwRet}{return}
\SetKwComment{Comment}{$\triangleright$\ }{}

\providecommand{\U}[1]{\protect\rule{.1in}{.1in}}
\allowdisplaybreaks[4]

\newtheorem{theorem}{Theorem}[section]

\newtheorem{lemma}[theorem]{Lemma}

\newtheorem{proposition}[theorem]{Proposition}
\newtheorem{remark}[theorem]{Remark}

\newenvironment{proof}[1][Proof]{\noindent \textbf{#1.} }{\  \rule{0.5em}{0.5em}}
\numberwithin{equation}{section}
\renewcommand \arraystretch{1.5}
\begin{document}

\title{Learning the Brownian Motion in Machine Learning for Parametric Parabolic PDEs: A Backward Quadrature Approach} 

\author{Qiang Han \thanks{School of Mathematics, Yangzhou University, Yangzhou 225002, PR China. Email: hanqiang@yzu.edu.cn.}
\and Shaolin Ji\thanks{Zhongtai Securities Institute for Financial Studies,
Shandong University, Jinan, Shandong 250100, PR China. Email: jsl@sdu.edu.cn}
\and Yanan Qi\thanks{Corresponding author. 
Zhongtai Securities Institute for Financial Studies,
Shandong University, Jinan, Shandong 250100, PR China. Email: 202112024@sdu.edu.cn}
}

\maketitle

\textbf{Abstract}.
We propose a backward quadrature learning (BQL) method for solving high dimensional parametric semilinear parabolic partial differential equations (PDEs).  To the best of our knowledge, this work is the first to extend the applicability of BSDE-based stochastic numerical methods to this class of PDEs.
Specifically, inspired by the approach of learning the random variables (LRV)~\cite{BeckerJentzenMullerWurstemberger2024}, BQL incorporates learnable quadrature for Brownian increments at each time level into the backward recursion of BSDEs, enabling  more efficient layerwise approximation of local conditional expectations. 
Furthermore, BQL introduces a continuation surrogate to facilitate the transfer of value-function information from subsequent to preceding time levels, avoiding nested quadrature while controlling the propagation of local approximation errors.
Our main theoretical contribution is the derivation of a global error bound that explicitly characterizes the decomposition of the total error, together with a proof of first-order mean-square convergence under suitable approximation and integrability assumptions.
Finally, numerical experiments demonstrate that BQL achieves competitive solution accuracy while substantially reducing the amortized computational cost per query. For the HJB benchmark across the tested dimensions, with $Q=7$ parameter queries, BQL achieves computational cost reductions of approximately 40-fold and 4-fold compared with P$^2$INN~\cite{ChoEtAl2024} and FBSNN~\cite{Raissi2018FBSNN}, respectively.



\addcontentsline{toc}{section}{\hspace*{1.8em}Abstract}

\section{Introduction}
\subsection{Problem setting and methodological development}
Parametric PDEs arise naturally in financial engineering, uncertainty quantification, optimal design, and many other applications in which model coefficients, boundary conditions, or physical properties vary over prescribed ranges \cite{CohenDeVore2015,BennerGugercinWillcox2015,QuarteroniManzoniNegri2016}. Consequently, practical computations often require repeated evaluations of PDE solutions or quantities of interest over a parameter domain rather than the solution of a single fixed problem \cite{BennerGugercinWillcox2015,QuarteroniManzoniNegri2016}. Such evaluations can be time-sensitive: in derivative pricing and risk management, changing market parameters require option values and risk sensitivities to be updated promptly \cite{GlauWunderlich2022}, while in engineering design and control, varying material properties or operating conditions require repeated evaluations of system responses \cite{BennerGugercinWillcox2015,QuarteroniManzoniNegri2016}. Efficient repeated evaluation is therefore essential in many-query settings. This motivates the development of solvers that can reuse information across parameter configurations and deliver rapid solutions over the entire prescribed parameter domain.

To make the problem setting concrete, we consider the following high dimensional parametric semilinear parabolic PDEs:
\begin{equation}
\left\{
\begin{aligned}
&\partial_t u^p+\mu_p\cdot D_xu^p+\frac{1}{2}\operatorname{Tr}\!\left(\sigma_p\sigma_p^{\top}D_x^2u^p\right)+f_p\!\left(\cdot,u^p,\sigma_p^{\top}D_xu^p\right)=0 \qquad \text{on }[0,T)\times\mathbb R^d,\\
&u^p(T,\cdot)=g_p(\cdot) \qquad \text{on }\mathbb R^d.
\end{aligned}
\right.
\label{eq:parametric_pde}
\end{equation}
Here, $p\in\mathcal P\subset\mathbb R^{d_p}$ is the model parameter, and $u^p$ is the corresponding solution. The drift and diffusion coefficients are given by $\mu_p:[0,T]\times\mathbb R^d\to\mathbb R^d$ and $\sigma_p:[0,T]\times\mathbb R^d\to\mathbb R^{d\times d}$. The function $f_p:[0,T]\times\mathbb R^d\times\mathbb R\times\mathbb R^d\to\mathbb R$ is the nonlinear generator, and $g_p:\mathbb R^d\to\mathbb R$ specifies the terminal condition. We use $D_x$ and $D_x^2$ for the spatial gradient and Hessian, $\cdot$ for the Euclidean inner product, ${}^\top$ for the transpose, and $\operatorname{Tr}$ for the matrix trace.

A straightforward approach to solving \eqref{eq:parametric_pde} is to treat each parameter configuration independently \cite{CohenDeVore2015,BennerGugercinWillcox2015,QuarteroniManzoniNegri2016}. However, this requires a separate numerical solve or model training for each parameter value, resulting in high computational costs in many-query settings \cite{BelbutePeresEtAl2021,ChoEtAl2024}.
To address this issue, classical parametric methods mainly include Galerkin-based reduced order methods, such as reduced basis (RB) methods and proper orthogonal decomposition (POD), with empirical interpolation techniques (EIM/DEIM) to accelerate nonlinear operator evaluation \cite{BennerGugercinWillcox2015,QuarteroniManzoniNegri2016,BarraultMadayNguyenPatera2004,ChaturantabutSorensen2010}, as well as stochastic methods based on polynomial expansions or parameter sampling, such as generalized polynomial chaos (GPC), stochastic Galerkin methods, and sparse-grid stochastic collocation \cite{XiuKarniadakis2002,BabuskaTemponeZouraris2004,BabuskaNobileTempone2007,NobileTemponeWebster2008,SchwabGittelson2011}. However, for high dimensional nonlinear parametric parabolic PDEs, these methods may incur substantial computational costs due to high fidelity simulations and nonlinear operator approximation, particularly when the parameter-to-solution map lacks sufficient regularity or sparsity for efficient approximation in the parameter space \cite{GreplMadayNguyenPatera2007,DrohmannHaasdonkOhlberger2012,BenaceurEhrlacherErnMeunier2018,CohenDeVore2015,SchwabGittelson2011,NobileTemponeWebster2008}. These limitations motivate alternative methods that directly approximate nonlinear parameter-to-solution maps without relying on prescribed low dimensional linear approximation spaces.

To address the limitations of traditional parametric methods in approximating complex nonlinear solution maps, neural networks have increasingly been used to construct unified solution representations over parameter domains, offering flexible nonlinear approximation for high dimensional parametric PDEs \cite{KutyniokPetersenRaslanSchneider2022,BhattacharyaHosseiniKovachkiStuart2021}. One class of methods directly learns the parameter-to-solution map by taking model parameters, time, and state variables as network inputs. Representative approaches include a Feynman-Kac-based statistical learning framework for high dimensional parametric Kolmogorov equations \cite{BernerDablanderGrohs2020}, a Monte Carlo(MC) and deep learning approach for jointly approximating PDE solutions and their gradients \cite{SabateVidalesSiskaSzpruch2022}, and the deep parametric PDE method, which incorporates parameter information into a unified neural network \cite{GlauWunderlich2022}. Within the physics-informed neural network (PINN) framework, HyperPINN,  P$^2$INN, and related conditioning approaches further address parameter conditioning and generalization across parameter configurations \cite{BelbutePeresEtAl2021,ChoEtAl2024,KassaiKoupaiEtAl2024}. Another class learns PDE solution operators as mappings between function spaces, including deep operator network (DeepONet), Fourier neural operator (FNO), and their physics-informed variants, such as PI-DeepONet and PINO \cite{LuEtAl2021DeepONet,LiEtAl2021FNO,WangWangPerdikaris2021,LiEtAl2024PINO}.

Although these methods enable information reuse across parameter configurations, achieving the desired accuracy may still incur substantial training or data generation costs. Directly parameterized neural networks must learn a unified solution map over the joint parameter-time-state domain, with performance depending on training sample coverage and network optimization. Parameterized PINNs additionally require repeated evaluations of PDE residuals and their derivatives at numerous collocation points, making their training costs dependent on both the number of collocation points and optimization iterations \cite{PenwardenZheNarayanKirby2023,BelbutePeresEtAl2021,ChoEtAl2024}. Data-driven neural operators require high fidelity datasets covering different parameter configurations or input functions, often involving multiple PDE solves. Although physics-informed operator learning reduces the need for paired solution data, it still entails evaluating physics-based constraints and iterative optimization \cite{LuEtAl2021DeepONet,LiEtAl2021FNO,WangWangPerdikaris2021,LiEtAl2024PINO}.

These limitations motivate alternative approaches that exploit the probabilistic structure of PDEs rather than directly learning solution maps or operators. For the high dimensional parametric semilinear parabolic PDEs considered here, the nonlinear Feynman-Kac formula establishes a connection with BSDEs, providing a probabilistic representation of their solutions \cite{PardouxPeng1992}. Compared with directly learning the parametric solution map over a high dimensional joint domain, BSDE-based probabilistic methods rely on stochastic trajectory simulation and conditional expectation evaluation. They avoid explicit high dimensional spatial grids and facilitate parallel computation through independent trajectory simulations. Moreover, the temporal structure of BSDEs enables backward schemes to decompose the computation into local tasks across successive time levels, providing a structured approach to high dimensional problems.
Based on this connection, deep BSDE methods combine stochastic trajectory simulation with neural network approximations to solve high dimensional PDEs \cite{EHanJentzen2017,HanJentzenE2018}, while deep backward dynamic programming (DBDP) and its multistep extension (MDBDP) exploit the temporal recursive structure of BSDEs to perform successive backward approximations at discrete time levels \cite{HurePhamWarin2020,GermainPhamWarin2022}. In addition, BSDE-based stochastic methods have also been applied to parametric linear parabolic PDEs, combining martingale representations with deep neural networks to jointly approximate PDE solutions and their gradients across parameter configurations \cite{SabateVidalesSiskaSzpruch2022}.

However, many parametric PDEs arising in practical applications also contain nonlinear terms. For example, in credit valuation adjustment~\cite{GermainPhamWarin2023} variations in parameters such as default intensity and volatility affect the nonlinear valuation equation. For such high dimensional problems, extending BSDE-based backward solvers remains challenging in three aspects: (i) \textit{Accuracy and efficiency of local conditional expectation approximation.} The nonlinear dependence of the generator on the unknown processes $Y$ and $Z$ increases the difficulty of approximating conditional expectations, making it challenging to balance approximation accuracy and computational efficiency; (ii) \textit{Accuracy and efficiency of information transfer across time levels.} Efficiently transferring value function information from subsequent time levels while controlling the propagation and accumulation of local approximation errors remains challenging; (iii) \textit{Training and generalization over the joint parameter-state space.} Introducing model parameters enlarges the input domain, increasing the difficulty of training sample coverage and model optimization, while imposing higher requirements on generalization across parameter configurations.

To address these challenges, inspired by the LRV approach \cite{BeckerJentzenMullerWurstemberger2024}, we explore the incorporation of learnable stochastic approximation into the time-discrete backward recursion of BSDEs. LRV introduces a novel learning paradigm for stochastic numerical approximation by treating realizations of random variables in MC estimators as trainable parameters, rather than relying solely on random sampling to determine these realizations. This idea is further extended to Brownian-motion-driven numerical approximations, where discretized Brownian increments are treated as learnable variables. Moreover, the combination of LRV with the Multilevel Picard (MLP) method demonstrates the potential of this learning paradigm for recursive stochastic estimation. Building on this pioneering idea, we propose BQL, a BSDE-based backward solution framework for high dimensional parametric semilinear parabolic PDEs.

Drawing on the core idea of the LRV learning paradigm, BQL adapts learnable stochastic approximation to the backward recursion of BSDEs. Specifically, starting from the terminal condition, BQL proceeds backward in time and learns separate sets of quadrature nodes for the Brownian increments involved in the $Z$ and $Y$ updates at each time level, which differs from LRV's direct learning of a discrete representation of Brownian motion. These nodes are optimized over joint parameter-state samples to approximate the corresponding local conditional expectations and are shared across different parameter-state configurations after training. Building on this mechanism, BQL further introduces a continuation surrogate to facilitate the transfer of value-function information from subsequent to preceding time levels, avoiding nested quadrature while controlling the propagation of local approximation errors. Overall, learnable quadrature and continuation approximation address local conditional expectation evaluation and information transfer across time levels, respectively, forming the proposed backward solution framework.

\subsection{Overview of this work}

The main contributions of this work are summarized as follows:
\begin{itemize}
\item To the best of our knowledge, we are the first to extend BSDE-based stochastic numerical methods to high dimensional parametric semilinear parabolic PDEs, establishing a new application setting for these methods.

\item By organically integrating the core idea of the LRV approach into the time-discrete backward recursion framework of BSDEs, we develop a learnable stochastic quadrature framework to address three key challenges: balancing approximation accuracy and computational efficiency in local conditional expectation evaluation, improving the efficiency of information transfer across time levels, and enhancing training and generalization over the joint parameter--state domain. Furthermore, we establish a global error bound that characterizes the decomposition of the total error and prove first-order mean-square convergence under suitable approximation and integrability assumptions.

\item Numerical experiments demonstrate that BQL achieves competitive solution accuracy while substantially reducing the amortized computational cost per query. Taking the HJB benchmark as an example, across the tested dimensions, BQL achieves approximately 40-fold and 4-fold reductions in cumulative computational costs compared with P$^2$INN and FBSNN, respectively, with $Q=7$ parameter queries.
\end{itemize}

The remainder of this paper is organized as follows. Section~2 introduces the problem setting for parametric semilinear parabolic PDEs, their probabilistic representation, the time discretization of the associated BSDEs, and the relevant function approximation tools. Section~3 presents the BQL method in detail, including learnable quadrature nodes, sequential backward training, and the continuation surrogate. Section~4 develops the global error analysis, local error bounds, and convergence analysis. Section~5 evaluates the solution accuracy and computational efficiency of the proposed method through numerical experiments on  high dimensional parametric semilinear
parabolic PDEs. Finally, Section~6 concludes the paper.

\section{Preliminaries}
\label{Preliminaries}

\subsection{Probabilistic Representation of Parametric PDEs}

In this section, we briefly review the probabilistic representation of parametric nonlinear parabolic PDEs and the conditional-expectation formulation used in BSDE-based numerical methods.

For each parameter $p\in\mathcal P$, consider the diffusion process $(X_t^p)_{0\le t\le T}$ defined by
\[
X_t^p=X_0^p+\int_0^t\mu_p(s,X_s^p)\,ds+\int_0^t\sigma_p(s,X_s^p)\,dW_s,\qquad 0\le t\le T,
\]
where $W$ is a $d$-dimensional standard Brownian motion defined on a filtered probability space $(\Omega,\mathcal F,(\mathcal F_t)_{0\le t\le T},\mathbb P)$, and $X_0^p$ is an $\mathcal F_0$-measurable random variable in $\mathbb R^d$.

Through the nonlinear Feynman--Kac representation~\cite{PardouxPeng1990,Peng1991,PardouxPeng1992}, the solution $u^p$ of the \eqref{eq:parametric_pde} admits the following BSDE representation:
\[
Y_t^p=g_p(X_T^p)+\int_t^T f_p(s,X_s^p,Y_s^p,Z_s^p)\,ds-\int_t^T Z_s^p\,dW_s,\qquad 0\le t\le T.
\]
Under suitable regularity assumptions,
\[
Y_t^p=u^p(t,X_t^p),\qquad Z_t^p=\sigma_p(t,X_t^p)^\top\nabla_xu^p(t,X_t^p).
\]

To numerically solve the BSDE system, we consider a time discretization $\pi=\{0=t_0<t_1<\cdots<t_N=T\}$, where $\Delta t_i=t_{i+1}-t_i$ and $\Delta W_i=W_{t_{i+1}}-W_{t_i}$ for $i=0,\ldots,N-1$. The forward diffusion process is approximated by the Euler scheme
\[
X_{i+1}^{p,\pi}=X_i^{p,\pi}+\mu_p(t_i,X_i^{p,\pi})\Delta t_i+\sigma_p(t_i,X_i^{p,\pi})\Delta W_i,\qquad i=0,\ldots,N-1.
\]

The backward component is discretized by the standard explicit conditional-expectation Euler scheme. Specifically, for $i=0,\ldots,N-1$, we have~\cite{BouchardTouzi2004, Zhang2004}
\begin{equation}
\begin{cases}
Y_i^{p,\pi}=\mathbb E_i\!\left[Y_{i+1}^{p,\pi}+f_p(t_i,X_i^{p,\pi},Y_{i+1}^{p,\pi},Z_i^{p,\pi})\Delta t_i\right],\\
Z_i^{p,\pi}=\mathbb E_i\!\left[Y_{i+1}^{p,\pi}\dfrac{\Delta W_i}{\Delta t_i}\right],
\end{cases}
\label{eq:YZ_conditional}
\end{equation}
with terminal condition $Y_N^{p,\pi}=g_p(X_N^{p,\pi})$, where $\mathbb E_i[\cdot]=\mathbb E[\cdot\mid\mathcal F_{t_i}]$.

The $Y$-recursion follows by taking conditional expectations in the integrated BSDE and approximating the generator integral. The $Z$-recursion is obtained from the Brownian-increment projection of the integrated BSDE, with the generator-covariation term and the variation of $Z$ over each time interval incorporated into the local time-discretization remainder. Neglecting these remainders yields the stated explicit Euler scheme.

\subsection{Feedforward Neural Networks as Function Approximators}

Classical universal approximation results establish that feedforward neural networks (FNN) with suitable nonlinear activation functions can approximate continuous functions on compact subsets of Euclidean spaces under appropriate assumptions~\cite{Cybenko1989,HornikStinchcombeWhite1989}. In this work, the parameter--state input is represented by the concatenated feature vector $\chi(p,x):=(p^\top,x^\top)^\top\in\mathbb R^{d_\chi}$, where $d_\chi=d_p+d$ when $p\in\mathbb R^{d_p}$. Consider a fully connected FNN with $L$ hidden layers. Setting $d_0=d_\chi$, let $d_1,\ldots,d_L$ denote the hidden-layer widths and $d_{\mathrm{out}}$ the output dimension. With $h_0=\chi(p,x)$, the hidden representations are defined recursively by
\[
h_{\ell+1}=\rho_{\mathrm{act}}(W_\ell h_\ell+b_\ell),\quad \ell=0,\ldots,L-1,
\]
where $W_\ell\in\mathbb R^{d_{\ell+1}\times d_\ell}$ and $b_\ell\in\mathbb R^{d_{\ell+1}}$ denote the weight matrix and bias vector of the $\ell$th hidden layer, respectively, and $\rho_{\mathrm{act}}$ is a componentwise nonlinear activation function. The output layer is defined as
\[
\mathcal M(p,x;\omega)=W_{\mathrm{out}}h_L+b_{\mathrm{out}},
\]
where $W_{\mathrm{out}}\in\mathbb R^{d_{\mathrm{out}}\times d_L}$ and $b_{\mathrm{out}}\in\mathbb R^{d_{\mathrm{out}}}$. The collection $\omega=\{W_\ell,b_\ell\}_{\ell=0}^{L-1}\cup\{W_{\mathrm{out}},b_{\mathrm{out}}\}$ denotes the trainable network parameters.

\subsection{Radial Basis Function Networks as Function Approximators}

Radial basis function (RBF) networks constitute a classical class of function approximators constructed from basis functions associated with a set of centers~\cite{ParkSandberg1991,Buhmann2003}. Using the feature vector $\chi(p,x)$ defined above, a scalar-valued RBF approximation with $N_c$ centers is
\[
\mathcal R(p,x)=\sum_{j=1}^{N_c}w_j\,\phi\!\left(\chi(p,x),c_j\right)+b,
\]
where $c_j\in\mathbb R^{d_\chi}$ denotes the $j$th center, $w_j\in\mathbb R$ its coefficient, and $b\in\mathbb R$ the bias term.

A commonly used choice is the isotropic Gaussian radial basis function~\cite{Buhmann2003},
\begin{equation}
\phi^{\mathrm G}\!\left(\chi(p,x),c;\gamma\right)=\exp\!\left(-\gamma\left\|\chi(p,x)-c\right\|_2^2\right),\qquad \gamma>0,
\label{eq:gaussian_rbf}
\end{equation}
where $\gamma$ is a scalar shape parameter shared by all feature coordinates.

A standard anisotropic Gaussian extension allows the shape parameter to vary across coordinates~\cite{FasshauerHickernellWozniakowski2012}:
\begin{equation}
\phi^{\mathrm{AG}}\!\left(\chi(p,x),c;\boldsymbol{\gamma}\right)=\exp\!\left(-\sum_{r=1}^{d_\chi}\gamma_r\left(\chi_r(p,x)-c_r\right)^2\right),\qquad \gamma_r>0,\quad r=1,\ldots,d_\chi,
\label{eq:anisotropic_rbf}
\end{equation}
where $\boldsymbol{\gamma}=(\gamma_1,\ldots,\gamma_{d_\chi})^\top$ contains the coordinate-dependent shape parameters, and $\chi_r(p,x)$ and $c_r$ denote the $r$th components of $\chi(p,x)$ and $c$, respectively. When $\gamma_1=\cdots=\gamma_{d_\chi}=\gamma$, \eqref{eq:anisotropic_rbf} reduces to the isotropic Gaussian basis in \eqref{eq:gaussian_rbf}.

\section{Backward Quadrature Learning Method}
\label{sec:proposed_method}

\subsection{Overview of the BQL Framework}
\label{subsec:global_architecture}

The time-discrete BSDE formulation in Section~\ref{Preliminaries} characterizes the $Z$ and $Y$ components through a sequence of conditional expectations. Building on this structure, we propose a BQL framework for solving high dimensional parametric nonlinear parabolic PDEs. Instead of introducing separate neural network approximations for the conditional expectation mappings, BQL directly learns the quadrature nodes used in their finite-sum approximations. These nodes are optimized over parameter--state samples and subsequently reused for different queries within the prescribed parameter domain. A continuation surrogate is introduced to transfer the learned value information between adjacent time levels while avoiding recursively nested quadrature evaluations.

At each time level $t_i$, $i=N-1,\ldots,0$, we introduce two sets of trainable quadrature parameters, $\vartheta_i^Z=(\vartheta_{i,m}^Z)_{m=1}^{M}$ and $\vartheta_i^Y=(\vartheta_{i,m}^Y)_{m=1}^{M}$, where $\vartheta_{i,m}^Z,\vartheta_{i,m}^Y\in\mathbb R^d$ and $M$ denotes the number of quadrature nodes in each set. The same quadrature learning mechanism is applied to the different integrands arising in the $Z$- and $Y$-recursions, with the corresponding parameters optimized separately. During training, MC estimates serve as reference targets for learning the quadrature parameters. Once optimized, the learned parameters $\widehat\vartheta_i^Z$ and $\widehat\vartheta_i^Y$ are fixed and reused without repeated MC sampling during evaluation.

The BQL framework proceeds sequentially backward in time. Given the continuation surrogate $\widehat y_{i+1}^{\mathrm{con}}$ from the subsequent time level, the $Z$-quadrature parameters are optimized first. The resulting $Z$-approximation is then used to optimize the $Y$-quadrature parameters, after which a new continuation surrogate $\widehat y_i^{\mathrm{con}}:\mathcal P\times\mathbb R^d\to\mathbb R$ is fitted to the learned $Y$-approximation. The backward dependency is therefore
\[
\widehat y_{i+1}^{\mathrm{con}}
\longrightarrow
\widehat\vartheta_i^Z
\longrightarrow
\widehat\vartheta_i^Y
\longrightarrow
\widehat y_i^{\mathrm{con}},
\]
where $i=N-1,\ldots,0$. The recursion is initialized by $\widehat y_N^{\mathrm{con}}(p,x)=g_p(x)$. 



\subsection{Quadrature Learning of BSDE Components}
\label{subsec:backward_quadrature}

Throughout the following $Z$- and $Y$-component updates, successor values at time $t_{i+1}$ are evaluated using the terminal function $g_p$ when $i=N-1$ and the continuation surrogate $\widehat y_{i+1}^{\mathrm{con}}$ otherwise. This convention applies to both MC reference evaluations and learned quadrature approximations.

\subsubsection{Training Samples and Learnable Quadrature Parameters}
\label{subsubsec:forward_paths}

We first generate $H$ forward training trajectories, referred to as base paths, where $H$ denotes the training batch size. The parameters $\{p^{(h)}\}_{h=1}^{H}$ are sampled independently and uniformly from $\mathcal P$, and the corresponding states are initialized by $X_0^{(h),\pi}:=X_0^{p^{(h)},\pi}=x_0$, where $x_0\in\mathbb R^d$ is the prescribed initial state. For each $h=1,\ldots,H$ and $i=0,\ldots,N-1$, an independent Brownian increment $\Delta W_i^{(h)}\sim\mathcal N(0,\Delta t_i I_d)$ is generated, where $I_d$ denotes the $d\times d$ identity matrix. The states are propagated according to the Euler--Maruyama scheme
\[
X_{i+1}^{(h),\pi}=X_i^{(h),\pi}+\mu_{p^{(h)}}\!\left(t_i,X_i^{(h),\pi}\right)\Delta t_i+\sigma_{p^{(h)}}\!\left(t_i,X_i^{(h),\pi}\right)\Delta W_i^{(h)}.
\]
The resulting parameter--state samples $\{(p^{(h)},X_i^{(h),\pi})\}_{h=1}^{H}$ are reused throughout the local training procedure at time level $t_i$.

For each component $\alpha\in\{Z,Y\}$, the trainable quadrature parameters $\vartheta_i^\alpha=(\vartheta_{i,m}^\alpha)_{m=1}^{M}$ introduced in Section~\ref{subsec:global_architecture} are represented by $M$ Brownian-increment-scale vectors in $\mathbb R^d$. They are initialized independently according to
\[
\vartheta_{i,m}^{\alpha,(0)}
\sim\mathcal N(0,\Delta t_i I_d),
\qquad m=1,\ldots,M,\quad \alpha\in\{Z,Y\},
\]
where the superscript $(0)$ denotes initialization. Unlike the independently sampled Brownian increments used to generate the base paths, the quadrature parameters are trainable and shared across all parameter--state samples at the same time level. Their Gaussian initialization is not imposed as a distributional constraint during optimization.

At each time level, the quadrature parameters are optimized by matching the corresponding finite-node approximations to MC reference estimates. Specifically, for $\alpha\in\{Z,Y\}$, the local training objective is
\[
\mathcal L_{i,\alpha}(\vartheta_i^\alpha)
=
\frac{1}{H}\sum_{h=1}^{H}
\left\|
\mathcal E_{i,M}^{\alpha,(h)}(\vartheta_i^\alpha)
-
\mathcal E_{i,K_\alpha}^{\alpha,(h)}
\right\|_2^2,
\]
where $\mathcal E_{i,K_\alpha}^{\alpha,(h)}$ denotes the MC reference estimate based on $K_\alpha$ samples, and $\mathcal E_{i,M}^{\alpha,(h)}(\vartheta_i^\alpha)$ denotes the corresponding learned quadrature approximation. These component-specific estimators are defined in the following subsections. For each component, the MC reference estimates are generated once at the current time level and held fixed during the subsequent optimization iterations.

The quadrature parameters are updated using the Adam optimizer \cite{Ruder2017,KingmaBa2015} with gradient-norm clipping. Denoting the optimization iteration by $r$, the update is written as
\[
\vartheta_i^{\alpha,(r+1)}
=
\operatorname{AdamStep}_{\eta_\alpha}
\left(
\vartheta_i^{\alpha,(r)},
\operatorname{clip}_{c_{\mathrm{clip}}}
\left(
\nabla_{\vartheta_i^\alpha}
\mathcal L_{i,\alpha}
\left(\vartheta_i^{\alpha,(r)}\right)
\right)
\right),
\]
where $\eta_\alpha$ denotes the learning rate, $\operatorname{clip}_{c_{\mathrm{clip}}}$ denotes gradient-norm clipping with threshold $c_{\mathrm{clip}}$, and $\operatorname{AdamStep}$ represents a standard Adam update, including its internal moment estimates. The $Z$-component parameters are optimized first. The resulting approximation is then fixed when optimizing the $Y$-component parameters. After training, the learned quadrature parameters $\widehat\vartheta_i^Z$ and $\widehat\vartheta_i^Y$ are retained for the subsequent backward evaluations.

\subsubsection{Learning the Z-Component}
\label{subsubsec:Z_component}

At time level $t_i$, the continuation surrogate $\widehat y_{i+1}^{\mathrm{con}}$ is assumed to be available from the subsequent time level. For each base sample $(p^{(h)},X_i^{(h),\pi})$, we generate $K_Z$ independent standard Gaussian samples $\zeta_{i,k}^{Z,(h)}\sim\mathcal N(0,I_d)$, $k=1,\ldots,K_Z$, and construct the corresponding successor states
\begin{equation}
X_{i+1,k}^{Z,(h),\pi}
=
X_i^{(h),\pi}
+
\mu_{p^{(h)}}\!\left(t_i,X_i^{(h),\pi}\right)\Delta t_i
+
\sigma_{p^{(h)}}\!\left(t_i,X_i^{(h),\pi}\right)
\sqrt{\Delta t_i}\,\zeta_{i,k}^{Z,(h)}.
\label{eq:Z_MC_successor}
\end{equation}
The corresponding successor values are given by $Y_{i+1,k}^{Z,(h)}=\widehat y_{i+1}^{\mathrm{con}}(p^{(h)},X_{i+1,k}^{Z,(h),\pi})$. The MC reference estimate for the $Z$-component is
\[
\mathcal E_{i,K_Z}^{Z,(h)}
=
\frac{1}{K_Z}
\sum_{k=1}^{K_Z}
Y_{i+1,k}^{Z,(h)}
\frac{\zeta_{i,k}^{Z,(h)}}{\sqrt{\Delta t_i}},
\qquad h=1,\ldots,H.
\]

To approximate this conditional expectation, we employ the trainable quadrature parameters $\vartheta_i^Z=(\vartheta_{i,m}^Z)_{m=1}^{M}$ introduced in Section~\ref{subsubsec:forward_paths}. Unlike the independently sampled MC increments, these parameters are shared across all $H$ base samples and updated through the quadrature learning procedure. For each $h=1,\ldots,H$ and $m=1,\ldots,M$, the successor state associated with $\vartheta_{i,m}^Z$ is
\[
X_{i+1,m}^{Z,(h),\pi}
=
X_i^{(h),\pi}
+
\mu_{p^{(h)}}\!\left(t_i,X_i^{(h),\pi}\right)\Delta t_i
+
\sigma_{p^{(h)}}\!\left(t_i,X_i^{(h),\pi}\right)
\vartheta_{i,m}^{Z}.
\]
Let $Y_{i+1,m}^{Z,(h)}=\widehat y_{i+1}^{\mathrm{con}}(p^{(h)},X_{i+1,m}^{Z,(h),\pi})$ denote the corresponding successor value. The resulting learned quadrature approximation is
\begin{equation}
\mathcal E_{i,M}^{Z,(h)}\!\left(\vartheta_i^Z\right)
=
\frac{1}{M}
\sum_{m=1}^{M}
Y_{i+1,m}^{Z,(h)}
\frac{\vartheta_{i,m}^{Z}}{\Delta t_i}.
\label{eq:Z_BQL_approximation}
\end{equation}

The quadrature parameters are trained by minimizing the discrepancy between the learned approximation and the MC reference estimates over the $H$ base samples. The corresponding optimization problem is
\begin{equation}
\min_{\vartheta_i^Z}
\mathcal L_{i,Z}\!\left(\vartheta_i^Z\right)
=
\min_{\vartheta_i^Z}
\frac{1}{H}
\sum_{h=1}^{H}
\left\|
\mathcal E_{i,K_Z}^{Z,(h)}
-
\mathcal E_{i,M}^{Z,(h)}\!\left(\vartheta_i^Z\right)
\right\|_2^2.
\label{eq:Z_node_optimization}
\end{equation}
Following the optimization procedure described in Section~\ref{subsubsec:forward_paths}, the parameters obtained after a finite number of iterations are denoted by $\widehat\vartheta_i^Z=(\widehat\vartheta_{i,m}^Z)_{m=1}^{M}$. These learned parameters are then fixed and define the $Z$-component approximation
\begin{equation}
\widehat z_i^{\mathrm{BQL}}(p,x)
=
\frac{1}{M}
\sum_{m=1}^{M}
\widehat y_{i+1}^{\mathrm{con}}
\left(
p,
x+\mu_p(t_i,x)\Delta t_i
+\sigma_p(t_i,x)\widehat\vartheta_{i,m}^{Z}
\right)
\frac{\widehat\vartheta_{i,m}^{Z}}{\Delta t_i}.
\label{eq:Z_BQL_function}
\end{equation}
In particular, for each base sample $(p^{(h)},X_i^{(h),\pi})$, we have $\widehat z_i^{\mathrm{BQL}}(p^{(h)},X_i^{(h),\pi})=\mathcal E_{i,M}^{Z,(h)}(\widehat\vartheta_i^Z)$.

The learned approximation targets the conditional expectation defining $Z_i^{p,\pi}$ in \eqref{eq:YZ_conditional}. Since this expectation is generally unavailable in closed form, the MC estimate provides a reference target during training, while the optimized quadrature parameters yield a reusable finite-sum approximation across different parameter--state inputs.

\subsubsection{Learning the Y-Component}
\label{subsubsec:Y_component}

Once the $Z$-component approximation $\widehat z_i^{\mathrm{BQL}}$ has been obtained, we proceed to the $Y$-component using the same quadrature learning mechanism, but with the integrand determined by the corresponding BSDE recursion. The optimized $Z$-component is kept fixed throughout this stage.

For each base sample $(p^{(h)},X_i^{(h),\pi})$, we independently generate $K_Y$ standard Gaussian samples $\zeta_{i,k}^{Y,(h)}\sim\mathcal N(0,I_d)$, $k=1,\ldots,K_Y$, and construct the successor states
\[
X_{i+1,k}^{Y,(h),\pi}
=
X_i^{(h),\pi}
+
\mu_{p^{(h)}}\!\left(t_i,X_i^{(h),\pi}\right)\Delta t_i
+
\sigma_{p^{(h)}}\!\left(t_i,X_i^{(h),\pi}\right)
\sqrt{\Delta t_i}\,\zeta_{i,k}^{Y,(h)}.
\]
The corresponding successor values are $Y_{i+1,k}^{Y,(h)}=\widehat y_{i+1}^{\mathrm{con}}(p^{(h)},X_{i+1,k}^{Y,(h),\pi})$. With the trained $Z$-component evaluated at the current base sample, the MC reference estimate for the $Y$-component is
\begin{equation}
\begin{aligned}
\mathcal E_{i,K_Y}^{Y,(h)}
=
\frac{1}{K_Y}
\sum_{k=1}^{K_Y}
\Big[
Y_{i+1,k}^{Y,(h)}
+\Delta t_i
f_{p^{(h)}}\!\left(
t_i,X_i^{(h),\pi},
Y_{i+1,k}^{Y,(h)},
\widehat z_i^{\mathrm{BQL}}\!\left(p^{(h)},X_i^{(h),\pi}\right)
\right)
\Big].
\end{aligned}
\label{eq:Y_mc_target}
\end{equation}

We then employ the trainable quadrature parameters $\vartheta_i^Y=(\vartheta_{i,m}^Y)_{m=1}^{M}$, which are optimized separately from $\vartheta_i^Z$. For each $h=1,\ldots,H$ and $m=1,\ldots,M$, the corresponding successor state is
\[
X_{i+1,m}^{Y,(h),\pi}
=
X_i^{(h),\pi}
+
\mu_{p^{(h)}}\!\left(t_i,X_i^{(h),\pi}\right)\Delta t_i
+
\sigma_{p^{(h)}}\!\left(t_i,X_i^{(h),\pi}\right)
\vartheta_{i,m}^{Y}.
\]
Let $Y_{i+1,m}^{Y,(h)}=\widehat y_{i+1}^{\mathrm{con}}(p^{(h)},X_{i+1,m}^{Y,(h),\pi})$ denote the associated successor value. The learned quadrature approximation of the $Y$-component is
\begin{equation}
\begin{aligned}
\mathcal E_{i,M}^{Y,(h)}\!\left(\vartheta_i^Y\right)
=
\frac{1}{M}
\sum_{m=1}^{M}
\Big[
Y_{i+1,m}^{Y,(h)}
+\Delta t_i
f_{p^{(h)}}\!\left(
t_i,X_i^{(h),\pi},
Y_{i+1,m}^{Y,(h)},
\widehat z_i^{\mathrm{BQL}}\!\left(p^{(h)},X_i^{(h),\pi}\right)
\right)
\Big].
\end{aligned}
\label{eq:Y_BQL_approximation}
\end{equation}

The $Y$-component quadrature parameters are trained by solving
\begin{equation}
\min_{\vartheta_i^Y}
\mathcal L_{i,Y}\!\left(\vartheta_i^Y\right)
=
\min_{\vartheta_i^Y}
\frac{1}{H}
\sum_{h=1}^{H}
\left|
\mathcal E_{i,K_Y}^{Y,(h)}
-
\mathcal E_{i,M}^{Y,(h)}\!\left(\vartheta_i^Y\right)
\right|^2.
\label{eq:Y_node_optimization}
\end{equation}
The MC reference estimates and the trained $Z$-component remain fixed during this optimization. Denoting the quadrature parameters obtained after the prescribed optimization iterations by $\widehat\vartheta_i^Y=(\widehat\vartheta_{i,m}^Y)_{m=1}^{M}$, we define the resulting $Y$-component approximation as
\begin{equation}
\begin{aligned}
\widehat y_i^{\mathrm{BQL}}(p,x)
=
\frac{1}{M}
\sum_{m=1}^{M}
\Bigg[
&\widehat y_{i+1}^{\mathrm{con}}
\left(
p,
x+\mu_p(t_i,x)\Delta t_i
+\sigma_p(t_i,x)\widehat\vartheta_{i,m}^{Y}
\right)
\\
&+\Delta t_i f_p
\Bigg(
t_i,x,
\widehat y_{i+1}^{\mathrm{con}}
\left(
p,
x+\mu_p(t_i,x)\Delta t_i
+\sigma_p(t_i,x)\widehat\vartheta_{i,m}^{Y}
\right),
\widehat z_i^{\mathrm{BQL}}(p,x)
\Bigg)
\Bigg].
\end{aligned}
\label{eq:Y_BQL_function}
\end{equation}
For each base sample $(p^{(h)},X_i^{(h),\pi})$, the approximation satisfies $\widehat y_i^{\mathrm{BQL}}(p^{(h)},X_i^{(h),\pi})=\mathcal E_{i,M}^{Y,(h)}(\widehat\vartheta_i^Y)$.

The two BSDE components therefore share the same quadrature learning mechanism, while their trainable parameters are fitted to different integrands. Once the $Y$-component has been constructed, its learned approximation provides the training targets for the continuation surrogate at time level $t_i$, as described in Section~\ref{subsec:continuation_surrogate}.

\subsection{Continuation Surrogate Approximation}
\label{subsec:continuation_surrogate}

The preceding quadrature learning procedure yields the approximations $\widehat z_i^{\mathrm{BQL}}$ and $\widehat y_i^{\mathrm{BQL}}$ at time level $t_i$. Although these approximations can be evaluated at arbitrary parameter--state inputs, the backward recursion at $t_{i-1}$ requires repeated evaluations of $\widehat y_i^{\mathrm{BQL}}$ at newly generated successor states. Directly substituting its finite-sum representation into earlier time-level updates would introduce recursively nested quadrature evaluations, whose computational cost may grow geometrically with the number of backward steps. To avoid this repeated expansion, we introduce a continuation surrogate $\widehat y_i^{\mathrm{con}}$ to approximate $\widehat y_i^{\mathrm{BQL}}$ and provide an explicitly evaluable representation for the subsequent backward recursion.

To construct the continuation training data, we generate $N_{\mathrm{con}}$ forward trajectories independently of the $H$ base paths used for quadrature learning. Specifically, the parameters $\{p^{(n)}\}_{n=1}^{N_{\mathrm{con}}}$ are sampled independently and uniformly from $\mathcal P$, with $X_{0,\mathrm{con}}^{(n),\pi}=x_0$. The corresponding states are propagated according to
\[
X_{j+1,\mathrm{con}}^{(n),\pi}
=
X_{j,\mathrm{con}}^{(n),\pi}
+
\mu_{p^{(n)}}\!\left(t_j,X_{j,\mathrm{con}}^{(n),\pi}\right)\Delta t_j
+
\sigma_{p^{(n)}}\!\left(t_j,X_{j,\mathrm{con}}^{(n),\pi}\right)
\Delta W_{j,\mathrm{con}}^{(n)},
\qquad j=0,\ldots,N-1,
\]
where $\Delta W_{j,\mathrm{con}}^{(n)}\sim\mathcal N(0,\Delta t_j I_d)$ are independent Brownian increments. These trajectories are generated once, and their states at the corresponding time levels are reused during the backward training procedure.

After the quadrature parameters $\widehat\vartheta_i^Z$ and $\widehat\vartheta_i^Y$ have been optimized, they are fixed to construct the continuation targets. For each sample $(p^{(n)},X_{i,\mathrm{con}}^{(n),\pi})$, we first evaluate $\widehat z_i^{\mathrm{BQL}}$ using \eqref{eq:Z_BQL_function} and then evaluate $\widehat y_i^{\mathrm{BQL}}$ using \eqref{eq:Y_BQL_function}. Thus, the continuation targets are obtained entirely from the learned quadrature approximations, without introducing additional MC reference estimates. The resulting dataset is
\[
\mathcal D_i^{\mathrm{con}}
=
\left\{
\left(
p^{(n)},
X_{i,\mathrm{con}}^{(n),\pi},
\widehat y_i^{\mathrm{BQL}}
\left(p^{(n)},X_{i,\mathrm{con}}^{(n),\pi}\right)
\right)
\right\}_{n=1}^{N_{\mathrm{con}}}.
\]
The target values are computed before continuation training and remain fixed throughout the optimization iterations.

To approximate the resulting parameter--state mapping, we employ a hybrid architecture combining anisotropic Gaussian RBF with a FNN. Let $\chi(p,x)$ denote the parameter--state feature mapping introduced in Section~\ref{Preliminaries}. At each time level, the features are standardized using the empirical means and standard deviations of the continuation samples, yielding a normalized representation $\widetilde\chi_i(p,x)$. Both components of the hybrid surrogate use this normalized representation. The RBF component is defined by
\[
\mathcal R_i^{\mathrm{AG}}(p,x)
=
\sum_{j=1}^{N_c}
w_{i,j}\,
\phi^{\mathrm{AG}}\!\left(
\widetilde\chi_i(p,x),
c_{i,j};
\boldsymbol{\gamma}_{i,j}
\right)
+b_i,
\]
where $N_c$ is the number of RBF centers, while $c_{i,j}$, $w_{i,j}$, and $\boldsymbol{\gamma}_{i,j}$ denote the center, coefficient, and positive coordinate-dependent shape parameters of the $j$th basis function at time level $t_i$, respectively. The centers are initialized by selecting samples from the normalized feature set and are subsequently optimized together with the other RBF parameters.

The complete continuation surrogate is represented by
\[
\widehat y_i^{\mathrm{con}}(p,x)
=
\mathcal R_i^{\mathrm{AG}}(p,x)
+
\mathcal M_i^{\mathrm{con}}(p,x;\omega_i^{\mathrm{con}}),
\]
where $\mathcal M_i^{\mathrm{con}}$ denotes the FNN component acting on $\widetilde\chi_i(p,x)$, and $\omega_i^{\mathrm{con}}$ represents its trainable parameters. The two components are fitted jointly, with the FNN providing an additive contribution to the RBF approximation.

With the learned quadrature parameters fixed, the continuation surrogate is trained by minimizing
\[
\mathcal L_{i,\mathrm{con}}
=
\frac{1}{N_{\mathrm{con}}}
\sum_{n=1}^{N_{\mathrm{con}}}
\left|
\widehat y_i^{\mathrm{BQL}}
\left(p^{(n)},X_{i,\mathrm{con}}^{(n),\pi}\right)
-
\widehat y_i^{\mathrm{con}}
\left(p^{(n)},X_{i,\mathrm{con}}^{(n),\pi}\right)
\right|^2.
\]
The trainable RBF parameters $\{c_{i,j},w_{i,j},\boldsymbol{\gamma}_{i,j}\}_{j=1}^{N_c}$, the bias $b_i$, and the FNN parameters $\omega_i^{\mathrm{con}}$ are jointly updated using Adam with gradient-norm clipping. The learned quadrature approximation $\widehat y_i^{\mathrm{BQL}}$ serves as the fixed training target, and its parameters are not updated during this stage. After training, $\widehat y_i^{\mathrm{con}}$ is fixed and used in place of the explicitly expanded quadrature approximation when evaluating successor values at the preceding time level $t_{i-1}$. This completes the local backward update and enables the sequential recursion to continue.

\subsection{Complete BQL Algorithm}
\label{subsec:complete_algorithm}

The preceding construction yields a sequential backward training procedure. At each time level, the $Z$-component quadrature parameters are optimized first, followed by the $Y$-component parameters. The resulting learned $Y$-approximation is then used to train the continuation surrogate for the next backward step. The MC reference estimates remain fixed during each component's optimization, while the continuation targets are constructed from the trained quadrature approximations. The complete offline procedure is summarized in Algorithm~\ref{alg:solver}.

\begin{algorithm}[htbp]
\small
\caption{BQL Algorithm for high dimensional parametric nonlinear parabolic PDEs}
\label{alg:solver}
\SetAlgoLined
\LinesNotNumbered
\SetAlgoLongEnd

\KwIn{PDE coefficients $\mu_p,\sigma_p,f_p,g_p$; time partition $\pi=\{t_i\}_{i=0}^{N}$; sample sizes $H,N_{\mathrm{con}},K_Z,K_Y$; quadrature node number $M$; RBF center number $N_c$; training iterations $R_Z,R_Y,R_{\mathrm{con}}$; learning rates $\eta_Z,\eta_Y,\eta_{\mathrm{con}}$; clipping threshold $c_{\mathrm{clip}}$.}

\KwOut{Trained quadrature parameters $\{\widehat\vartheta_i^Z,\widehat\vartheta_i^Y\}_{i=0}^{N-1}$ and continuation surrogates $\{\widehat y_i^{\mathrm{con}}\}_{i=0}^{N}$.}

\underline{\textbf{Forward generation of training trajectories:}}

\For{$h=1,\ldots,H$}{
Sample $p^{(h)}$ uniformly from $\mathcal P$ and set $X_0^{(h),\pi}=x_0$\;
\For{$i=0,\ldots,N-1$}{
Generate $\Delta W_i^{(h)}\sim\mathcal N(0,\Delta t_i I_d)$\;
Compute $X_{i+1}^{(h),\pi}\leftarrow X_i^{(h),\pi}+\mu_{p^{(h)}}(t_i,X_i^{(h),\pi})\Delta t_i+\sigma_{p^{(h)}}(t_i,X_i^{(h),\pi})\Delta W_i^{(h)}$\;
}}

\underline{\textbf{Sequential backward training:}}
Set $\widehat y_N^{\mathrm{con}}(p,x)=g_p(x)$\;

\For{$i=N-1,\ldots,0$}{

\underline{\textbf{$Z$-component quadrature learning:}}

Initialize $\vartheta_{i,m}^{Z,(0)}\sim\mathcal N(0,\Delta t_i I_d)$ independently for $m=1,\ldots,M$\;

Apply Algorithm~\ref{alg:Z_component} to construct the fixed references $\{\mathcal E_{i,K_Z}^{Z,(h)}\}_{h=1}^{H}$ and the parameterized evaluations $\mathcal E_{i,M}^{Z,(h)}(\vartheta_i^Z)$\;

Initialize an Adam optimizer with learning rate $\eta_Z$\;

\For{$r=0,\ldots,R_Z-1$}{
Evaluate $\mathcal L_{i,Z}(\vartheta_i^{Z,(r)})=\frac{1}{H}\sum_{h=1}^{H}\left\|\mathcal E_{i,K_Z}^{Z,(h)}-\mathcal E_{i,M}^{Z,(h)}(\vartheta_i^{Z,(r)})\right\|_2^2$\;

$G_i^{Z,(r)}\leftarrow\operatorname{clip}_{c_{\mathrm{clip}}}\left(\nabla_{\vartheta_i^Z}\mathcal L_{i,Z}(\vartheta_i^{Z,(r)})\right)$\;

$\vartheta_i^{Z,(r+1)}\leftarrow\operatorname{AdamStep}_{\eta_Z}\left(\vartheta_i^{Z,(r)},G_i^{Z,(r)}\right)$\;
}
\underline{\textbf{$Y$-component quadrature learning:}}

Initialize $\vartheta_{i,m}^{Y,(0)}\sim\mathcal N(0,\Delta t_i I_d)$ independently for $m=1,\ldots,M$\;

Apply Algorithm~\ref{alg:Y_component} using $\widehat z_i^{\mathrm{BQL}}$ to construct the fixed references $\{\mathcal E_{i,K_Y}^{Y,(h)}\}_{h=1}^{H}$ and the parameterized evaluations $\mathcal E_{i,M}^{Y,(h)}(\vartheta_i^Y)$\;

Initialize an Adam optimizer with learning rate $\eta_Y$\;

\For{$r=0,\ldots,R_Y-1$}{
Evaluate $\mathcal L_{i,Y}(\vartheta_i^{Y,(r)})=\frac{1}{H}\sum_{h=1}^{H}\left|\mathcal E_{i,K_Y}^{Y,(h)}-\mathcal E_{i,M}^{Y,(h)}(\vartheta_i^{Y,(r)})\right|^2$\;

$G_i^{Y,(r)}\leftarrow\operatorname{clip}_{c_{\mathrm{clip}}}\left(\nabla_{\vartheta_i^Y}\mathcal L_{i,Y}(\vartheta_i^{Y,(r)})\right)$\;

$\vartheta_i^{Y,(r+1)}\leftarrow\operatorname{AdamStep}_{\eta_Y}\left(\vartheta_i^{Y,(r)},G_i^{Y,(r)}\right)$\;
}
\underline{\textbf{Continuation surrogate training:}}

Apply Algorithm~\ref{alg:surrogate_component} to fit $\widehat y_i^{\mathrm{con}}$ using the stored continuation samples and $\widehat y_i^{\mathrm{BQL}}$\;

Store $\widehat\vartheta_i^Z$, $\widehat\vartheta_i^Y$, and $\widehat y_i^{\mathrm{con}}$ for subsequent evaluations\;
}

\end{algorithm}

\begin{algorithm}[htbp]
\caption{MC Reference and Learned Quadrature Evaluations for the $Z$-Component}
\label{alg:Z_component}
\SetAlgoLined
\LinesNotNumbered
\SetAlgoLongEnd

\KwIn{Time level $t_i$; base samples $\{(p^{(h)},X_i^{(h),\pi})\}_{h=1}^{H}$; successor continuation function $\widehat y_{i+1}^{\mathrm{con}}$; trainable quadrature parameters $\vartheta_i^Z$; sample sizes $K_Z,M$.}

\KwOut{Fixed MC references $\{\mathcal E_{i,K_Z}^{Z,(h)}\}_{h=1}^{H}$ and parameterized quadrature evaluations $\{\mathcal E_{i,M}^{Z,(h)}(\vartheta_i^Z)\}_{h=1}^{H}$.}

\For{$h=1,\ldots,H$}{

\underline{\textbf{MC reference:}}

\For{$k=1,\ldots,K_Z$}{
Generate $\zeta_{i,k}^{Z,(h)}\sim\mathcal N(0,I_d)$\;

$X_{i+1,k}^{Z,(h),\pi}\leftarrow X_i^{(h),\pi}+\mu_{p^{(h)}}(t_i,X_i^{(h),\pi})\Delta t_i+\sigma_{p^{(h)}}(t_i,X_i^{(h),\pi})\sqrt{\Delta t_i}\,\zeta_{i,k}^{Z,(h)}$\;

$Y_{i+1,k}^{Z,(h)}\leftarrow\widehat y_{i+1}^{\mathrm{con}}(p^{(h)},X_{i+1,k}^{Z,(h),\pi})$\;
}

$\mathcal E_{i,K_Z}^{Z,(h)}\leftarrow\frac{1}{K_Z}\sum_{k=1}^{K_Z}Y_{i+1,k}^{Z,(h)}\frac{\zeta_{i,k}^{Z,(h)}}{\sqrt{\Delta t_i}}$\;

Fix $\mathcal E_{i,K_Z}^{Z,(h)}$ as the training reference\;

\underline{\textbf{Learnable quadrature evaluator:}}

\For{$m=1,\ldots,M$}{
Define $X_{i+1,m}^{Z,(h),\pi}(\vartheta_i^Z)\leftarrow X_i^{(h),\pi}+\mu_{p^{(h)}}(t_i,X_i^{(h),\pi})\Delta t_i+\sigma_{p^{(h)}}(t_i,X_i^{(h),\pi})\vartheta_{i,m}^{Z}$\;

Define $Y_{i+1,m}^{Z,(h)}(\vartheta_i^Z)\leftarrow\widehat y_{i+1}^{\mathrm{con}}(p^{(h)},X_{i+1,m}^{Z,(h),\pi}(\vartheta_i^Z))$\;
}

Define $\mathcal E_{i,M}^{Z,(h)}(\vartheta_i^Z)\leftarrow\frac{1}{M}\sum_{m=1}^{M}Y_{i+1,m}^{Z,(h)}(\vartheta_i^Z)\frac{\vartheta_{i,m}^{Z}}{\Delta t_i}$\;
}

\end{algorithm}

\begin{algorithm}[htbp]
\caption{MC Reference and Learned Quadrature Evaluations for the $Y$-Component}
\label{alg:Y_component}
\SetAlgoLined
\LinesNotNumbered
\SetAlgoLongEnd

\KwIn{Time level $t_i$; base samples $\{(p^{(h)},X_i^{(h),\pi})\}_{h=1}^{H}$; successor continuation function $\widehat y_{i+1}^{\mathrm{con}}$; trained $Z$-approximation $\widehat z_i^{\mathrm{BQL}}$; trainable quadrature parameters $\vartheta_i^Y$; sample sizes $K_Y,M$.}

\KwOut{Fixed MC references $\{\mathcal E_{i,K_Y}^{Y,(h)}\}_{h=1}^{H}$ and parameterized quadrature evaluations $\{\mathcal E_{i,M}^{Y,(h)}(\vartheta_i^Y)\}_{h=1}^{H}$.}

\For{$h=1,\ldots,H$}{

Set $z_i^{(h)}\leftarrow\widehat z_i^{\mathrm{BQL}}(p^{(h)},X_i^{(h),\pi})$ and keep it fixed during the $Y$-component training\;

\underline{\textbf{MC reference:}}

\For{$k=1,\ldots,K_Y$}{
Generate $\zeta_{i,k}^{Y,(h)}\sim\mathcal N(0,I_d)$\;

$X_{i+1,k}^{Y,(h),\pi}\leftarrow X_i^{(h),\pi}+\mu_{p^{(h)}}(t_i,X_i^{(h),\pi})\Delta t_i+\sigma_{p^{(h)}}(t_i,X_i^{(h),\pi})\sqrt{\Delta t_i}\,\zeta_{i,k}^{Y,(h)}$\;

$Y_{i+1,k}^{Y,(h)}\leftarrow\widehat y_{i+1}^{\mathrm{con}}(p^{(h)},X_{i+1,k}^{Y,(h),\pi})$\;
}

$\mathcal E_{i,K_Y}^{Y,(h)}\leftarrow\frac{1}{K_Y}\sum_{k=1}^{K_Y}\left[Y_{i+1,k}^{Y,(h)}+\Delta t_i f_{p^{(h)}}(t_i,X_i^{(h),\pi},Y_{i+1,k}^{Y,(h)},z_i^{(h)})\right]$\;

Fix $\mathcal E_{i,K_Y}^{Y,(h)}$ as the training reference\;

\underline{\textbf{Learnable quadrature evaluator:}}

\For{$m=1,\ldots,M$}{
Define $X_{i+1,m}^{Y,(h),\pi}(\vartheta_i^Y)\leftarrow X_i^{(h),\pi}+\mu_{p^{(h)}}(t_i,X_i^{(h),\pi})\Delta t_i+\sigma_{p^{(h)}}(t_i,X_i^{(h),\pi})\vartheta_{i,m}^{Y}$\;

Define $Y_{i+1,m}^{Y,(h)}(\vartheta_i^Y)\leftarrow\widehat y_{i+1}^{\mathrm{con}}(p^{(h)},X_{i+1,m}^{Y,(h),\pi}(\vartheta_i^Y))$\;
}

Define $\mathcal E_{i,M}^{Y,(h)}(\vartheta_i^Y)\leftarrow\frac{1}{M}\sum_{m=1}^{M}\left[Y_{i+1,m}^{Y,(h)}(\vartheta_i^Y)+\Delta t_i f_{p^{(h)}}(t_i,X_i^{(h),\pi},Y_{i+1,m}^{Y,(h)}(\vartheta_i^Y),z_i^{(h)})\right]$\;
}

\end{algorithm}

\begin{algorithm}[htbp]
\caption{Construction and Training of the Continuation Surrogate}
\label{alg:surrogate_component}
\SetAlgoLined
\LinesNotNumbered
\SetAlgoLongEnd

\KwIn{Time level $t_i$; continuation samples $\{(p^{(n)},X_{i,\mathrm{con}}^{(n),\pi})\}_{n=1}^{N_{\mathrm{con}}}$; trained quadrature parameters $\widehat\vartheta_i^Z,\widehat\vartheta_i^Y$; approximations $\widehat z_i^{\mathrm{BQL}},\widehat y_i^{\mathrm{BQL}}$; feature map $\chi$; RBF center number $N_c$; training iterations $R_{\mathrm{con}}$; learning rate $\eta_{\mathrm{con}}$; clipping threshold $c_{\mathrm{clip}}$.}

\KwOut{Fitted continuation surrogate $\widehat y_i^{\mathrm{con}}(p,x)$.}

\underline{\textbf{Continuation target construction:}}

\For{$n=1,\ldots,N_{\mathrm{con}}$}{
Set $q_n\leftarrow(p^{(n)},X_{i,\mathrm{con}}^{(n),\pi})$\;

Evaluate $\widehat z_i^{\mathrm{BQL}}(q_n)$ and set $y_n^{\mathrm{tar}}\leftarrow\widehat y_i^{\mathrm{BQL}}(q_n)$ using the trained quadrature parameters\;
}

Construct $\mathcal D_i^{\mathrm{con}}=\{(q_n,y_n^{\mathrm{tar}})\}_{n=1}^{N_{\mathrm{con}}}$ and fix all targets\;

\underline{\textbf{Hybrid surrogate initialization:}}

Compute and normalize $\{\chi(q_n)\}_{n=1}^{N_{\mathrm{con}}}$ to obtain $\widetilde\chi_i$\;

Initialize $N_c$ RBF centers from distinct normalized samples\;

Initialize the RBF and MLP parameters, using softplus for positive effective RBF shape parameters\;

Construct
$\widehat y_i^{\mathrm{con}}(p,x)
=\sum_{j=1}^{N_c}w_{i,j}\phi^{\mathrm{AG}}\!\left(\widetilde\chi_i(p,x),c_{i,j};\boldsymbol{\gamma}_{i,j}\right)
+b_i+\mathcal M_i^{\mathrm{con}}(p,x;\omega_i^{\mathrm{con}})$\;

\underline{\textbf{Joint surrogate training:}}

Initialize Adam with learning rate $\eta_{\mathrm{con}}$\;

\For{$r=0,\ldots,R_{\mathrm{con}}-1$}{
Evaluate
$\mathcal L_{i,\mathrm{con}}
=\frac{1}{N_{\mathrm{con}}}
\sum_{n=1}^{N_{\mathrm{con}}}
\left|y_n^{\mathrm{tar}}-\widehat y_i^{\mathrm{con}}(q_n)\right|^2$\;

Compute gradients, clip each parameter gradient tensor to $c_{\mathrm{clip}}$, and jointly update the RBF and MLP parameters using Adam\;
}

Fix and retain $\widehat y_i^{\mathrm{con}}$ for the next backward time level\;

\end{algorithm}


\section{Theoretical Analysis}
\label{sec:theoretical_analysis}
This section analyzes how temporal discretization, learnable finite-node quadrature, and continuation-surrogate fitting contribute to the error of the proposed BQL method. We first establish a global structural error estimate to clarify the relationship between the total error and its local components. We then analyze the local BQL errors and continuation-approximation error separately, and combine the resulting estimates to derive a quantitative global error bound and the corresponding consistency result. Detailed stability arguments, supplementary local-risk analysis, and architecture-dependent continuation-approximation estimates are provided in Appendices~\ref{app:stability}--\ref{app:architecture}, respectively. The theoretical guarantees are established under the standing regularity assumptions specified below, and their applicability is therefore restricted to PDEs whose associated BSDE generators satisfy these conditions.

\subsection{Analytical Framework and Main Error Estimate}
\label{subsec:analytical_setting}
At each time level, the same BQL quadrature nodes are used across parameter--state inputs. We therefore assess their approximation errors with respect to the joint distribution of the parameter and the corresponding forward state. Let $P\sim\nu_{\mathcal P}$ be a random parameter taking values in $\mathcal P$, where $\nu_{\mathcal P}$ is a prescribed probability measure on $\mathcal P$. We assume that $P$ is $\mathcal F_0$-measurable and independent
of the Brownian motion $W$. The implementation in Section~\ref{sec:proposed_method} uses the uniform distribution on $\mathcal P$, while the analysis allows a general choice of $\nu_{\mathcal P}$. For $i=0,\ldots,N$, let $X_i^{P,\pi}$ denote the Euler approximation of the forward state at time $t_i$ associated with the sampled parameter $P$, and define $Q_i=(P,X_i^{P,\pi})$ and $\rho_i=\mathcal L(Q_i)$, where $\mathcal L(Q_i)$ denotes the probability distribution of $Q_i$. For any measurable function $\varphi$ that is square integrable with respect to $\rho_i$, define $\|\varphi\|_{L^2(\rho_i)}^2=\mathbb E[|\varphi(Q_i)|^2]$, the expectation is taken with respect to $\rho_i$.

We make the following standard assumptions for the error analysis in this section:
\begin{enumerate}
\item[(i)] $\mathbb E[|X_0^P|^2]<\infty$. 

\item[(ii)]  $(\mu_p,\sigma_p)$ are uniformly $1/2$-Hölder continuous with respect to the time variable $t$. There exist constants $C_X,L_X>0$, independent of $p$, for all $p\in\mathcal P$, $t\in[0,T]$, and $x,x'\in\mathbb R^d$, such that
\[
\sup_{\substack{p\in\mathcal P,\,t\in[0,T]}}
\left(|\mu_p(t,0)|+\|\sigma_p(t,0)\|\right)\le C_X,  \quad |\mu_p(t,x)-\mu_p(t,x')|+\|\sigma_p(t,x)-\sigma_p(t,x')\|
\le L_X|x-x'|.\]

\item[(iii)] 
$f_p$ is uniformly $1/2$-Hölder continuous with respect to the time variable $t$. There exist constants $C_f,L_f>0$, independent of $p$, for all $p\in\mathcal P$, $t\in[0,T]$, $x,x',z,z'\in\mathbb R^d$, and $y,y'\in\mathbb R$, such that
\[
\sup_{p\in\mathcal P,\,t\in[0,T]}|f_p(t,0,0,0)|\le C_f, \quad f_p(t,x,y,z)-f_p(t,x',y',z')|
\le L_f\bigl(|x-x'|+|y-y'|+|z-z'|\bigr).
\]

\item[(iv)] There exist constants $C_g,L_g>0$, independent of $p$, such that
\[
\sup_{p\in\mathcal P}|g_p(0)|\le C_g,\qquad |g_p(x)-g_p(x')|\le L_g|x-x'|,\quad p\in\mathcal P,\quad x,x'\in\mathbb R^d.
\]
\end{enumerate}

The standing assumptions above specify the class of parametric semilinear parabolic PDEs for which the theoretical error estimates in this section are established. In particular, the uniform Lipschitz continuity of the generator with respect to the state and BSDE variables is required for the stability analysis. Consequently, the theoretical convergence results apply to problems satisfying these regularity conditions and do not cover general nonlinear generators with superlinear growth.

Throughout the theoretical analysis, $C$ denotes a generic positive constant, whose value may change from line to line. Unless otherwise specified, it may depend on $T$, $d$, the uniform constants in the standing assumptions, and $\mathbb{E}[|X_0^P|^2]$, but is independent of the time partition, the time index, the parameter value, and the trained approximations.

We now estimate the error of the numerical solutions $(Y_i^{P},Z_i^{P})_{0\leq i \leq N}$ and the BQL approximations $(\widehat y_i^{\mathrm{con}}(Q_i),\widehat z_i^{\mathrm{BQL}}(Q_i))_{0\leq i \leq N}$.
To compare these quantities in a common parameter--state representation, we express the discrete BSDE solution in terms of deterministic functions. Specifically, for $i=0,\ldots,N-1$, we have $Y_i^{P,\pi}=y_i^\pi(Q_i)$ and $Z_i^{P,\pi}=z_i^\pi(Q_i)$, where $y_i^\pi:\mathcal P\times\mathbb R^d\to\mathbb R$ and $z_i^\pi:\mathcal P\times\mathbb R^d\to\mathbb R^d$ are deterministic functions. At the terminal time, $Y_N^{P,\pi}=y_N^\pi(Q_N)$, where $y_N^\pi(p,x)=g_p(x)$.
For notational convenience, we write parameter--state functions in terms of $q=(p,x)$. In particular, $\widehat z_i^{\mathrm{BQL}}$ and $\widehat y_i^{\mathrm{BQL}}$ retain their definitions in \eqref{eq:Z_BQL_function} and \eqref{eq:Y_BQL_function}, respectively. The corresponding one-step Euler transition is defined by
\[
\chi_i(q,\zeta)=x+\mu_p(t_i,x)\Delta t_i+\sigma_p(t_i,x)\sqrt{\Delta t_i}\,\zeta,
\]
where $\zeta\sim\mathcal N(0,I_d)$ is independent of the current parameter--state pair, and $\mathbb E_\zeta$ denotes expectation with respect to $\zeta$.

To separate the local quadrature errors from the error inherited through the continuation surrogate, we introduce two auxiliary functions. Given the fixed surrogate $\widehat y_{i+1}^{\mathrm{con}}$, define
\begin{equation}
\bar z_i(q)=\mathbb E_\zeta\left[\widehat y_{i+1}^{\mathrm{con}}\!\left(p,\mathcal X_i(q,\zeta)\right)\frac{\zeta}{\sqrt{\Delta t_i}}\right].
\label{eq:auxiliary_Z_function}
\end{equation}
With the trained approximation $\widehat z_i^{\mathrm{BQL}}$ fixed, define
\begin{equation}
\bar y_i(q)=\mathbb E_\zeta\left[\widehat y_{i+1}^{\mathrm{con}}\!\left(p,\mathcal X_i(q,\zeta)\right)+f_p\!\left(t_i,x,\widehat y_{i+1}^{\mathrm{con}}\!\left(p,\mathcal X_i(q,\zeta)\right),\widehat z_i^{\mathrm{BQL}}(q)\right)\Delta t_i\right].
\label{eq:auxiliary_Y_function}
\end{equation}
The auxiliary functions $\bar z_i$ and $\bar y_i$ are introduced solely for error analysis. The differences between $(\bar z_i,\bar y_i)$ and $(\widehat z_i^{\mathrm{BQL}},\widehat y_i^{\mathrm{BQL}})$ isolate the local quadrature errors, while the additional discrepancy between $\widehat y_i^{\mathrm{BQL}}$ and its continuation surrogate is treated separately.
For $i=0,\ldots,N-1$, define
\[
r_{i,Z}:=\|\bar z_i-\widehat z_i^{\mathrm{BQL}}\|_{L^2(\rho_i)}^2,\qquad r_{i,Y}:=\|\bar y_i-\widehat y_i^{\mathrm{BQL}}\|_{L^2(\rho_i)}^2,\qquad r_{i,\mathrm{con}}:=\|\widehat y_i^{\mathrm{BQL}}-\widehat y_i^{\mathrm{con}}\|_{L^2(\rho_i)}^2.
\]
The quantities $r_{i,Z}$ and $r_{i,Y}$ measure the local errors incurred by replacing the exact expectations in \eqref{eq:auxiliary_Z_function} and \eqref{eq:auxiliary_Y_function} with the learned finite-node quadrature rules, while $r_{i,\mathrm{con}}$ measures the error incurred by fitting the continuation surrogate to the learned $Y$-BQL function.

\begin{theorem}[Total error estimate]
\label{thm:total_error}
Under the standing conditions in Section~\ref{subsec:analytical_setting}, for sufficiently fine uniform partitions $\pi$ of $[0,T]$, we have
\begin{equation}
\begin{aligned}
\mathcal E_{\mathrm{total}}(\pi)&:=\max_{0\le i\le N}\mathbb E\!\left[|Y_{t_i}^P-\widehat y_i^{\mathrm{con}}(Q_i)|^2\right]+\mathbb E\!\left[\sum_{i=0}^{N-1}\int_{t_i}^{t_{i+1}}|Z_t^P-\widehat z_i^{\mathrm{BQL}}(Q_i)|^2\,\mathrm dt\right]\\
&\le C|\pi|+C\sum_{i=0}^{N-1}\left[\Delta t_i r_{i,Z}+\frac{r_{i,Y}+r_{i,\mathrm{con}}}{\Delta t_i}\right],
\label{eq:total_structural_error}
\end{aligned}
\end{equation}
where $r_{i,Z}$, $r_{i,Y}$, and $r_{i,\mathrm{con}}$ are the local errors defined above. The constant $C$ is independent of the time partition and the trained approximations.
\end{theorem}

\begin{proof}
Throughout the proof, $C$ denotes a generic positive constant, whose value may change from line to line and which satisfies the uniformity convention stated above.

We first decompose the total error by inserting the exact time-discrete BSDE solution $(y_i^\pi(Q_i),z_i^\pi(Q_i))$. Specifically, for $i=0,\ldots,N$ and $t\in[t_i,t_{i+1})$, we have
\[
Y_{t_i}^P-\widehat y_i^{\mathrm{con}}(Q_i)
=(Y_{t_i}^P-y_i^\pi(Q_i))+(y_i^\pi(Q_i)-\widehat y_i^{\mathrm{con}}(Q_i)),
\]
and
\[
Z_t^P-\widehat z_i^{\mathrm{BQL}}(Q_i)
=(Z_t^P-z_i^\pi(Q_i))+(z_i^\pi(Q_i)-\widehat z_i^{\mathrm{BQL}}(Q_i)).
\]
These decompositions separate the total error into two contributions: the time-discretization error and the BQL approximation error relative to the exact time-discrete solution. We first control the former using the classical time-discretization estimate for BSDEs:
\begin{equation}
\max_{0\le i\le N}
\mathbb E\left[
|Y_{t_i}^{P}-y_i^\pi(Q_i)|^2
\right]
+\mathbb E\left[
\sum_{i=0}^{N-1}
\int_{t_i}^{t_{i+1}}
|Z_t^{P}-z_i^\pi(Q_i)|^2\,\mathrm dt
\right]
\le C|\pi|.
\label{eq:time_discretization_bound}
\end{equation}
Here, $C$ satisfies the uniformity convention stated above. A detailed proof of this estimate is provided in Appendix~\ref{app:time_discretization}.

We next estimate the BQL approximation error, starting with the $Z$-error component. Using the auxiliary function $\bar z_i$, we decompose the error as
\[
z_i^\pi-\widehat z_i^{\mathrm{BQL}}
=(z_i^\pi-\bar z_i)
+(\bar z_i-\widehat z_i^{\mathrm{BQL}}).
\label{eq:Z_error_decomposition}
\]
The second term is the local $Z$-BQL quadrature error, measured by $r_{i,Z}$. The first term represents the error propagated from the subsequent value approximation.
To quantify this propagated contribution, define, for $i=0,\ldots,N-1$,
\[
\mathcal V_i^{\mathrm{fluc}}
:=\mathbb E[
| U_{i+1} - \mathbb E_i[U_{i+1}] |^2],
\]
where $U_{i+1}:=(y_{i+1}^{\pi}-\widehat y_{i+1}^{\mathrm{con}})(Q_{i+1})$.
The subsequent $Y$-error $U_{i+1}$ can be decomposed as conditional mean $\mathbb E[U_{i+1}]$ and fluctuation $U_{i+1}-\mathbb E[U_{i+1}]$.
The first term is predictable given the current information, whereas the second represents the remaining random fluctuation. Since the Brownian increment has zero conditional mean, only the fluctuation component can contribute to the current $Z$-error through its interaction with the Brownian increment. Thus, $\mathcal V_i^{\mathrm{fluc}}$ measures the mean-square magnitude of this fluctuation.

One-step $Z$-error established in Lemma~\ref{lem:app_Z_error} therefore yields
\[
\Delta t_i
\|z_i^\pi-\widehat z_i^{\mathrm{BQL}}\|_{L^2(\rho_i)}^2
\le
2\mathcal V_i^{\mathrm{fluc}}
+2\Delta t_i r_{i,Z}.
\]
Thus, the current $Z$-error is controlled by the local quadrature error $r_{i,Z}$ and the propagated conditional fluctuation $\mathcal V_i^{\mathrm{fluc}}$.

It remains to bound $\mathcal V_i^{\mathrm{fluc}}$. Combining the Lipschitz continuity of the generator with the conditional variance decomposition and Young's inequality yields the one-step $Y$-error established in Lemma~\ref{lem:app_Y_error}:
\[
\mathbb E\left[
|(y_i^\pi-\widehat y_i^{\mathrm{con}})(Q_i)|^2
\right]
+c_0\mathcal V_i^{\mathrm{fluc}}\le
(1+C\Delta t_i)
\mathbb E\left[
|(y_{i+1}^\pi-\widehat y_{i+1}^{\mathrm{con}})(Q_{i+1})|^2
\right]+
C\Delta t_i r_{i,Z}
+\frac{C}{\Delta t_i}
(r_{i,Y}+r_{i,\mathrm{con}}),
\]
where $c_0>0$ is independent of the time partition. This estimate serves two purposes: dropping the nonnegative fluctuation term gives a closed recursion for the one-step $Y$-error, while retaining it provides the cumulative control needed for the $Z$-error.

We first give one-step $Y$-error. Since $\mathcal V_i^{\mathrm{fluc}}\ge0$, dropping the fluctuation term yields
\begin{equation}
\begin{aligned}
&\mathbb E\left[
|(y_i^\pi-\widehat y_i^{\mathrm{con}})(Q_i)|^2
\right]\le
(1+C\Delta t_i)
\mathbb E\left[
|(y_{i+1}^\pi-\widehat y_{i+1}^{\mathrm{con}})(Q_{i+1})|^2
\right]+
C\Delta t_i r_{i,Z}
+\frac{C}{\Delta t_i}
(r_{i,Y}+r_{i,\mathrm{con}}).
\end{aligned}
\label{eq:one_step_Y_error}
\end{equation}
The right-hand side contains only the propagated $Y$-error and the three local error contributions. Since $y_N^\pi=\widehat y_N^{\mathrm{con}}=g_p$, the terminal $Y$-error vanishes. Backward iteration and the discrete Gronwall inequality give
\begin{equation}
\max_{0\le i\le N}
\mathbb E\left[
|(y_i^\pi-\widehat y_i^{\mathrm{con}})(Q_i)|^2
\right]\le
C\sum_{i=0}^{N-1}
\left[
\Delta t_i r_{i,Z}
+\frac{r_{i,Y}+r_{i,\mathrm{con}}}{\Delta t_i}
\right].
\label{eq:global_Y_error}
\end{equation}
Hence, the global $Y$-error component is controlled by the accumulated local approximation errors.

We next bound $\mathcal V_i^{\mathrm{fluc}}$, which appears in the one-step $Z$-error. Retaining the fluctuation term in the one-step $Y$-error and summing over $i=0,\ldots,N-1$ gives
\begin{align*}
c_0\sum_{i=0}^{N-1}\mathcal V_i^{\mathrm{fluc}}
\le{}&
\sum_{i=0}^{N-1}
\Big[
(1+C\Delta t_i)
\mathbb E\left[
|(y_{i+1}^\pi-\widehat y_{i+1}^{\mathrm{con}})(Q_{i+1})|^2
\right]-
\mathbb E\left[
|(y_i^\pi-\widehat y_i^{\mathrm{con}})(Q_i)|^2
\right]
\Big]\\
&+
C\sum_{i=0}^{N-1}
\left[
\Delta t_i r_{i,Z}
+\frac{r_{i,Y}+r_{i,\mathrm{con}}}{\Delta t_i}
\right].
\end{align*}
The consecutive $Y$-error terms telescope, apart from the contribution associated with the stability factor $C\Delta t_i$. Since the terminal error vanishes and all $Y$-error terms are nonnegative, we obtain
\[
c_0\sum_{i=0}^{N-1}\mathcal V_i^{\mathrm{fluc}}
\le{}
CT\max_{0\le j\le N}
\mathbb E\left[
|(y_j^\pi-\widehat y_j^{\mathrm{con}})(Q_j)|^2
\right]+
C\sum_{i=0}^{N-1}
\left[
\Delta t_i r_{i,Z}
+\frac{r_{i,Y}+r_{i,\mathrm{con}}}{\Delta t_i}
\right].
\]
Substituting the global $Y$-error estimate \eqref{eq:global_Y_error} into this inequality yields
\[
\sum_{i=0}^{N-1}\mathcal V_i^{\mathrm{fluc}}
\le
C\sum_{i=0}^{N-1}
\left[
\Delta t_i r_{i,Z}
+\frac{r_{i,Y}+r_{i,\mathrm{con}}}{\Delta t_i}
\right].
\]
Consequently, the accumulated conditional fluctuation is controlled by the same local errors.

Returning to the one-step $Z$-error, summing over all time levels and applying the preceding bound gives
\begin{equation}
\begin{aligned}
\sum_{i=0}^{N-1}\Delta t_i
\|z_i^\pi-\widehat z_i^{\mathrm{BQL}}\|_{L^2(\rho_i)}^2
\le
C\sum_{i=0}^{N-1}
\left[
\Delta t_i r_{i,Z}
+\frac{r_{i,Y}+r_{i,\mathrm{con}}}{\Delta t_i}
\right].
\end{aligned}
\label{eq:global_Z_error}
\end{equation}
We have therefore controlled the error of the BQL approximation relative to the exact time-discrete solution.

Finally, applying $|a+b|^2\le 2|a|^2+2|b|^2$ to the error decompositions at the beginning of the proof, we combine the time-discretization estimate \eqref{eq:time_discretization_bound} with the global $Y$ and $Z$-error estimates in \eqref{eq:global_Y_error} and \eqref{eq:global_Z_error}, respectively. This yields \eqref{eq:total_structural_error}.
\end{proof}

\textbf{Remark}.
The error estimate in Theorem~\ref{thm:total_error} is analogous in its overall structure to classical error estimates for nonparametric BSDE-based numerical schemes. In particular, Germain, Pham, and Warin~\cite{GermainPhamWarin2022} established an approximation error estimate for the MDBDP scheme (Theorem 3.3 therein). Let $\mathcal X$ denote the continuous forward diffusion process and $X_i$ its Euler--Maruyama approximation. Denote by $(\widehat{\mathcal U}_i,\widehat{\mathcal Z}_i)$ the neural network approximations of the BSDE components $(Y_{t_i},Z_{t_i})$ produced by MDBDP. Their error estimate takes the form
\begin{align}
&\sup_{0\leq i\leq N}\mathbb E\left[
|Y_{t_i}-\widehat{\mathcal U}_i(X_i)|^2\right]
+\mathbb E\left[
\sum_{i=0}^{N-1}\int_{t_i}^{t_{i+1}}
|Z_s-\widehat{\mathcal Z}_i(X_i)|_2^2\,\mathrm ds
\right]
\nonumber\\
&\quad\leq C\left(
\mathbb E[|g(\mathcal X_T)-g(X_N)|^2]
+|\pi|+\varepsilon^Z(\pi)
+\sum_{i=0}^{N-1}
(\varepsilon_i^{1,y}+\Delta t_i\varepsilon_i^{1,z})
\right),
\label{eq:MDBDP_error_comparison}
\end{align}
where $\varepsilon_Z(\pi)$ denotes the time-regularity error of the BSDE component $Z$, while $\varepsilon_i^y$ and $\varepsilon_i^z$ denote the corresponding neural-network approximation errors. 
Both estimates separate the time-discretization error from the errors introduced by the numerical approximation of the backward recursion. In the BQL framework, however, the latter are further decomposed into local quadrature errors and the continuation-surrogate error, yielding the distinct contributions $r_{i,Z}$, $r_{i,Y}$, and $r_{i,\mathrm{con}}$ in Theorem~\ref{thm:total_error}. This comparison shows that the proposed estimate follows the established stability-based error-analysis framework for BSDE numerical schemes, while explicitly accounting for the additional approximation mechanisms introduced by BQL.

The estimate separates the total error into the time-discretization contribution and three local approximation contributions. In particular, the $Z$-quadrature error is weighted by $\Delta t_i$, whereas the $Y$-quadrature and continuation errors are weighted by $1/\Delta t_i$. Thus, the local errors need not enter the global estimate at the same scale; their time-step dependence must be quantified before a convergence rate can be concluded. Section~\ref{subsec:local_error_control} establishes these local bounds, and Section~\ref{subsec:convergence} combines them to obtain the global convergence result.

\subsection{Local Quadrature and Continuation Errors}
\label{subsec:local_error_control}

The global structural estimate in Theorem~\ref{thm:total_error} expresses the total error in terms of the local $Z$-BQL, $Y$-BQL, and continuation errors. This section establishes quantitative bounds for these three errors. We first analyze the two BQL quadrature errors, separating finite-node approximation from generalization and optimization errors, and then examine the continuation-surrogate error in terms of approximation capacity and training residuals. These estimates provide the local error bounds needed for the subsequent convergence analysis.

\subsubsection{Local BQL Quadrature Error}
\label{subsec:local_BQL_error}

We first establish quantitative bounds for $r_{i,Z}$ and $r_{i,Y}$. At time $t_i$, the continuation surrogate $\widehat y_{i+1}^{\mathrm{con}}$ is held fixed, and the trained approximation $\widehat z_i^{\mathrm{BQL}}$ is also held fixed in the $Y$-stage. For $q=(p,x)$ and a Brownian-increment-scale node $\vartheta\in\mathbb R^d$, define the successor-state map $\mathcal X_i^{\mathrm{node}}(q,\vartheta):=x+\mu_p(t_i,x)\Delta t_i+\sigma_p(t_i,x)\vartheta$.
For $\vartheta\sim\mathcal N(0,\Delta t_i I_d)$, define the two quadrature integrands 
\[
F_i^Z(q,\vartheta):=\widehat y_{i+1}^{\mathrm{con}}\!\left(p,\mathcal X_i^{\mathrm{node}}(q,\vartheta)\right)\frac{\vartheta}{\Delta t_i},
\]
and
\[
F_i^Y(q,\vartheta):=\widehat y_{i+1}^{\mathrm{con}}\!\left(p,\mathcal X_i^{\mathrm{node}}(q,\vartheta)\right)+\Delta t_i f_p\!\left(t_i,x,\widehat y_{i+1}^{\mathrm{con}}\!\left(p,\mathcal X_i^{\mathrm{node}}(q,\vartheta)\right),\widehat z_i^{\mathrm{BQL}}(q)\right).
\]
Their exact Gaussian expectations coincide with the auxiliary functions in Section~\ref{subsec:analytical_setting}. Specifically, set $m_i^Z:=\bar z_i$ and $m_i^Y:=\bar y_i$, so that $m_i^\alpha(q)=\mathbb E_\vartheta[F_i^\alpha(q,\vartheta)]$ for $\alpha\in\{Z,Y\}$. This scaling is consistent with Section~\ref{subsec:backward_quadrature}: a standard Gaussian sample $\zeta$ corresponds to the increment-scale node $\vartheta=\sqrt{\Delta t_i}\,\zeta$.

For $\alpha\in\{Z,Y\}$, let $\vartheta_i^\alpha=(\vartheta_{i,1}^\alpha,\ldots,\vartheta_{i,M}^\alpha)\in(\mathbb R^d)^M$ be a candidate configuration of $M$ quadrature nodes. Define the finite-node approximation and its population risk by
\[
\mathcal E_{i,M}^\alpha(q;\vartheta_i^\alpha):=\frac1M\sum_{m=1}^M F_i^\alpha(q,\vartheta_{i,m}^\alpha),\qquad R_i^\alpha(\vartheta_i^\alpha):=\left\|m_i^\alpha-\mathcal E_{i,M}^\alpha(\cdot;\vartheta_i^\alpha)\right\|_{L^2(\rho_i)}^2.
\]
Let $\widehat\vartheta_i^\alpha$ denote the trained configuration from Section~\ref{subsec:backward_quadrature}. Its finite-node evaluations coincide with the corresponding BQL approximations, and therefore $r_{i,\alpha}=R_i^\alpha(\widehat\vartheta_i^\alpha)$.

To quantify the discrepancy between the training objective and the population risk, let $\vartheta_{i,k}^{\alpha,\mathrm{ref}}\sim\mathcal N(0,\Delta t_i I_d)$, $k=1,\ldots,K_\alpha$, be independent reference nodes, and define $m_{i,K_\alpha}^{\alpha,\mathrm{ref}}(q):=\frac1{K_\alpha}\sum_{k=1}^{K_\alpha}F_i^\alpha(q,\vartheta_{i,k}^{\alpha,\mathrm{ref}})$.
The population training risk is
\[
R_{i,K_\alpha}^\alpha(\vartheta_i^\alpha):=\mathbb E_{\mathrm{ref}}\!\left[\left\|m_{i,K_\alpha}^{\alpha,\mathrm{ref}}-\mathcal E_{i,M}^\alpha(\cdot;\vartheta_i^\alpha)\right\|_{L^2(\rho_i)}^2\right],
\]
where the expectation is over the reference nodes, with the candidate configuration and learned stage functions held fixed. For the sampled parameter--state pairs $q_i^{(h)}:=(p^{(h)},X_i^{(h),\pi})$, $h=1,\ldots,H$, and their realized reference targets, the empirical objective is
\[
\widehat R_{i,H,K_\alpha}^\alpha(\vartheta_i^\alpha):=\frac1H\sum_{h=1}^H\left|m_{i,K_\alpha}^{\alpha,\mathrm{ref},(h)}(q_i^{(h)})-\mathcal E_{i,M}^\alpha(q_i^{(h)};\vartheta_i^\alpha)\right|^2.
\]
Define the generalization discrepancy  and optimization residual by
\[
\varepsilon_{i,\alpha}^{\mathrm{gen}}:=\sup_{\vartheta_i^\alpha\in(\mathbb R^d)^M}\left|\widehat R_{i,H,K_\alpha}^\alpha(\vartheta_i^\alpha)-R_{i,K_\alpha}^\alpha(\vartheta_i^\alpha)\right|,
\qquad 
\delta_{i,\alpha}^{\mathrm{opt}}:=\widehat R_{i,H,K_\alpha}^\alpha(\widehat\vartheta_i^\alpha)-\inf_{\vartheta_i^\alpha\in(\mathbb R^d)^M}\widehat R_{i,H,K_\alpha}^\alpha(\vartheta_i^\alpha)\ge0.
\]
The generalization discrepancy accounts for parameter--state sampling and reference-target sampling, whereas the optimization residual measures the gap between the attained empirical risk and its infimum. Since the node class is unrestricted, finiteness of $\varepsilon_{i,\alpha}^{\mathrm{gen}}$ is an additional condition.

\begin{proposition}[Local BQL error bound]
\label{prop:local_BQL_error}
Under the standing conditions in Section~\ref{subsec:analytical_setting}, suppose that the risks are well defined, that $\varepsilon_{i,Z}^{\mathrm{gen}}$ and $\varepsilon_{i,Y}^{\mathrm{gen}}$ are finite, and that there exist constants $L_{\mathrm{con}},C_{\mathrm{con}}>0$, independent of the time levels, discretization levels, and network sizes, such that
\begin{equation}
|\widehat y_j^{\mathrm{con}}(p,x)-\widehat y_j^{\mathrm{con}}(p,x')|\le L_{\mathrm{con}}|x-x'|,\qquad |\widehat y_j^{\mathrm{con}}(p,x)|\le C_{\mathrm{con}}(1+|x|),
\label{eq:uniform_continuation_regularity}
\end{equation}
for all relevant time levels $j$, parameters $p\in\mathcal P$, and states $x,x'\in\mathbb R^d$.
Then there exists a constant $C>0$, independent of $i$, $M$, and $\pi$, such that
\begin{equation}
r_{i,Z}\le\frac{C}{M\Delta t_i}+2\varepsilon_{i,Z}^{\mathrm{gen}}+\delta_{i,Z}^{\mathrm{opt}},\qquad r_{i,Y}\le\frac{C\Delta t_i}{M}+2\varepsilon_{i,Y}^{\mathrm{gen}}+\delta_{i,Y}^{\mathrm{opt}},
\label{eq:local_BQL_error_bound}
\end{equation}
for $i=0,\ldots,N-1$.
\end{proposition}

\begin{proof}
The proof proceeds in two steps. We first compare the trained configurations with an ideal finite-node benchmark. We then estimate the integrated conditional variances to obtain the time-step dependence.

\medskip
\noindent\textit{Step 1: Comparison with the finite-node benchmark.}

For $\alpha\in\{Z,Y\}$, define the integrated conditional variance
\[
V_i^\alpha:=\mathbb E\!\left[\mathbb E_\vartheta\!\left[\left|F_i^\alpha(Q_i,\vartheta)-m_i^\alpha(Q_i)\right|^2\right]\right].
\]
By the reference-risk identity and the finite-node benchmark established in Appendix~\ref{app:quadrature_identities}, we have
\[
R_{i,K_\alpha}^\alpha(\vartheta_i^\alpha)=R_i^\alpha(\vartheta_i^\alpha)+\frac{V_i^\alpha}{K_\alpha},
\qquad
\inf_{\vartheta_i^\alpha}R_i^\alpha(\vartheta_i^\alpha)\le\frac{V_i^\alpha}{M}.
\]
Since $V_i^\alpha/K_\alpha$ is independent of the candidate configuration, it cancels in the population excess-risk comparison. By the definitions of the generalization discrepancy and optimization residual, it follows that
\begin{equation}
r_{i,\alpha}\le\frac{V_i^\alpha}{M}+2\varepsilon_{i,\alpha}^{\mathrm{gen}}+\delta_{i,\alpha}^{\mathrm{opt}}.
\label{eq:BQL_error_variance_bound}
\end{equation}
It remains to establish the time-step dependence of $V_i^Z$ and $V_i^Y$.

\medskip
\noindent\textit{Step 2: Time-step dependence of the conditional variances.}

We now estimate $V_i^Z$ and $V_i^Y$ separately. The key observation is that the $Z$-integrand contains the factor $\vartheta/\Delta t_i$, whereas the $Y$-integrand depends on the Gaussian increment through the continuation function and the driver. The uniform regularity condition \eqref{eq:uniform_continuation_regularity} allows us to quantify these different scalings.

For the $Z$-integrand, the conditional variance is bounded by the second moment. Writing $\vartheta=\sqrt{\Delta t_i}\,\zeta$, where $\zeta\sim\mathcal N(0,I_d)$, and using the linear-growth bound in \eqref{eq:uniform_continuation_regularity}, we obtain
\begin{equation}
\begin{aligned}
V_i^Z
&=\mathbb E\mathbb E_\zeta\left[
\left|F_i^Z(Q_i,\sqrt{\Delta t_i}\,\zeta)-m_i^Z(Q_i)\right|^2
\right]=\mathbb E\mathbb E_\zeta\left[
\left|F_i^Z(Q_i,\sqrt{\Delta t_i}\,\zeta)\right|^2
\right]
-\mathbb E\left[|m_i^Z(Q_i)|^2\right]\\
&\leq \mathbb E\mathbb E_\zeta\left[
\left|F_i^Z(Q_i,\sqrt{\Delta t_i}\,\zeta)\right|^2
\right]=\frac{1}{\Delta t_i}\mathbb E\mathbb E_\zeta\left[
\left|\widehat y_{i+1}^{\mathrm{con}}\left(P,\mathcal X_i^{\mathrm{node}}(Q_i,\sqrt{\Delta t_i}\,\zeta)\right)\right|^2
|\zeta|^2
\right]\\
&\leq\frac{C}{\Delta t_i}\mathbb E\mathbb E_\zeta\left[
\left(1+\left|\mathcal X_i^{\mathrm{node}}(Q_i,\sqrt{\Delta t_i}\,\zeta)\right|\right)^2
|\zeta|^2
\right]\\
&\leq\frac{C}{\Delta t_i}\mathbb E\mathbb E_\zeta\left[
\left(1+\left|\mathcal X_i^{\mathrm{node}}(Q_i,\sqrt{\Delta t_i}\,\zeta)\right|^2\right)
|\zeta|^2
\right].
\end{aligned}
\label{eq:Z_variance_second_moment}
\end{equation}
The linear-growth conditions on the coefficients imply
\[
|\mathcal X_i^{\mathrm{node}}(q,\sqrt{\Delta t_i}\,\zeta)|^2\le C(1+|x|^2)(1+|\zeta|^2),
\]
where $C$ is independent of $i$ and $\pi$. The uniform second-moment bound for the Euler states and the finite second and fourth moments of the Gaussian variable therefore yield
\[
V_i^Z\le\frac{C}{\Delta t_i}.
\]

For the $Y$-integrand, fix $q=(p,x)$ and let $\zeta'$ be an independent copy of $\zeta$. Since the trained approximation $\widehat z_i^{\mathrm{BQL}}(q)$ is fixed, the two integrand evaluations differ only through the continuation value in the driver and the continuation term itself. By the Lipschitz condition in assumption (iii) of Section~\ref{subsec:analytical_setting}, applied with the $x$ and $z$-arguments fixed, and the uniform Lipschitz bound in  \eqref{eq:uniform_continuation_regularity},
\[
\left|F_i^Y(q,\sqrt{\Delta t_i}\,\zeta)-F_i^Y(q,\sqrt{\Delta t_i}\,\zeta')\right|\leq (1+L_f\Delta t_i)L_{\mathrm{con}}\sqrt{\Delta t_i}
\left|\sigma_p(t_i,x)(\zeta-\zeta')\right|.
\]
Using the independent-copy variance identity, we obtain
\begin{equation}
\begin{aligned}
V_i^Y
&=\frac{1}{2}\mathbb E\mathbb E_{\zeta,\zeta'}\left[
\left|F_i^Y(Q_i,\sqrt{\Delta t_i}\,\zeta)
-F_i^Y(Q_i,\sqrt{\Delta t_i}\,\zeta')\right|^2
\right]\\
&\leq \frac{1}{2}(1+L_f\Delta t_i)^2L_{\mathrm{con}}^2\Delta t_i
\mathbb E\mathbb E_{\zeta,\zeta'}\left[
\left|\sigma_P(t_i,X_i^{P,\pi})(\zeta-\zeta')\right|^2
\right]\\
&=(1+L_f\Delta t_i)^2L_{\mathrm{con}}^2\Delta t_i
\mathbb E\left[\|\sigma_P(t_i,X_i^{P,\pi})\|_F^2\right]\\
&\leq C\Delta t_i.
\end{aligned}
\label{eq:Y_variance_bound}
\end{equation}
where $\|\cdot\|_{\mathrm F}$ denotes the Frobenius norm. The last inequality follows from the coefficient growth condition and the uniform second-moment bound for the Euler states. Hence,
\[
V_i^Z\le\frac{C}{\Delta t_i},\qquad V_i^Y\le C\Delta t_i.
\]
Substituting these estimates into \eqref{eq:BQL_error_variance_bound} proves \eqref{eq:local_BQL_error_bound}.
\end{proof}

The two BQL quadrature stages have different time-step variance scalings. Multiplying the $Z$-error bound by $\Delta t_i$ and the $Y$-error bound by $(\Delta t_i)^{-1}$ gives
\begin{equation}
\Delta t_i r_{i,Z}+\frac{r_{i,Y}}{\Delta t_i}\le\frac{C}{M}+\Delta t_i\left(2\varepsilon_{i,Z}^{\mathrm{gen}}+\delta_{i,Z}^{\mathrm{opt}}\right)+\frac{2\varepsilon_{i,Y}^{\mathrm{gen}}+\delta_{i,Y}^{\mathrm{opt}}}{\Delta t_i}.
\label{eq:weighted_local_BQL_error}
\end{equation}
Summing the finite-node contributions in \eqref{eq:weighted_local_BQL_error} over all time levels gives an overall bound of $CN/M$, while the remaining terms are the accumulated weighted training residuals.

This completes the control of the local BQL quadrature errors in the global structural estimate. We next analyze the continuation error $r_{i,\mathrm{con}}$, which arises from fitting the continuation surrogate to the learned $Y$-BQL function.

\subsubsection{Continuation-Surrogate Error}
\label{subsec:continuation_error}

Unlike the BQL quadrature errors, the continuation error arises from fitting a surrogate to the learned $Y$-BQL function. We decompose it into function-class approximation error and errors due to finite-sample training and numerical optimization.
After completing the BQL stage at time $t_i$, we fix the resulting function $\widehat y_i^{\mathrm{BQL}}$ and fit a continuation surrogate to the mapping $q\mapsto\widehat y_i^{\mathrm{BQL}}(q)$ under the parameter-state distribution $\rho_i$.

\paragraph{Continuation-surrogate error decomposition.}

Let $\mathcal H_i^{\mathrm{con}}$ denote the selected continuation-surrogate class. For any candidate $v\in\mathcal H_i^{\mathrm{con}}$, define the population risk
\[
R_i^{\mathrm{con}}(v):=\left\|\widehat y_i^{\mathrm{BQL}}-v\right\|_{L^2(\rho_i)}^2.
\]
In particular, the trained continuation surrogate satisfies $r_{i,\mathrm{con}}=R_i^{\mathrm{con}}(\widehat y_i^{\mathrm{con}})$. Define the best-approximation error
\begin{equation}
\eta_{i,\mathrm{con}}^{\mathrm{app}}:=\inf_{v\in\mathcal H_i^{\mathrm{con}}}\left(R_i^{\mathrm{con}}(v)\right)^{1/2}.
\label{eq:continuation_best_approximation}
\end{equation}
This quantity characterizes the smallest achievable $L^2(\rho_i)$ approximation error within the selected class, independently of the particular training outcome. Define $\eta_{\mathrm{con},N}^{\mathrm{app}}:=\max_{0\le i<N}\eta_{i,\mathrm{con}}^{\mathrm{app}}$.

Let $\{Q_{i,n}^{\mathrm{con}}\}_{n=1}^{N_{\mathrm{con}}}$ denote the continuation-training samples, assumed to be independent of the samples used in the preceding BQL optimization. The empirical continuation risk is
\[
\widehat R_{i,N_{\mathrm{con}}}^{\mathrm{con}}(v):=\frac1{N_{\mathrm{con}}}\sum_{n=1}^{N_{\mathrm{con}}}\left|\widehat y_i^{\mathrm{BQL}}(Q_{i,n}^{\mathrm{con}})-v(Q_{i,n}^{\mathrm{con}})\right|^2.
\]
Define the generalization and optimization residuals by
\[\varepsilon_{i,\mathrm{con}}^{\mathrm{gen}}:=\sup_{v\in\mathcal H_i^{\mathrm{con}}}\left|\widehat R_{i,N_{\mathrm{con}}}^{\mathrm{con}}(v)-R_i^{\mathrm{con}}(v)\right|, \qquad 
\delta_{i,\mathrm{con}}^{\mathrm{opt}}:=\widehat R_{i,N_{\mathrm{con}}}^{\mathrm{con}}(\widehat y_i^{\mathrm{con}})-\inf_{v\in\mathcal H_i^{\mathrm{con}}}\widehat R_{i,N_{\mathrm{con}}}^{\mathrm{con}}(v)\ge0.
\]
The generalization residual measures the discrepancy between empirical and population risks over the selected class, whereas the optimization residual measures the gap between the attained empirical risk and its infimum. As in the BQL analysis, the generalization residual is treated as an abstract finite quantity unless additional complexity and boundedness assumptions are imposed on $\mathcal H_i^{\mathrm{con}}$.

\begin{proposition}[Continuation-surrogate error decomposition]
\label{prop:continuation_error_decomposition}
Assume that $\varepsilon_{i,\mathrm{con}}^{\mathrm{gen}}<\infty$. Then
\begin{equation}
r_{i,\mathrm{con}}\le(\eta_{i,\mathrm{con}}^{\mathrm{app}})^2+2\varepsilon_{i,\mathrm{con}}^{\mathrm{gen}}+\delta_{i,\mathrm{con}}^{\mathrm{opt}}.
\label{eq:continuation_surrogate_error_bound}
\end{equation}
\end{proposition}

\begin{proof}
The proof is provided in Appendix~\ref{app:continuation-risk-comparison}.
\end{proof}

Proposition~\ref{prop:continuation_error_decomposition} reduces the continuation-error analysis to controlling the best-approximation error within the selected function class. In this work, we adopt an additive hybrid continuation surrogate consisting of a Gaussian RBF branch and an FNN branch, whose outputs are added to form the continuation approximation. The two branches are jointly trained. We next establish a best-approximation error bound for the resulting hybrid continuation class under suitable regularity and moment conditions.

\paragraph{Hybrid approximation bound.}
Let $S_i$ denote the network size bound for the FNN branch, controlling both the total number of weights and the number of computation units, and let $N_c$ denote the number of centers in the Gaussian RBF branch. Denote the admissible additive hybrid continuation function class by $\mathcal H_i^{\mathrm{con}}=\mathcal H_i^{\mathrm{hyb}}(S_i,N_c)$, and let $\eta_{i,\mathrm{con}}^{\mathrm{app}}$ be its best-approximation error defined in \eqref{eq:continuation_best_approximation}.

For the convergence analysis, we consider a sequence of uniform partitions $\pi_N=\{t_i=iT/N\}_{i=0}^N$. Throughout this paragraph, we use the notation introduced in the preceding sections, with its dependence on $N$ understood when $\pi=\pi_N$. In particular, $\widehat y_i^{\mathrm{BQL}}$, $Q_i$, $S_i$, and $\eta_{i,\mathrm{con}}^{\mathrm{app}}$ correspond to $\widehat y_{i,N}^{\mathrm{BQL}}$, $Q_{i,N}$, $S_{i,N}$, and $\eta_{i,N,\mathrm{con}}^{\mathrm{app}}$ in Appendix~\ref{app:architecture}, respectively.

Let $d_q=d_p+d$. Under the uniform regularity and moment assumptions stated in Appendix~\ref{app:architecture}, for every $R\geq1$ and sufficiently large $S_i$, the hybrid approximation error satisfies
\begin{equation}
(\eta_{i,\mathrm{con}}^{\mathrm{app}})^2
\leq C(1+R)^4
\left(\frac{S_i}{\log(1+S_i)}\right)^{-2/d_q}
+C_\tau R^{-\delta},
\qquad R\geq1,
\label{eq:hybrid_approximation_bound}
\end{equation}
where $C$ and $C_\tau$ are independent of $N$, $i$, $R$, $S_i$, and $N_c$. The proof is provided in Appendix~\ref{app:architecture}.

\subsection{Consistency and Convergence Analysis}
\label{subsec:convergence}

We now combine the structural estimate in Theorem~\ref{thm:total_error} with the local error bounds from Section~\ref{subsec:local_error_control} to establish the convergence rate of the proposed scheme. We first derive a global error estimate in terms of the quadrature-node number, the continuation approximation error, and the accumulated training residuals. We then use the architecture-dependent approximation bounds in Appendix~\ref{app:architecture} to obtain sufficient conditions for preserving the mean-square order of the time discretization.

Consider uniform partitions $\pi_N$ with $\Delta t=T/N$. To collect the effects of practical training, define
\begin{equation}
\mathcal R_{\mathrm{train},N}:=\sum_{i=0}^{N-1}\left[\Delta t\left(2\varepsilon_{i,Z}^{\mathrm{gen}}+\delta_{i,Z}^{\mathrm{opt}}\right)+\frac{2\varepsilon_{i,Y}^{\mathrm{gen}}+\delta_{i,Y}^{\mathrm{opt}}+2\varepsilon_{i,\mathrm{con}}^{\mathrm{gen}}+\delta_{i,\mathrm{con}}^{\mathrm{opt}}}{\Delta t}\right].
\label{eq:accumulated_training_residual}
\end{equation}

\begin{theorem}[Mean-square convergence rate]
\label{thm:practical_global_error}
Suppose that the standing assumptions of Section~\ref{subsec:analytical_setting} and the regularity conditions required for Proposition~\ref{prop:local_BQL_error} hold uniformly along the sequence of discretizations and learned functions. Assume that the relevant risks and training residuals are finite. Furthermore, suppose that there exist constants $c_M>0$, $C_{\mathrm{app}}>0$, and $C_{\mathrm{train}}>0$, independent of $N$, such that
\begin{equation}
M\ge c_MN^2,\qquad (\eta_{\mathrm{con},N}^{\mathrm{app}})^2\le C_{\mathrm{app}}N^{-3},\qquad \mathcal R_{\mathrm{train},N}\le C_{\mathrm{train}}N^{-1}.
\label{eq:first_order_accuracy_scaling}
\end{equation}
Then the approximation $(\widehat y_i^{\mathrm{con}},\widehat z_i^{\mathrm{BQL}})_{0\le i<N}$ satisfies
\begin{equation}
\max_{0\le i\le N}\mathbb E\!\left[|Y_{t_i}^{P}-\widehat y_i^{\mathrm{con}}(Q_i)|^2\right]+\sum_{i=0}^{N-1}\mathbb E\int_{t_i}^{t_{i+1}}|Z_t^{P}-\widehat z_i^{\mathrm{BQL}}(Q_i)|^2\,\mathrm dt\le C\Delta t,
\label{eq:uniform_convergence_bound}
\end{equation}
where $C$ is independent of $N$ and $M$.
\end{theorem}

\begin{proof}
Substituting the local error bounds \eqref{eq:local_BQL_error_bound} and \eqref{eq:continuation_surrogate_error_bound} into the structural estimate \eqref{eq:total_structural_error} gives
\begin{equation}
\mathcal E_{\mathrm{total}}(\pi)\le C\left[\frac{1}{N}+\frac{N}{M}+N^2(\eta_{\mathrm{con},N}^{\mathrm{app}})^2+\mathcal R_{\mathrm{train},N}\right].
\label{eq:practical_convergence_bound}
\end{equation}
Here, $N/M$ results from accumulating the finite-node quadrature errors over all time levels, while the continuation approximation contribution follows from
\[
\sum_{i=0}^{N-1}\frac{(\eta_{i,\mathrm{con}}^{\mathrm{app}})^2}{\Delta t}\le\frac{N^2}{T}(\eta_{\mathrm{con},N}^{\mathrm{app}})^2.
\]
The remaining generalization and optimization contributions are collected in $\mathcal R_{\mathrm{train},N}$.

By \eqref{eq:first_order_accuracy_scaling}, each term on the right-hand side of \eqref{eq:practical_convergence_bound} is bounded by a constant multiple of $N^{-1}$. Consequently, $\mathcal E_{\mathrm{total}}(\pi)\le C/N=C\Delta t$, where the fixed factor $T^{-1}$ is absorbed into the constant. This proves \eqref{eq:uniform_convergence_bound}.
\end{proof}

\paragraph{Sufficient conditions for the hybrid architecture.}
We now translate the continuation-approximation requirement in \eqref{eq:first_order_accuracy_scaling} into sufficient conditions on the components of the hybrid architecture adopted in BQL. Accordingly, we require
\begin{equation}
(\eta_{\mathrm{con},N}^{\mathrm{app}})^2
\leq C_{\mathrm{app}}N^{-3},
\label{eq:continuation_approximation_requirement}
\end{equation}
where $C_{\mathrm{app}}>0$ is independent of $N$. We next derive sufficient conditions for this requirement using the hybrid approximation bound established in \eqref{eq:hybrid_approximation_bound}.

To control the truncation error, choose $R=N^{3/\delta}$. Then $R\geq1$ and $C_\tau R^{-\delta}=C_\tau N^{-3}$. Moreover, since $(1+R)^4\leq16R^4=16N^{12/\delta}$, the hybrid approximation bound \eqref{eq:hybrid_approximation_bound} gives
\begin{equation}
(\eta_{i,\mathrm{con}}^{\mathrm{app}})^2\leq CN^{12/\delta}\left(\frac{S_{i}}{\log(1+S_{i})}\right)^{-2/d_q}+CN^{-3}.
\label{eq:hybrid_approximation_after_truncation}
\end{equation}
Here, $C$ and $C_\tau$ are independent of $N$ and $i$.

To ensure that the first term on the right-hand side is also of order $N^{-3}$, it suffices to require
\[
\left(\frac{S_{i}}{\log(1+S_{i})}\right)^{-2/d_q}\leq cN^{-(3+12/\delta)},
\]
or equivalently,
\begin{equation}
\frac{S_{i}}{\log(1+S_{i})}\geq cN^{\frac{3d_q}{2}(1+4/\delta)}.
\label{eq:S_i,N}
\end{equation}
A sufficient condition for this inequality is
\begin{equation}
\min_{0\leq i<N}S_{i}\geq C_SN^{\frac{3d_q}{2}(1+4/\delta)}\log(1+N),
\label{eq:hybrid_network_size_condition}
\end{equation}
where $C_S>0$ is sufficiently large and independent of $N$. Indeed, the function $s\mapsto s/\log(1+s)$ is increasing for $s>0$. 
Substituting this estimate\eqref{eq:S_i,N} into \eqref{eq:hybrid_approximation_after_truncation}, we obtain
\[
(\eta_{i,\mathrm{con}}^{\mathrm{app}})^2
\leq CN^{12/\delta}N^{-3(1+4/\delta)}+CN^{-3}\leq CN^{-3},
\]
where the constant $C$ is independent of $N$ and $i$. Taking the maximum over all time levels yields
\begin{equation}
(\eta_{\mathrm{con},N}^{\mathrm{app}})^2:=\max_{0\leq i<N}(\eta_{i,\mathrm{con}}^{\mathrm{app}})^2\leq CN^{-3}.
\label{eq:uniform_hybrid_approximation}
\end{equation}
Thus, the sufficient condition \eqref{eq:continuation_approximation_requirement} is satisfied.

Combining this approximation condition with $M\geq c_MN^2$ and $\mathcal R_{\mathrm{train},N}\leq C_{\mathrm{train}}N^{-1}$, Theorem~\ref{thm:practical_global_error} gives
\[
\mathcal E_{\mathrm{total}}(\pi_N)\leq C\Delta t.
\]
Hence, under the stated uniform regularity, approximation, and training assumptions, the hybrid architecture preserves the mean-square convergence order of the underlying BSDE time discretization. More precisely, the resulting global squared error is bounded by $C\Delta t$, with a constant independent of the time-step number $N$.

\section{Numerical Experiments}
\label{sec:numerical_experiments}

In this section, we evaluate the proposed BQL method from three complementary perspectives: the approximation accuracy of the resulting PDE solutions, the computational efficiency in repeated parameter queries, and the effectiveness of the learned quadrature nodes in approximating conditional expectations.

All numerical experiments were implemented in Python 3.12.3 on a Linux workstation equipped with an Intel Xeon Silver 4314 CPU at 2.40 GHz, 256 GB of system memory, and an NVIDIA GeForce RTX 4090 GPU with 24 GB of memory. TensorFlow 2.21.0 and PyTorch 2.12.0 were used according to the corresponding implementations.

Let $N_{\mathrm{test}}$ denote the total number of scalar solution values included in the common test set, and let $u^{(r)}$ and $\widehat u^{(r)}$, $r=1,\ldots,N_{\mathrm{test}}$, denote the corresponding reference and numerical solutions. We report the root-mean-square error (RMSE), relative $L^2$ error ($\mathrm{Rel.}\ L^2$), and $L^\infty$ error,
\[
\mathrm{RMSE}=\left(\frac{1}{N_{\mathrm{test}}}\sum_{r=1}^{N_{\mathrm{test}}}\left|\widehat u^{(r)}-u^{(r)}\right|^2\right)^{1/2}, 
\quad 
\mathrm{Rel.}\ L^2=\frac{\left(\sum_{r=1}^{N_{\mathrm{test}}}\left|\widehat u^{(r)}-u^{(r)}\right|^2\right)^{1/2}}{\left(\sum_{r=1}^{N_{\mathrm{test}}}\left|u^{(r)}\right|^2\right)^{1/2}},
\]
and
\[
L^\infty=\max_{1\le r\le N_{\mathrm{test}}}\left|\widehat u^{(r)}-u^{(r)}\right|.
\]

The evaluation metrics for quadrature accuracy, including the mean absolute error (MAE) and the improvement factor (IF), together with their detailed computational procedures, are provided in Appendix~\ref{app:quadrature_evaluation}.

\subsection{Common Experimental Setup}
\label{subsec:common_setup}

\paragraph{Common Numerical Settings.}
All numerical experiments are conducted for state dimensions $d\in\{10,20,50,100\}$. We use $N=10$ uniform time steps, a training batch size of $H=256$, and $M=128$ learnable quadrature nodes. The continuation surrogate uses 200 RBF centers and an MLP with hidden-layer widths $[64,64,32]$ and ReLU activations. The $Z$- and $Y$-quadrature components are each trained for 400 epochs, and the continuation surrogate is also trained for 400 epochs. All components are optimized using Adam with an initial learning rate of $10^{-3}$.

During training, the parameter vectors are sampled uniformly from the corresponding parameter domains $\mathcal P$. For each tested parameter value or parameter combination, 20 independent trajectories are generated. Accordingly, $N_{\mathrm{path}}=20$ and $N_{\mathrm{test}}=N_{\mathrm{path}}N=20\times10=200$ for each test parameter configuration. Other experiment-specific settings, including the terminal time, initial state, MC sample size, and test parameter values, are specified separately for each example.

For the local quadrature accuracy comparison, the test parameter values or parameter combinations are chosen to be the same as those used for the solution-accuracy evaluation. For each test parameter configuration, $N_{\mathrm{map}}=20$ independent trajectories are used to evaluate the quadrature errors.

\paragraph{Common Evaluation Protocols.}

\begin{enumerate}
\item[(i)] \noindent\textbf{Approximation accuracy.}
To assess approximation accuracy, we compare the proposed BQL method with P$^2$INN~\cite{ChoEtAl2024} and FBSNN~\cite{Raissi2018FBSNN} on the same test set for each benchmark. For each state dimension, the error statistics are summarized over all tested parameter values or parameter combinations. 

\item[(ii)]
\noindent\textbf{Computational efficiency.}
To evaluate computational efficiency in repeated-query settings, we consider $Q$ parameter queries. BQL and P$^2$INN are trained once over the entire parameter domain, and their trained models are reused for all queries. In contrast, FBSNN is trained independently for each queried parameter value or parameter combination. Accordingly, the cumulative computational cost consists of one training cost and $Q$ per-query solution costs for the parametric methods, and $Q$ training costs and the corresponding per-query solution costs for FBSNN.

\item[(iii)]
\noindent\textbf{Quadrature accuracy.}
To evaluate the effectiveness of learning the quadrature nodes, we compare BQL with untrained Gaussian MC quadrature under the same node budget $M=128$. 
\end{enumerate}

\subsection{Example 1}
\label{subsec:HJB}
We consider a $d$-dimensional parametric Hamilton--Jacobi--Bellman (HJB) equation~\cite{Raissi2018FBSNN}, where the generic parameter $p$ is specialized to the scalar diffusion parameter $\sigma\in\mathcal P$ and $u(t,x;\sigma):=u^p(t,x)\big|_{p=\sigma}$:
\begin{equation}
\left\{
\begin{aligned}
&\frac{\partial u}{\partial t}(t,x;\sigma)
+\frac{\sigma^2}{2}\Delta_x u(t,x;\sigma)
-\left\|\nabla_x u(t,x;\sigma)\right\|_2^2
=0,
\qquad (t,x)\in[0,T)\times\mathbb R^d,\\
&u(T,x;\sigma)
=
g(x)
=
\log\!\left(
\frac{1}{2}
+\frac{1}{2}\|x\|_2^2
\right),
\qquad x\in\mathbb R^d.
\end{aligned}
\right.
\label{eq:HJB_parameterized}
\end{equation}
The objective is to approximate the solution map \((t,x,\sigma)\mapsto u(t,x;\sigma)\) over the parameter domain \(\mathcal P=[1.27,1.55]\). It is worth noting that the nonlinear term \(-\|\nabla_x u\|^2\) leads to a quadratic-growth BSDE generator. Specifically, under the standard correspondence \(z=\sigma\nabla_x u\), the generator is given by \(f_\sigma(z)=-\|z\|^2/\sigma^2\), which is locally Lipschitz continuous but does not satisfy the global Lipschitz condition imposed in Section 4.1. Therefore, this HJB benchmark falls outside the scope of the theoretical convergence results established in Section 4. We include it as an additional numerical test to examine the practical performance of BQL on nonlinear problems beyond the current theoretical assumptions.

The Cole--Hopf transformation yields
\begin{equation}
u(t,x;\sigma)
=
-\frac{\sigma^2}{2}
\log\!\left(
\mathbb E_{\zeta\sim\mathcal N(0,I_d)}
\left[
\exp\!\left(
-\frac{2}{\sigma^2}
g\!\left(
x+\sigma\sqrt{T-t}\,\zeta
\right)
\right)
\right]
\right).
\label{eq:HJB_exact_representation}
\end{equation}
The reference solution is computed from \eqref{eq:HJB_exact_representation} by MC integration with $N_{\mathrm{ref}}=10^5$ independent Gaussian samples.


\subsubsection{Example-specific Setup}
\label{HJB:experimental_setup}
We set the terminal time to $T=0.5$ and the initial state to $x_0=(0,\ldots,0)^\top\in\mathbb R^d$. Unless otherwise stated, $K_Z=K_Y=1024$ MC samples are used to construct the conditional-expectation targets.

The diffusion amplitude $\sigma$ is treated as a variable parameter that controls the intensity of Brownian fluctuations in the underlying stochastic dynamics. The classical HJB benchmark uses $\sigma=\sqrt{2}$~\cite{Raissi2018FBSNN}. We consider the compact parameter domain $\mathcal P=[1.27,1.55]$, which contains $\sqrt{2}$ and extends approximately $10\%$ below and above it. Testing is performed at $\sigma\in\{1.30,1.34,1.38,1.41,1.45,1.48,1.52\}$ with $N_p=7$.

\subsubsection{Comparison of Approximation Accuracy}
\label{subsec:HJB_accuracy_comparison}
As shown in Table~\ref{HJB_table1}, the three methods achieve comparable accuracy in most cases, although their relative performance varies with the state dimension and error metric. BQL yields the smallest relative $L^2$ error for $d=10$ and the smallest relative $L^2$ error and RMSE for $d=20$, while FBSNN performs better for several higher-dimensional metrics. Overall, BQL remains competitive and outperforms the baselines in several configurations.

\begin{table}[!htbp]
\centering
\small
\caption{Minimum and maximum values of the relative $L^2$ error, RMSE, and $L^\infty$ error for BQL, P$^2$INN, and FBSNN across the tested parameter values for different state dimensions.}
\label{HJB_table1}
\setlength{\tabcolsep}{5pt}
\renewcommand{\arraystretch}{1.15}
\resizebox{\textwidth}{!}{
\begin{tabular}{c c c c c c c c}
\toprule
\textbf{$d$} &
\textbf{Method} &
\textbf{Rel. $L^2$ error (min)} &
\textbf{Rel. $L^2$ error (max)} &
\textbf{RMSE (min)} &
\textbf{RMSE (max)} &
\textbf{$L^\infty$ error (min)} &
\textbf{$L^\infty$ error (max)} \\
\midrule

\multirow{3}{*}{$10$}
& BQL
& \textbf{0.011838}
& \textbf{0.022385}
& \textbf{0.016388}
& \textbf{0.037533}
& 0.109452
& 0.469557 \\
& P$^2$INN
& 0.020974
& 0.031269
& 0.029434
& 0.052428
& \textbf{0.049694}
& \textbf{0.078227} \\
& FBSNN
& 0.035330
& 0.043288
& 0.057757
& 0.059926
& 0.155626
& 0.180865 \\
\midrule

\multirow{3}{*}{$20$}
& BQL
& \textbf{0.007343}
& \textbf{0.016487}
& \textbf{0.017236}
& 0.036880
& \textbf{0.058567}
& 0.123200 \\
& P$^2$INN
& 0.014993
& 0.018245
& 0.034536
& 0.037923
& 0.092185
& 0.106737 \\
& FBSNN
& 0.014200
& 0.017540
& 0.032815
& \textbf{0.036458}
& 0.083464
& \textbf{0.091539} \\
\midrule

\multirow{3}{*}{$50$}
& BQL
& 0.012261
& 0.017329
& 0.040609
& 0.054807
& 0.248245
& 0.476937 \\
& P$^2$INN
& 0.024834
& 0.029594
& 0.075719
& 0.098012
& 0.363888
& 0.439252 \\
& FBSNN
& \textbf{0.004695}
& \textbf{0.006309}
& \textbf{0.014374}
& \textbf{0.019331}
& \textbf{0.034399}
& \textbf{0.055351} \\
\midrule

\multirow{3}{*}{$100$}
& BQL
& 0.010512
& 0.018260
& 0.039736
& 0.070876
& 0.202676
& 0.523967 \\
& P$^2$INN
& 0.038584
& 0.062377
& 0.143520
& 0.251457
& 0.509962
& 0.562124 \\
& FBSNN
& \textbf{0.007556}
& \textbf{0.011163}
& \textbf{0.029327}
& \textbf{0.044405}
& \textbf{0.112885}
& \textbf{0.202146} \\
\bottomrule
\end{tabular}
}
\end{table}

\begin{figure}[!htbp]
\centering
\begin{minipage}[t]{0.48\textwidth}
\centering
\includegraphics[width=\linewidth]{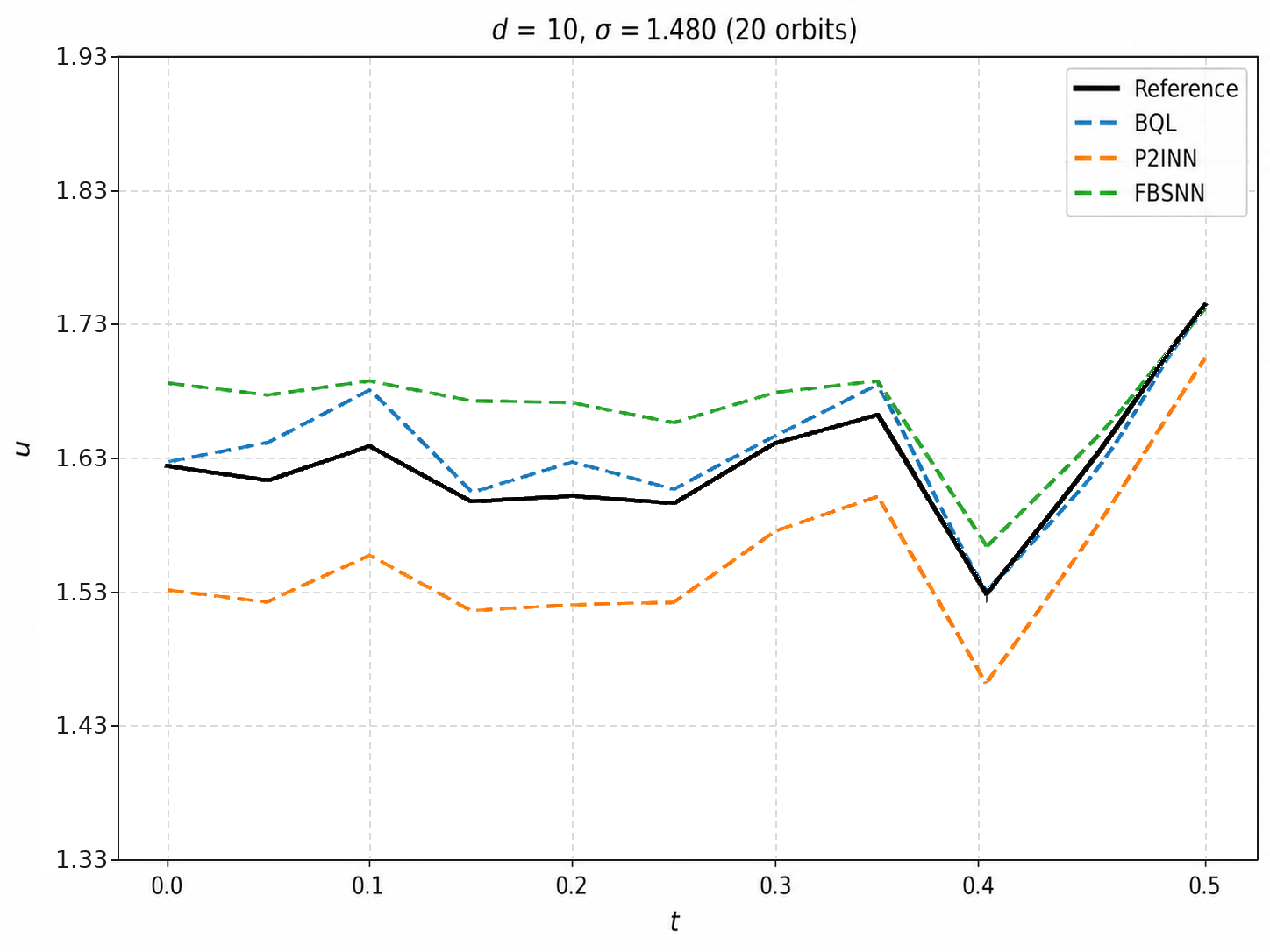}
\end{minipage}
\hfill
\begin{minipage}[t]{0.48\textwidth}
\centering
\includegraphics[width=\linewidth]{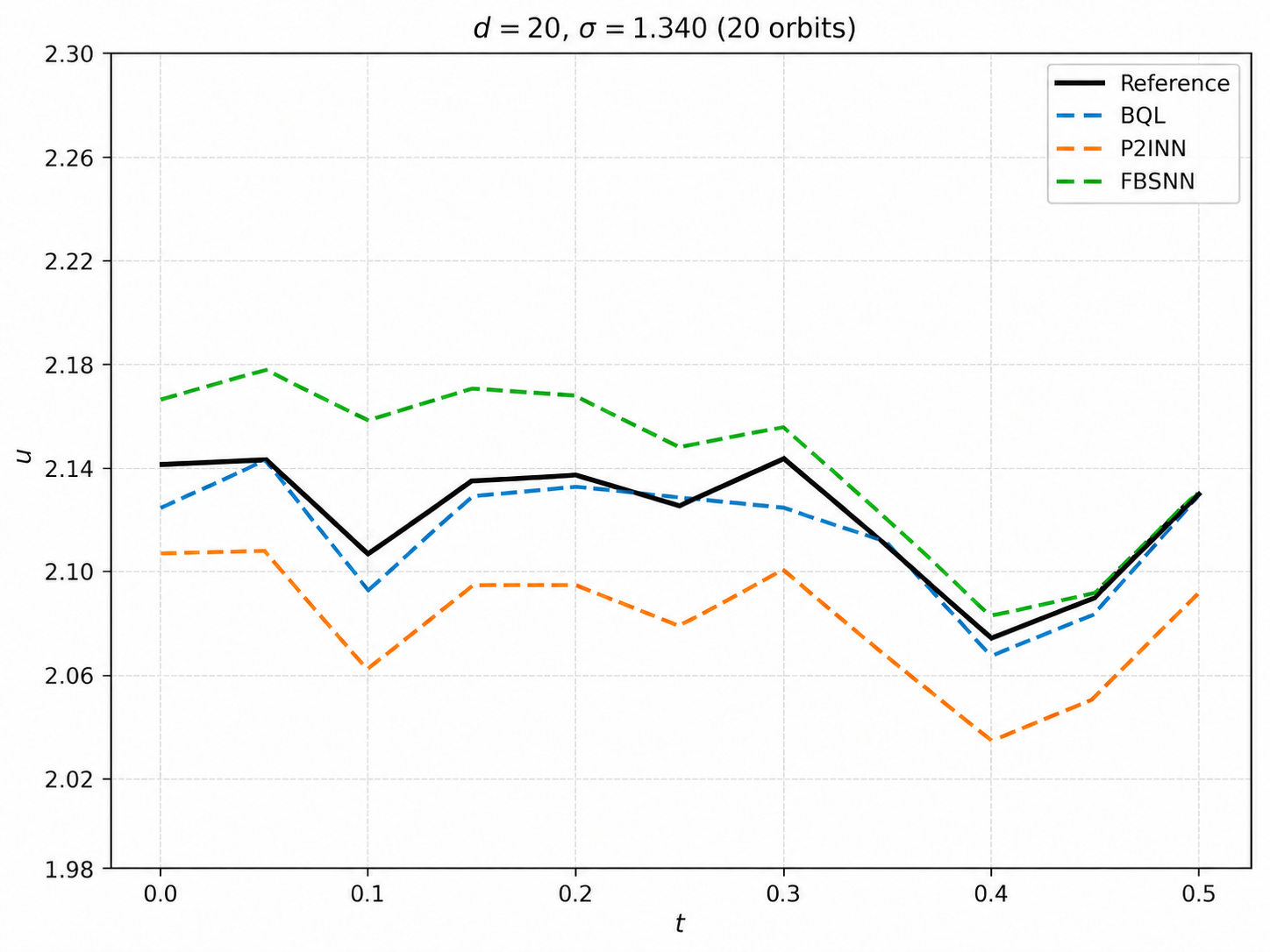}
\end{minipage}

\vspace{0.3cm}

\begin{minipage}[t]{0.48\textwidth}
\centering
\includegraphics[width=\linewidth]{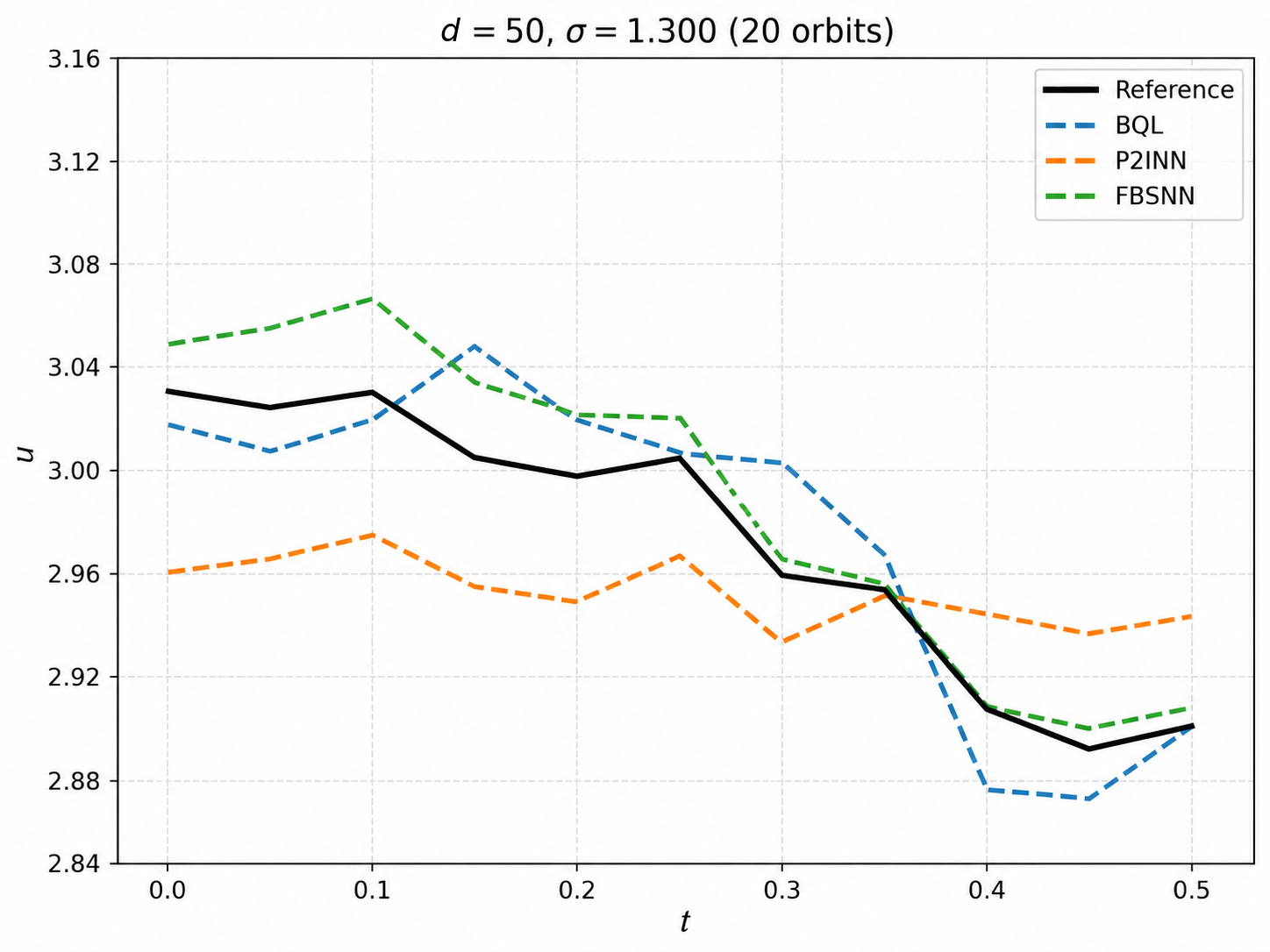}
\end{minipage}
\hfill
\begin{minipage}[t]{0.48\textwidth}
\centering
\includegraphics[width=\linewidth]{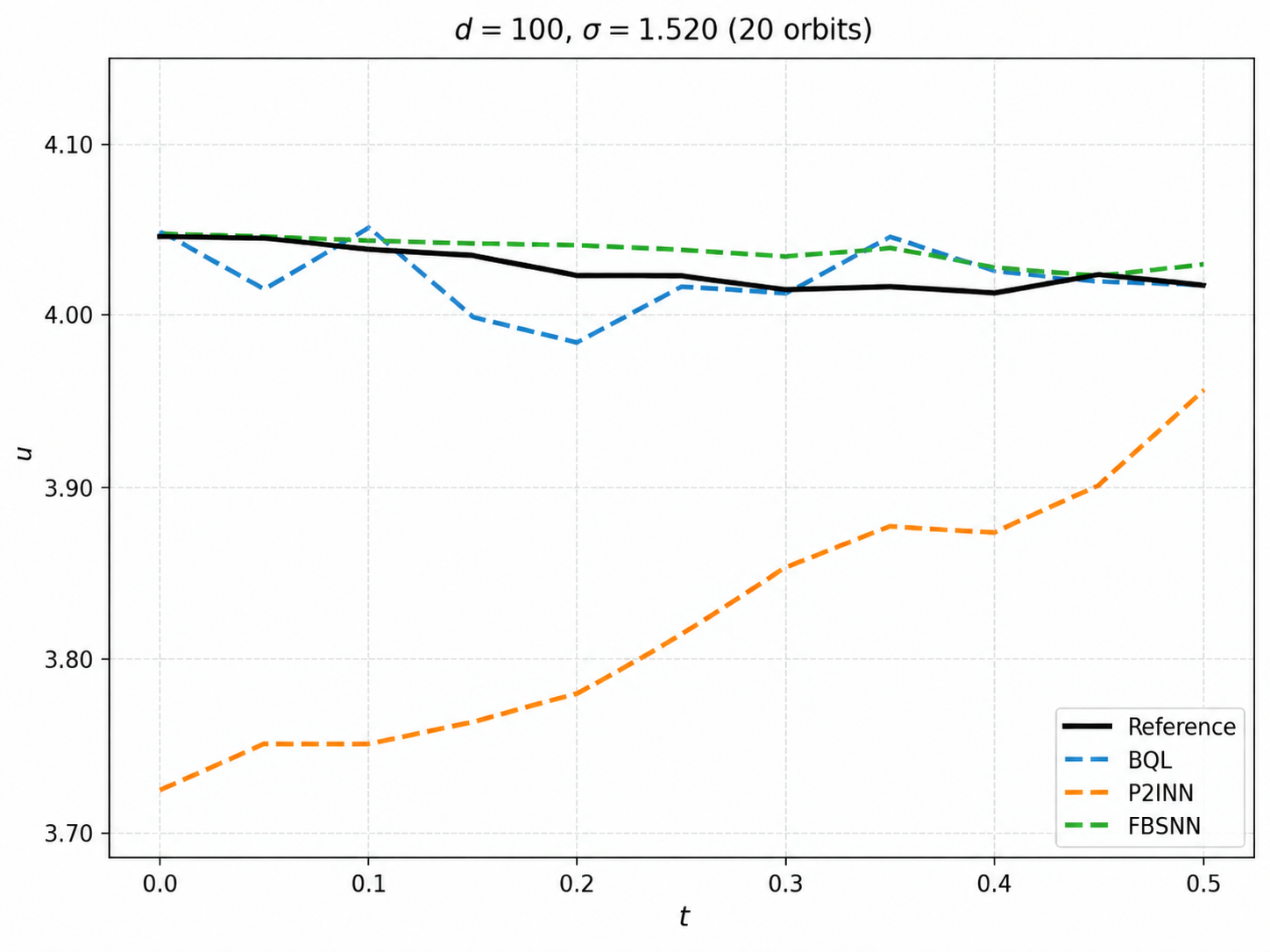}
\end{minipage}

\caption{Comparison of the reference solutions and numerical solutions obtained by BQL, P$^2$INN, and FBSNN for different state dimensions and diffusion parameters: (a) $d=10$, $\sigma=1.480$; (b) $d=20$, $\sigma=1.340$; (c) $d=50$, $\sigma=1.300$; (d) $d=100$, $\sigma=1.520$. At each time point, the plotted values are obtained by averaging the solution evaluations along 20 independent sample trajectories.}
\label{HJB_figure1}
\end{figure}

Figure~\ref{HJB_figure1} further illustrates the solution profiles obtained by the three methods at representative diffusion parameters for $d=10$, $20$, $50$, and $100$. BQL generally captures the temporal evolution of the reference solutions with comparable accuracy to the FBSNN methods, while P$^2$INN exhibits more pronounced deviations in the higher-dimensional case. Together with the results in Table~\ref{HJB_table1}, these comparisons demonstrate that BQL maintains a competitive level of approximation accuracy across different state dimensions, providing a reliable basis for efficient repeated parameter queries.

\subsubsection{Computational Efficiency in the Many-Query Setting}
\label{HJB:time_comparison}

Table~\ref{HJB_table2} shows that BQL achieves the lowest cumulative cost, primarily due to its lower one-time training cost and parameter reuse. For $Q=7$, its cost is about $4.0$--$5.8$ times lower than FBSNN and $44.5$--$208.0$ times lower than P$^2$INN, despite the shorter online solution times of the baselines.

\begin{table}[htbp]
\centering
\small
\caption{Comparison of training time, per-query solution time, and cumulative computational costs for BQL, P$^2$INN, and FBSNN with $Q=7$ parameter queries. All computational times are reported in seconds. The last column reports the cumulative cost relative to that of BQL.}
\label{HJB_table2}
\setlength{\tabcolsep}{4pt}
\resizebox{\textwidth}{!}{
\begin{tabular}{cclccccc}
\toprule
\textbf{$d$} &
\textbf{Method} &
\textbf{Training regime} &
\textbf{Training time} &
\textbf{Solution time per query} &
\textbf{Cumulative time } &
\textbf{Cost ratio to BQL} \\
\midrule

\multirow{3}{*}{$10$}
& BQL
& Parametric
& 452.576486
& 0.003052
& \textbf{452.597850}
& 1.000 \\
& P$^2$INN
& Parametric
& 20139.535884
& 0.000331
& 20139.538201
& 44.498 \\
& FBSNN
& Non-parameter
& 258.532542
& 0.000374
& 1809.730412
& 3.999 \\
\midrule

\multirow{3}{*}{$20$}
& BQL
& Parametric
& 453.422724
& 0.002975
& \textbf{453.443549}
& 1.000 \\
& P$^2$INN
& Parametric
& 27554.359541
& 0.000330
& 27554.361851
& 60.767 \\
& FBSNN
& Non-parameter
& 274.385364
& 0.000368
& 1920.700124
& 4.236 \\
\midrule

\multirow{3}{*}{$50$}
& BQL
& Parametric
& 456.201761
& 0.002977
& \textbf{456.222600}
& 1.000 \\
& P$^2$INN
& Parametric
& 51499.136693
& 0.000327
& 51499.138982
& 112.882 \\
& FBSNN
& Non-parameter
& 310.551435
& 0.000387
& 2173.862754
& 4.765 \\
\midrule

\multirow{3}{*}{$100$}
& BQL
& Parametric
& 456.587821
& 0.002981
& \textbf{456.608688} 
& 1.000 \\
& P$^2$INN
& Parametric
& 94980.169845
& 0.002292
& 94980.185889
& 208.012 \\
& FBSNN
& Non-parameter
& 378.188033
& 0.000362
& 2647.318765
& 5.798 \\
\bottomrule
\end{tabular}
}
\end{table}

\begin{figure}[htbp]
\centering

\begin{minipage}[t]{0.48\textwidth}
\centering
\includegraphics[height=5.2cm]{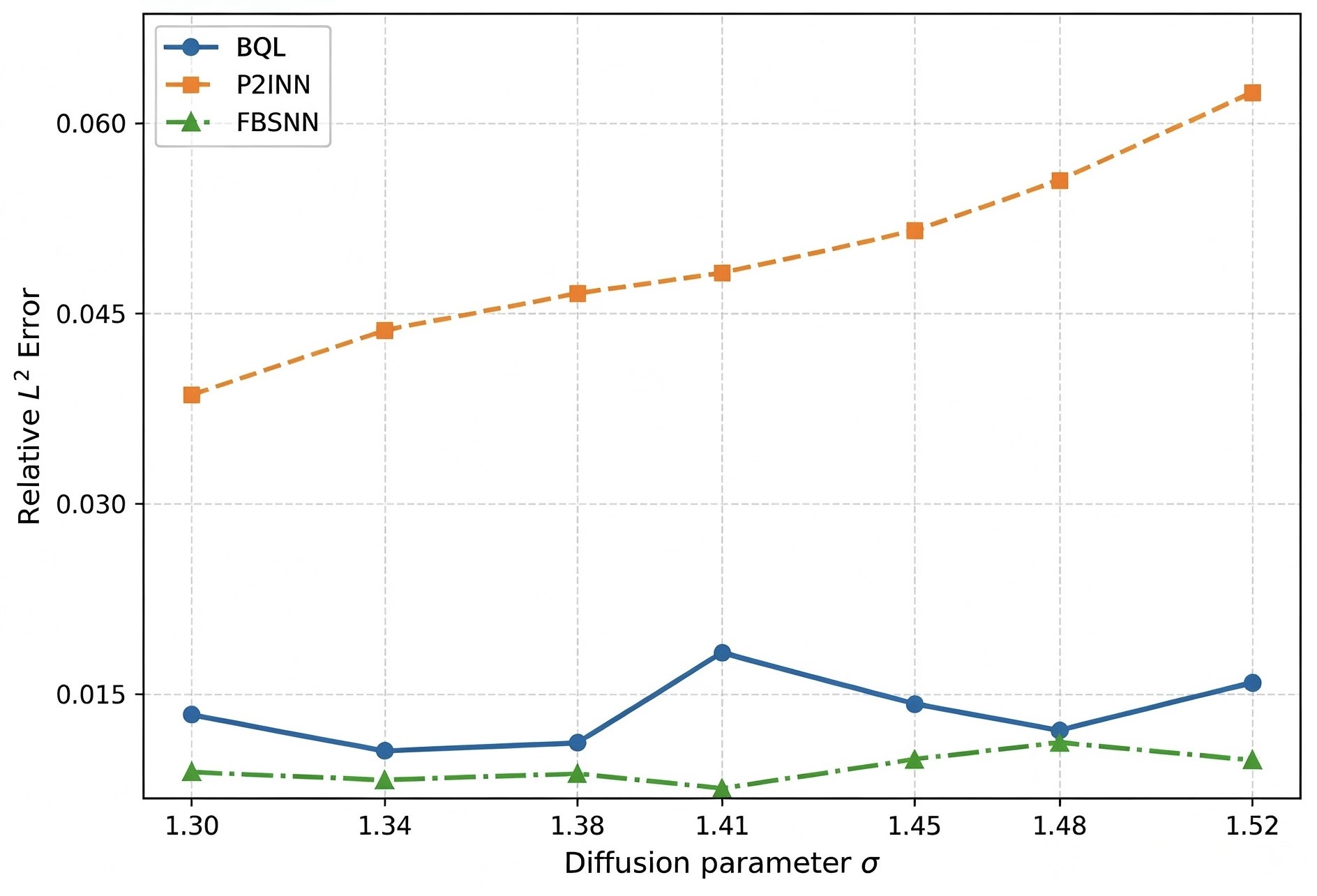}
\vspace{1mm}
\end{minipage}
\hfill
\begin{minipage}[t]{0.48\textwidth}
\centering
\includegraphics[height=5.2cm]{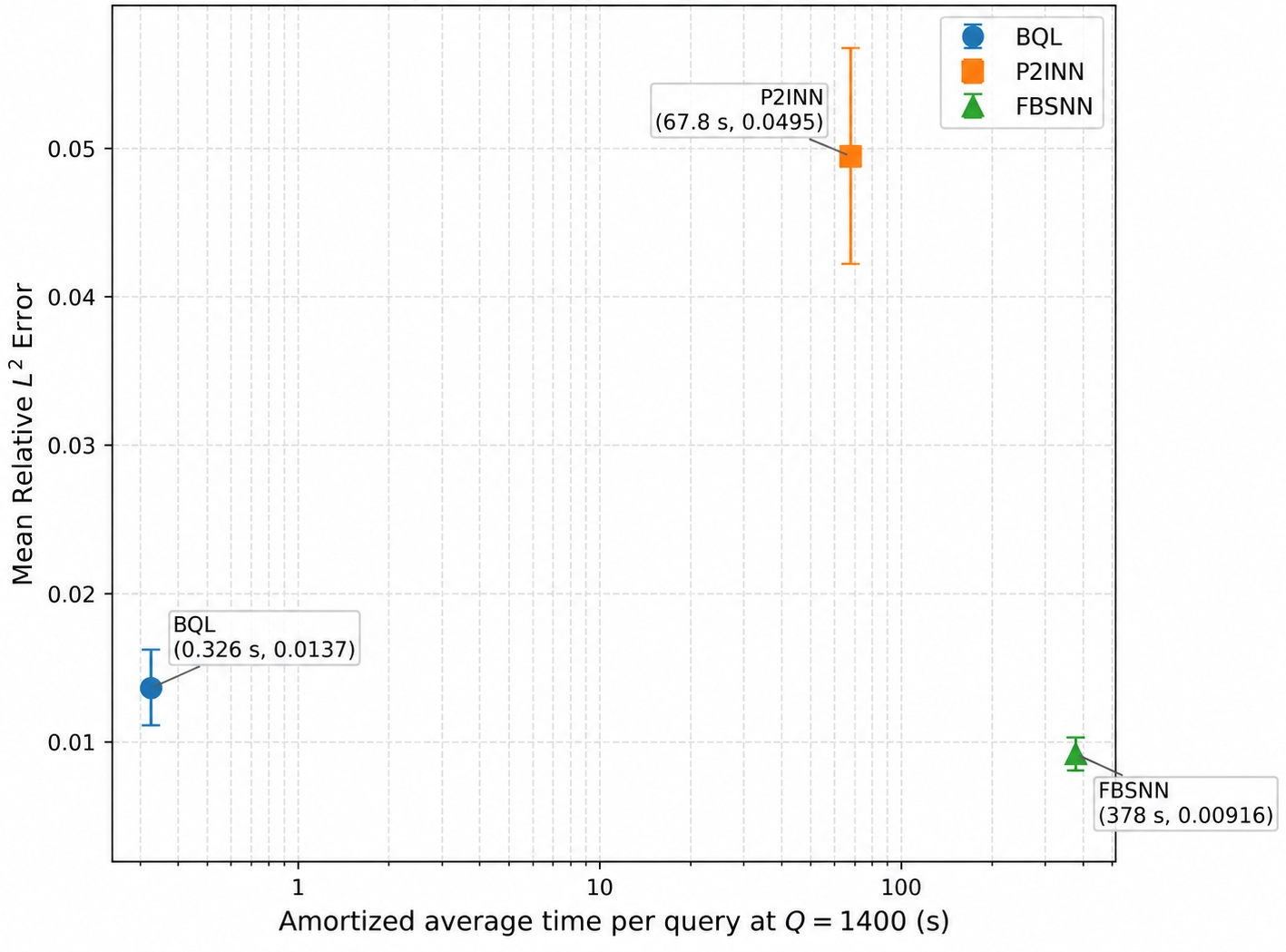}
\vspace{1mm}
\end{minipage}

\vspace{0.3cm}

\makebox[\textwidth][c]{%
\begin{minipage}[t]{0.62\textwidth}
\centering
\includegraphics[width=\linewidth]{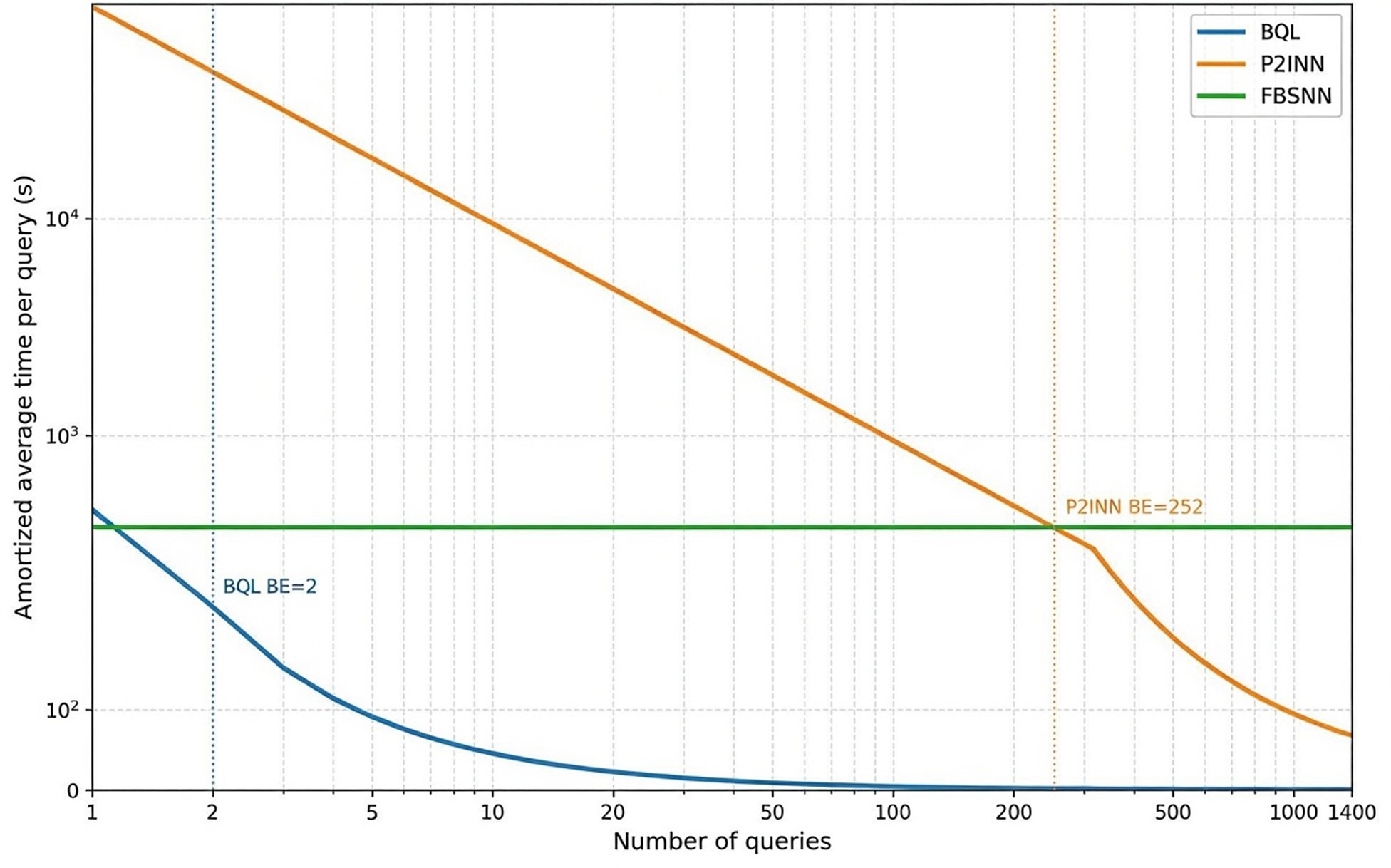}
\vspace{1mm}
\end{minipage}%
}

\caption{Comparison of approximation accuracy and computational efficiency of BQL, P$^2$INN, and FBSNN for $d=100$: (a) relative $L^2$ errors across different parameters; (b) accuracy--efficiency comparison for $Q=1400$ parameter queries; and (c) amortized average computational time per query as a function of the number of parameter queries $Q$.}
\label{HJB_figure2}

\end{figure}

Figure~\ref{HJB_figure2}(a) compares the relative $L^2$ errors across different diffusion parameters. BQL consistently achieves lower errors than P$^2$INN, while maintaining accuracy comparable to FBSNN. At $Q=1400$, as shown in Fig.~\ref{HJB_figure2}(b), the mean relative $L^2$ errors are $0.0137$, $0.0495$, and $0.00916$ for BQL, P$^2$INN, and FBSNN, respectively. Despite its slightly higher error than FBSNN, BQL requires only $0.326$ s per query, compared with $67.8$ s for P$^2$INN and $378$ s for FBSNN, making it approximately $208$ and $1160$ times faster, respectively. Figure~\ref{HJB_figure2}(c) further illustrates the computational advantage of BQL as the number of parameter queries increases. For $d=100$, BQL becomes more efficient than FBSNN when $Q>2$, whereas P$^2$INN reaches this break-even point only when $Q>252$. Moreover, BQL maintains a substantially lower amortized cost than P$^2$INN throughout the tested query range. Similar qualitative trends are observed for the other tested dimensions. Overall, BQL achieves a favorable balance between approximation accuracy and computational efficiency, making it particularly advantageous for repeated parameter queries.

\subsubsection{Effectiveness of Learning the Quadrature Nodes}
\label{HJB:Compare MC}

\begin{figure}[htbp]
\centering
\begin{minipage}[t]{0.48\textwidth}
\centering
\includegraphics[width=\linewidth]{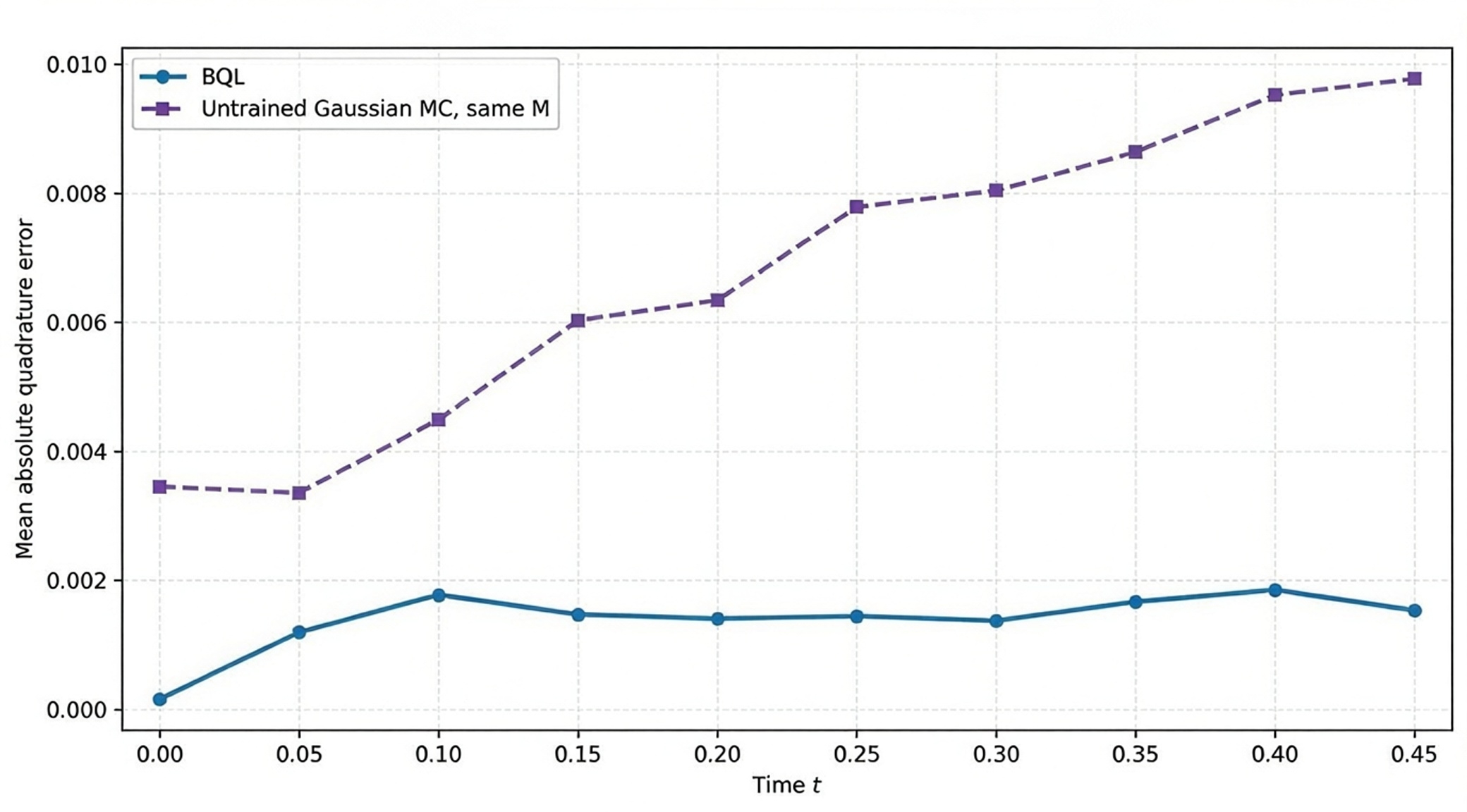}
\end{minipage}
\hfill
\begin{minipage}[t]{0.48\textwidth}
\centering
\includegraphics[width=\linewidth]{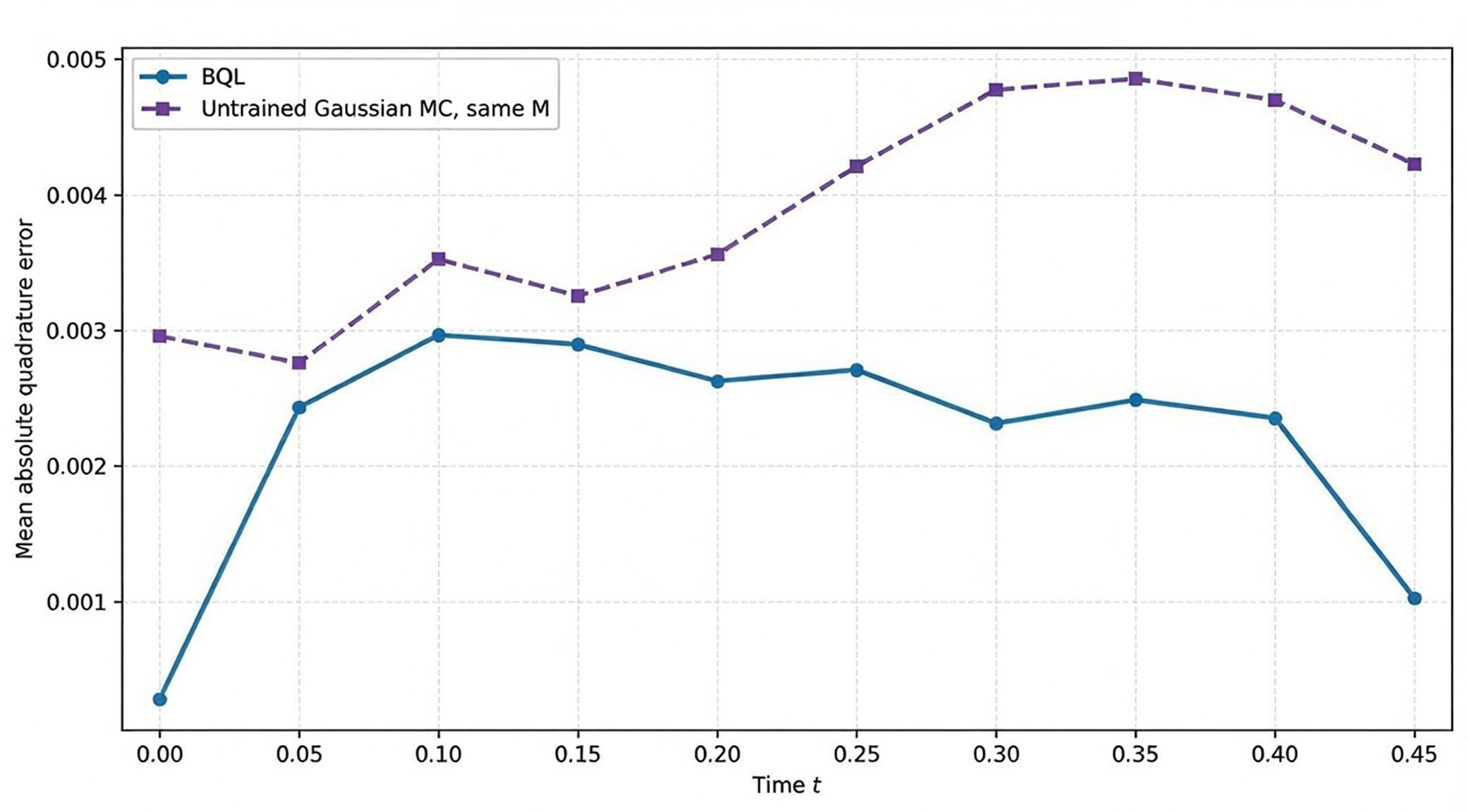}
\end{minipage}

\caption{Comparison of MAEs between BQL and untrained Gaussian MC quadrature for $d=20$ (left) and $d=100$ (right). At each time point, the MAEs are averaged over the $N_p$ test parameter values.}
\label{HJB_figure3-1}
\end{figure}

\begin{figure}

\begin{minipage}[t]{0.96\textwidth}
\centering
\includegraphics[width=\linewidth]{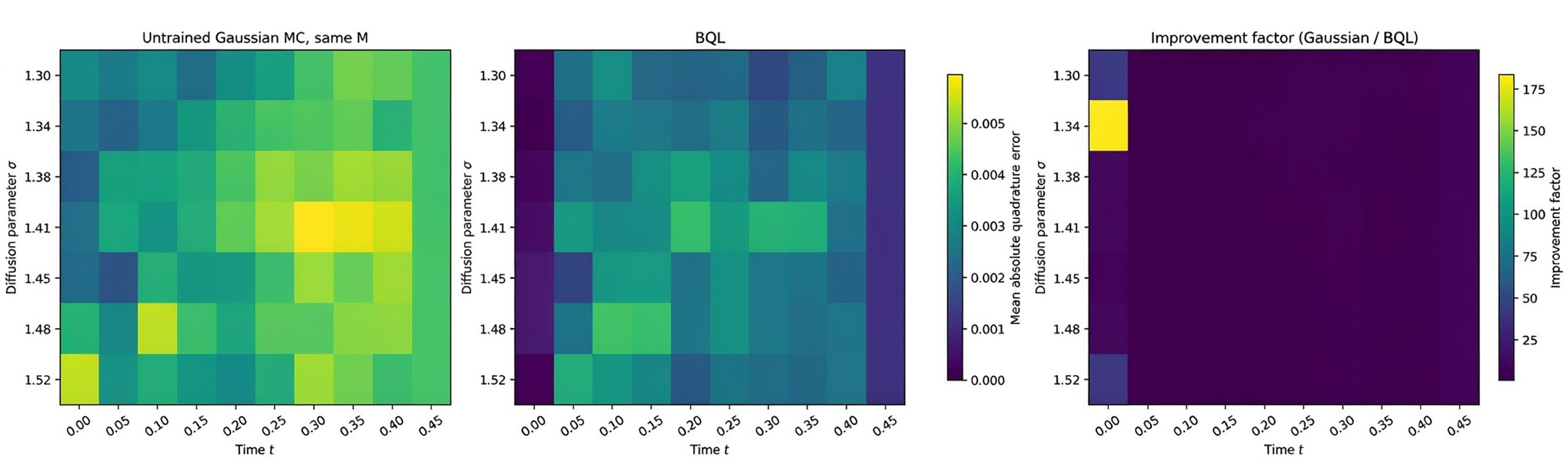}
\end{minipage}

\vspace{0.25cm}

\begin{minipage}[t]{0.96\textwidth}
\centering
\includegraphics[width=\linewidth]{HJB_pictures/HJB_MC2_d100.jpg}
\end{minipage}

\caption{Heatmap comparison of MAEs obtained by untrained Gaussian MC quadrature and BQL across different diffusion-parameter values $\sigma$ and discrete time points $t$. The upper and lower rows correspond to $d=20$ and $d=100$, respectively. From left to right, the panels show the results obtained by untrained Gaussian MC quadrature, BQL, and the corresponding improvement factors (Gaussian/BQL).}
\label{HJB_figure3-2}
\end{figure}

As shown in Figure~\ref{HJB_figure3-1}, BQL achieves lower mean absolute quadrature errors at all tested time points for both $d=20$ and $d=100$. For $d=20$, the Gaussian MC error is approximately 5.8 times that of BQL at $t=0.30$. For $d=100$, the corresponding ratios are approximately 2.1 and 4.2 at $t=0.30$ and $t=0.45$, respectively. Figure~\ref{HJB_figure3-2} further confirms this advantage across most parameter--time combinations. These results demonstrate that, under a fixed node budget, optimizing node locations enables BQL to use the available quadrature nodes more effectively and substantially improve conditional-expectation approximation accuracy.

\section{Conclusion}
\label{sec:conclusion}

In this paper, we proposed the BQL method for solving high dimensional parametric semilinear parabolic PDEs through BSDEs. Inspired by the LRV approach, BQL incorporates learnable stochastic quadrature into the backward recursion of BSDEs and introduces a continuation surrogate to facilitate information transfer across time levels while avoiding nested quadrature. Together, these components provide a unified framework for local conditional expectation approximation and backward information transfer over the joint parameter--state domain.

From a theoretical perspective, we established a global error bound that characterizes the contributions of local quadrature, continuation approximation, and training errors. Based on this estimate, we proved first-order mean-square convergence under suitable approximation and integrability assumptions. Numerical experiments on three high dimensional parametric PDE benchmarks, with state dimensions ranging from 10 to 100, demonstrated that BQL achieves competitive solution accuracy and effectively approximates local conditional expectations. Moreover, BQL substantially reduces the amortized computational cost per query in many-query settings, with clear cost advantages over FBSNN and P$^2$INN when the number of queries exceeds 5.

Future work will focus on extending the theoretical analysis to broader continuation architectures, improving approximation and training strategies over higher-dimensional parameter--state domains, and enhancing the scalability of BQL for more challenging nonlinear parametric PDEs.

\appendix
\renewcommand{\thesubsection}{\thesection.\arabic{subsection}}
\renewcommand{\thesubsubsection}{\thesubsection.\arabic{subsubsection}}
\renewcommand{\theequation}{\thesection-\arabic{equation}}

\section{Time Discretization and Backward Stability}
\label{app:stability}
This appendix supplies the time-discretization estimate and the one-step stability estimates used in Theorem~\ref{thm:total_error}. The evaluation law, auxiliary functions, and local errors are those of Section~\ref{subsec:analytical_setting}. The trained functions are held fixed throughout. The proof of Theorem~\ref{thm:total_error}, including the global summation argument, is given below.

\subsection{Time-discretization estimate}
\label{app:time_discretization}
For each fixed parameter $p\in\mathcal P$ and deterministic initial state $x_0\in\mathbb R^d$, the time-discretization error estimate for Markovian BSDEs, applied to the explicit backward recursion \eqref{eq:YZ_conditional}~\cite{Zhang2004, BenderZhang2008}, yields, under the uniform regularity assumptions in Section~4.1 and for a sufficiently fine uniform partition $\pi$,
\begin{equation}
\begin{aligned}
&\max_{0\le i\le N}\mathbb E^{p,x_0}\!\left[|Y_{t_i}^{p,x_0}-Y_i^{p,x_0,\pi}|^2\right]+\sum_{i=0}^{N-1}\mathbb E^{p,x_0}\!\left[\int_{t_i}^{t_{i+1}}|Z_t^{p,x_0}-Z_i^{p,x_0,\pi}|^2\,dt\right]\le C(1+|x_0|^2)|\pi|.
\label{eq:time_discretization_fixed}
\end{aligned}
\end{equation}
Here, $\mathbb E^{p,x_0}$ denotes expectation for the system with parameter $p$ and initial state $x_0$. The estimate follows from Bender and Zhang~\cite[Theorem 2.3]{BenderZhang2008}, specialized to the decoupled forward equation. In this case, the Markovian iteration stabilizes, and the iteration error vanishes. The constant $C$ is independent of $p$, $x_0$, and $\pi$ under the uniform regularity assumptions in Section~\ref{subsec:analytical_setting}.

We now extend this fixed-parameter estimate to the parameterized system. Since $P$ and $X_0^P$ are $\mathcal F_0$-measurable, the uniform estimate \eqref{eq:time_discretization_fixed} can be applied conditionally at $(p,x_0)=(P,X_0^P)$, yielding almost surely,
\[
\begin{aligned}
&\max_{0\le i\le N}\mathbb E\!\left[|Y_{t_i}^P-Y_i^{P,\pi}|^2\mid P,X_0^P\right]+\sum_{i=0}^{N-1}\mathbb E\!\left[\left.\int_{t_i}^{t_{i+1}}|Z_t^P-Z_i^{P,\pi}|^2\,dt\right|P,X_0^P\right]\le C(1+|X_0^P|^2)|\pi|.
\end{aligned}
\]
Taking expectations and using the tower property of conditional expectation, we obtain
\[
\begin{aligned}
&\mathbb E\!\left[\max_{0\le i\le N}\mathbb E\!\left[|Y_{t_i}^P-Y_i^{P,\pi}|^2\mid P,X_0^P\right]\right]+\sum_{i=0}^{N-1}\mathbb E\!\left[\int_{t_i}^{t_{i+1}}|Z_t^P-Z_i^{P,\pi}|^2\,dt\right]\le C\bigl(1+\mathbb E|X_0^P|^2\bigr)|\pi|.
\end{aligned}
\]
Moreover, since the maximum of expectations does not exceed the expectation of the maximum, we have
\[
\begin{aligned}
&\max_{0\le i\le N}\mathbb E\!\left[|Y_{t_i}^P-Y_i^{P,\pi}|^2\right]+\sum_{i=0}^{N-1}\mathbb E\!\left[\int_{t_i}^{t_{i+1}}|Z_t^P-Z_i^{P,\pi}|^2\,dt\right] \le C\bigl(1+\mathbb E|X_0^P|^2\bigr)|\pi|.
\end{aligned}
\]

Finally, the initial-state moment assumption $\mathbb E|X_0^P|^2<\infty$ allows the factor $1+\mathbb E|X_0^P|^2$ to be absorbed into the constant. Consequently, the parameterized time-discretization error satisfies
\begin{equation}
\begin{aligned}
&\max_{0\le i\le N}\mathbb E\!\left[|Y_{t_i}^P-Y_i^{P,\pi}|^2\right]+\sum_{i=0}^{N-1}\mathbb E\!\left[\int_{t_i}^{t_{i+1}}|Z_t^P-Z_i^{P,\pi}|^2\,dt\right]\le C|\pi|,
\end{aligned}
\end{equation}
where $C$ satisfies the uniformity convention stated in Section~\ref{subsec:analytical_setting}.

\subsection{One-step Y and Z estimates}
\label{app:one_step_stability}

We now establish the one-step $Y$ and $Z$ estimates used in the proof of Theorem~\ref{thm:total_error}. The auxiliary functions $\bar y_i$ and $\bar z_i$ and the local errors $r_{i,Z}$, $r_{i,Y}$, and $r_{i,\mathrm{con}}$ are those defined in Section~\ref{subsec:analytical_setting}. Recall the definition of conditional fluctuation quantity
\begin{equation}
\mathcal V_i^{\mathrm{fluc}}
:=\mathbb E[
| U_{i+1} - \mathbb E_i[U_{i+1}] |^2],
\label{eq:app_conditional_fluctuation}
\end{equation}
where $U_{i+1}:=(y_{i+1}^{\pi}-\widehat y_{i+1}^{\mathrm{con}})(Q_{i+1})$.
By the orthogonality property of conditional expectation, we have
\begin{equation}
\mathbb E[|\mathbb E_i\left[
U_{i+1}\right]|^2]=
\mathbb E\left[|U_{i+1}|^2\right]
-\mathcal V_i^{\mathrm{fluc}}.
\label{eq:app_conditional_decomposition}
\end{equation}

\begin{lemma}[One-step $Z$-error]
\label{lem:app_Z_error}
For each $i<N$ and every $\lambda>0$,
\begin{equation}
\Delta t_i
\|z_i^\pi-\bar z_i\|_{L^2(\rho_i)}^2
\le
\mathcal V_i^{\mathrm{fluc}},
\qquad
\Delta t_i
\|z_i^\pi-\widehat z_i^{\mathrm{BQL}}\|_{L^2(\rho_i)}^2
\le
(1+\lambda)\mathcal V_i^{\mathrm{fluc}}
+
(1+\lambda^{-1})\Delta t_i r_{i,Z}.
\label{eq:app_Z_stability}
\end{equation}
\end{lemma}

\begin{proof}
Let $\zeta_i=\frac{\Delta W_i}{\sqrt{\Delta t_i}}$.
Conditional on $\mathcal F_{t_i}$, the coordinates of $\zeta_i$ have zero mean and identity covariance. By the conditional-expectation representations of $z_i^\pi$ and $\bar z_i$, we have
\[
\sqrt{\Delta t_i}\,(z_i^\pi-\bar z_i)(Q_i)
=\mathbb E_i[U_{i+1}\zeta_i].
\]
Since $\mathbb E_i[\zeta_i]=0$, we can subtract the conditional mean of $U_{i+1}$ to obtain
\[
\sqrt{\Delta t_i}\,(z_i^\pi-\bar z_i)(Q_i)
=\mathbb E_i[\widetilde U_{i+1}\zeta_i],
\]
where $\widetilde U_{i+1}:=U_{i+1}-\mathbb E_i[U_{i+1}]$.
Since the coordinates of $\zeta_i$ form a conditionally orthonormal system, the conditional Bessel inequality yields
\[
|\mathbb E_i[\widetilde U_{i+1}\zeta_i]|^2
\le \mathbb E_i[|\widetilde U_{i+1}|^2].
\]
Taking unconditional expectations on both sides and applying the tower property of conditional expectation, we obtain
\[
\Delta t_i\|z_i^\pi-\bar z_i\|_{L^2(\rho_i)}^2
\le \mathbb E[|\widetilde U_{i+1}|^2]
=\mathcal V_i^{\mathrm{fluc}}.
\]
Finally, inserting $\bar z_i$ between $z_i^\pi$ and $\widehat z_i^{\mathrm{BQL}}$ and applying Young's inequality yields the second estimate in \eqref{eq:app_Z_stability}.
\end{proof}

\begin{lemma}[One-step $Y$-error]
\label{lem:app_Y_error}
Under the standing conditions in Section~\ref{subsec:analytical_setting}, there exist $C,c_0,h_0>0$, independent of $i$, the partition, and the trained functions, such that, for $\Delta t_i\le h_0$,
\begin{equation}
\begin{aligned}
&\mathbb E\left[
|(y_i^\pi-\widehat y_i^{\mathrm{con}})(Q_i)|^2
\right]
+c_0\mathcal V_i^{\mathrm{fluc}}\le
(1+C\Delta t_i)
\mathbb E\left[
|(y_{i+1}^\pi-\widehat y_{i+1}^{\mathrm{con}})(Q_{i+1})|^2
\right]+
C\Delta t_i r_{i,Z}
+\frac{C}{\Delta t_i}
(r_{i,Y}+r_{i,\mathrm{con}}).
\end{aligned}
\label{eq:app_refined_stability}
\end{equation}
\end{lemma}

\begin{proof}
Set $A_i:=(y_i^\pi-\bar y_i)(Q_i)$.
Subtracting the exact and auxiliary $Y$ updates gives
\[
A_i=\mathbb E_i[U_{i+1}]+\Delta t_i b_i,
\]
where
\[
b_i:=\mathbb E_i\Big[
f_P\big(t_i,X_i^{P,\pi},y_{i+1}^\pi(Q_{i+1}),z_i^\pi(Q_i)\big)
-f_P\big(t_i,X_i^{P,\pi},\widehat y_{i+1}^{\mathrm{con}}(Q_{i+1}),\widehat z_i^{\mathrm{BQL}}(Q_i)\big)
\Big].
\]
By conditional Jensen's inequality and the Lipschitz continuity of the generator,
\[
\mathbb E|b_i|^2
\le C_b\left(
\mathbb E|U_{i+1}|^2
+\|z_i^\pi-\widehat z_i^{\mathrm{BQL}}\|_{L^2(\rho_i)}^2
\right),
\]
where $C_b:=2L_f^2$.

For $\gamma>0$, applying Young's inequality to the above decomposition of $A_i$ gives
\[
\mathbb E|A_i|^2
\le (1+\gamma\Delta t_i)\mathbb E|\mathbb E_i[U_{i+1}]|^2
+\left(\frac{\Delta t_i}{\gamma}+(\Delta t_i)^2\right)\mathbb E|b_i|^2.
\]
Using the preceding Lipschitz estimate, we further obtain
\[
\begin{aligned}
\mathbb E|A_i|^2
\le{}&
(1+\gamma\Delta t_i)\mathbb E|\mathbb E_i[U_{i+1}]|^2+C_b\left(\frac{\Delta t_i}{\gamma}+(\Delta t_i)^2\right)
\left(
\mathbb E|U_{i+1}|^2
+\|z_i^\pi-\widehat z_i^{\mathrm{BQL}}\|_{L^2(\rho_i)}^2
\right).
\end{aligned}
\]
By the conditional decomposition \eqref{eq:app_conditional_decomposition},
\[
\mathbb E|\mathbb E_i[U_{i+1}]|^2
=\mathbb E|U_{i+1}|^2-\mathcal V_i^{\mathrm{fluc}}.
\]
Moreover, Lemma~\ref{lem:app_Z_error} with $\lambda=1$ yields
\[
\|z_i^\pi-\widehat z_i^{\mathrm{BQL}}\|_{L^2(\rho_i)}^2
\le
\frac{2}{\Delta t_i}\mathcal V_i^{\mathrm{fluc}}+2r_{i,Z}.
\]
Substituting these estimates gives
\[
\begin{aligned}
\mathbb E|A_i|^2
\le{}&
(1+\gamma\Delta t_i)
\left(\mathbb E|U_{i+1}|^2-\mathcal V_i^{\mathrm{fluc}}\right)+C_b\left(\frac{\Delta t_i}{\gamma}+(\Delta t_i)^2\right)
\left[
\mathbb E|U_{i+1}|^2
+\frac{2}{\Delta t_i}\mathcal V_i^{\mathrm{fluc}}
+2r_{i,Z}
\right].
\end{aligned}
\]
Collecting the coefficients of $\mathbb E|U_{i+1}|^2$, $\mathcal V_i^{\mathrm{fluc}}$, and $r_{i,Z}$, respectively, gives
\begin{equation}
\begin{aligned}
\mathbb E|A_i|^2
\le{}&
\left[
1+\gamma\Delta t_i
+C_b\left(\frac{\Delta t_i}{\gamma}+(\Delta t_i)^2\right)
\right]\mathbb E|U_{i+1}|^2-\left[
1-\frac{2C_b}{\gamma}
+(\gamma-2C_b)\Delta t_i
\right]\mathcal V_i^{\mathrm{fluc}}\\
&+2C_b\left(\frac{\Delta t_i}{\gamma}+(\Delta t_i)^2\right)r_{i,Z}.
\end{aligned}
\label{eq:app_auxiliary_Y}
\end{equation}

Choose $\gamma>2C_b$ and restrict $\Delta t_i\le1$. With $c_*:=1-\frac{2C_b}{\gamma}>0$, we obtain
\[
\mathbb E|A_i|^2+c_*\mathcal V_i^{\mathrm{fluc}}
\le
(1+C\Delta t_i)\mathbb E|U_{i+1}|^2
+C\Delta t_i r_{i,Z}.
\]
Here, $C>0$ denotes a generic constant, whose value may change from line to line. It is independent of $i$, the time partition, and the trained functions, but may depend on the fixed choice of $\gamma$ and $C_b$.

The remaining discrepancy between the full value error and $A_i$ is
\[
\begin{aligned}
D_i
&:=(y_i^\pi-\widehat y_i^{\mathrm{con}})(Q_i)-A_i=(\bar y_i-\widehat y_i^{\mathrm{BQL}})(Q_i)
+(\widehat y_i^{\mathrm{BQL}}-\widehat y_i^{\mathrm{con}})(Q_i).
\end{aligned}
\]
Using $|a+b|^2\le2|a|^2+2|b|^2$ and the definitions of $r_{i,Y}$ and $r_{i,\mathrm{con}}$, we obtain
\[
\begin{aligned}
\mathbb E|D_i|^2
&\le
2\|\bar y_i-\widehat y_i^{\mathrm{BQL}}\|_{L^2(\rho_i)}^2
+2\|\widehat y_i^{\mathrm{BQL}}-\widehat y_i^{\mathrm{con}}\|_{L^2(\rho_i)}^2=2(r_{i,Y}+r_{i,\mathrm{con}}).
\end{aligned}
\]
Since $(y_i^\pi-\widehat y_i^{\mathrm{con}})(Q_i)=A_i+D_i$, Young's inequality with parameter $\Delta t_i$ gives
\[
\begin{aligned}
\mathbb E\left[|(y_i^\pi-\widehat y_i^{\mathrm{con}})(Q_i)|^2\right]
&\le (1+\Delta t_i)\mathbb E|A_i|^2
+\left(1+\frac{1}{\Delta t_i}\right)\mathbb E|D_i|^2\\
&\le (1+\Delta t_i)\mathbb E|A_i|^2
+2\left(1+\frac{1}{\Delta t_i}\right)(r_{i,Y}+r_{i,\mathrm{con}}).
\end{aligned}
\]
For $\Delta t_i\le1$, we have $1+\Delta t_i^{-1}\le2\Delta t_i^{-1}$, and hence
\[
\begin{aligned}
\mathbb E\left[|(y_i^\pi-\widehat y_i^{\mathrm{con}})(Q_i)|^2\right]
\le
(1+\Delta t_i)\mathbb E|A_i|^2
+\frac{4}{\Delta t_i}(r_{i,Y}+r_{i,\mathrm{con}}).
\end{aligned}
\]
Combining this estimate with the preceding bound for $\mathbb E|A_i|^2+c_*\mathcal V_i^{\mathrm{fluc}}$, we obtain
\[
\begin{aligned}
&\mathbb E\left[|(y_i^\pi-\widehat y_i^{\mathrm{con}})(Q_i)|^2\right]
+(1+\Delta t_i)c_*\mathcal V_i^{\mathrm{fluc}}\\
&\le
(1+\Delta t_i)(1+C\Delta t_i)\mathbb E|U_{i+1}|^2
+C(1+\Delta t_i)\Delta t_i r_{i,Z}
+\frac{4}{\Delta t_i}(r_{i,Y}+r_{i,\mathrm{con}}).
\end{aligned}
\]
Since $\Delta t_i\le1$, we have
$(1+\Delta t_i)(1+C\Delta t_i)\le1+C'\Delta t_i$
and $(1+\Delta t_i)\Delta t_i\le2\Delta t_i$, where $C'$ is independent of the time partition. Moreover, $(1+\Delta t_i)c_*\ge c_*$. Enlarging the generic constant $C$ therefore yields
\begin{equation}
\begin{aligned}
&\mathbb E\left[|(y_i^\pi-\widehat y_i^{\mathrm{con}})(Q_i)|^2\right]
+c_*\mathcal V_i^{\mathrm{fluc}}\le
(1+C\Delta t_i)\mathbb E|U_{i+1}|^2
+C\Delta t_i r_{i,Z}
+\frac{C}{\Delta t_i}(r_{i,Y}+r_{i,\mathrm{con}}).
\end{aligned}
\end{equation}
Thus, \eqref{eq:app_refined_stability} holds with $c_0=c_*$, where $c_0$ is independent of $i$, the time partition, and the trained functions.
\end{proof}

Dropping the nonnegative term $c_0\mathcal V_i^{\mathrm{fluc}}$ in \eqref{eq:app_refined_stability} gives the one-step value-error estimate \eqref{eq:one_step_Y_error} used in the proof of Theorem~\ref{thm:total_error}. Retaining this term yields the cumulative fluctuation estimate needed to close the global $Z$-error bound.

\section{Supplementary Local-Risk Arguments}
\label{app:local_proofs}
Section~\ref{subsec:local_BQL_error} already proves the local BQL bound, including the regularity and variance estimates. This appendix records only the two reference-sampling identities used there and the continuation-risk comparison deferred from Section~\ref{subsec:continuation_error}. All stage functions are fixed when the corresponding risks are evaluated.

\subsection{Finite-node benchmark and reference-risk identity}
\label{app:quadrature_identities}

Fix $i\in\{0,\ldots,N-1\}$ and $\alpha\in\{Z,Y\}$. We use the notation $m_i^\alpha$, $V_i^\alpha$, and $\rho_i$ introduced in Section~\ref{subsec:local_BQL_error}. For a candidate configuration $\vartheta_i^\alpha=(\vartheta_{i,1}^\alpha,\ldots,\vartheta_{i,M}^\alpha)\in(\mathbb R^d)^M$, the finite-node approximation and its population risk are
\[
\mathcal E_{i,M}^\alpha(q;\vartheta_i^\alpha):=\frac{1}{M}\sum_{m=1}^M F_i^\alpha(q,\vartheta_{i,m}^\alpha),\qquad R_i^\alpha(\vartheta_i^\alpha):=\|m_i^\alpha-\mathcal E_{i,M}^\alpha(\cdot;\vartheta_i^\alpha)\|_{L^2(\rho_i)}^2.
\]

\paragraph{Finite-node benchmark.}
We first establish an upper bound for the best population risk achievable with $M$ quadrature nodes by comparing it with the risk of an auxiliary MC configuration. Let $\vartheta_{i,1}^{\alpha,\mathrm{aux}},\ldots,\vartheta_{i,M}^{\alpha,\mathrm{aux}}$ be independent random vectors with distribution $\mathcal N(0,\Delta t_i I_d)$, independent of $Q_i$, and write $\vartheta_i^{\alpha,\mathrm{aux}}=(\vartheta_{i,1}^{\alpha,\mathrm{aux}},\ldots,\vartheta_{i,M}^{\alpha,\mathrm{aux}})$. For each fixed $q$, independence and the MC variance identity \cite[Chapter~2]{Owen2013} yield
\[
\begin{aligned}
&\mathbb E_{\mathrm{aux}}\!\left[|\mathcal E_{i,M}^\alpha(q;\vartheta_i^{\alpha,\mathrm{aux}})-m_i^\alpha(q)|^2\right]=\frac{1}{M}\mathbb E_\vartheta\!\left[|F_i^\alpha(q,\vartheta)-m_i^\alpha(q)|^2\right].
\end{aligned}
\]
Indeed, the centered summands have zero mean, so all cross terms vanish after taking expectation over the independent auxiliary nodes. Integrating with respect to $\rho_i$ and interchanging expectations gives
\[
\mathbb E_{\mathrm{aux}}\!\left[R_i^\alpha(\vartheta_i^{\alpha,\mathrm{aux}})\right]=\frac{V_i^\alpha}{M}.
\]
Since the infimum over deterministic candidate configurations cannot exceed the average risk of these random configurations, we obtain
\[
\inf_{\vartheta_i^\alpha\in(\mathbb R^d)^M}R_i^\alpha(\vartheta_i^\alpha)\le\mathbb E_{\mathrm{aux}}\!\left[R_i^\alpha(\vartheta_i^{\alpha,\mathrm{aux}})\right]=\frac{V_i^\alpha}{M}.
\]
Thus, the best deterministic $M$-node configuration achieves a population risk no larger than the MC benchmark $V_i^\alpha/M$. 

\paragraph{Reference-risk identity.}
We next relate the population training risk based on finitely many reference nodes to the population risk defined using the exact conditional mean. Recall the reference target $ m_{i,K_\alpha}^{\alpha,\mathrm{ref}}$ defined in Section \ref{subsec:local_BQL_error}, constructed from $K_\alpha$ independent Gaussian reference vectors with distribution $\mathcal N(0,\Delta t_i I_d)$. For each fixed $q$, this target is unbiased for $m_i^\alpha(q)$, and
\[
\mathbb E_{\mathrm{ref}}\!\left[| m_{i,K_\alpha}^{\alpha,\mathrm{ref}}(q)-m_i^\alpha(q)|^2\right]=\frac{1}{K_\alpha}\mathbb E_\vartheta\!\left[|F_i^\alpha(q,\vartheta)-m_i^\alpha(q)|^2\right].
\]
For a fixed candidate configuration $\vartheta_i^\alpha$, the finite-node approximation is independent of the reference vectors. We therefore expand the squared loss around the exact conditional mean $m_i^\alpha(q)$ and take expectation over the reference vectors. The cross term vanishes by the unbiasedness of the reference target, giving
\[
\begin{aligned}
&\mathbb E_{\mathrm{ref}}\!\left[|\mathcal E_{i,M}^\alpha(q;\vartheta_i^\alpha)- m_{i,K_\alpha}^{\alpha,\mathrm{ref}}(q)|^2\right]=|\mathcal E_{i,M}^\alpha(q;\vartheta_i^\alpha)-m_i^\alpha(q)|^2+\frac{1}{K_\alpha}\mathbb E_\vartheta\!\left[|F_i^\alpha(q,\vartheta)-m_i^\alpha(q)|^2\right].
\end{aligned}
\]
Integrating with respect to $\rho_i$ yields the reference-risk identity
\[
R_{i,K_\alpha}^\alpha(\vartheta_i^\alpha)=R_i^\alpha(\vartheta_i^\alpha)+\frac{V_i^\alpha}{K_\alpha}.
\]
The identity shows that the reference-node sampling introduces an additional variance term $V_i^\alpha/K_\alpha$ into the population training risk. Since this term is independent of the candidate configuration, it shifts all candidate risks by the same amount and therefore cancels in the excess-risk comparison:
\[
R_{i,K_\alpha}^\alpha(\vartheta_i^\alpha)-\inf_{\widetilde\vartheta_i^\alpha\in(\mathbb R^d)^M}R_{i,K_\alpha}^\alpha(\widetilde\vartheta_i^\alpha)=R_i^\alpha(\vartheta_i^\alpha)-\inf_{\widetilde\vartheta_i^\alpha\in(\mathbb R^d)^M}R_i^\alpha(\widetilde\vartheta_i^\alpha).
\]

\paragraph{Local quadrature error bound.}
We now combine the preceding two results with the generalization discrepancy and optimization residual defined in the main text. By the definition of $\varepsilon_{i,\alpha}^{\mathrm{gen}}$, the population and empirical risks satisfy 
\[
\begin{aligned}
R_{i,K_\alpha}^\alpha(\widehat\vartheta_i^\alpha)
\le \widehat R_{i,H,K_\alpha}^\alpha(\widehat\vartheta_i^\alpha)+\varepsilon_{i,\alpha}^{\mathrm{gen}},\qquad
\inf_{\vartheta_i^\alpha\in(\mathbb R^d)^M}
\widehat R_{i,H,K_\alpha}^\alpha(\vartheta_i^\alpha)
\le
\inf_{\vartheta_i^\alpha\in(\mathbb R^d)^M}
R_{i,K_\alpha}^\alpha(\vartheta_i^\alpha)+\varepsilon_{i,\alpha}^{\mathrm{gen}}.
\end{aligned}
\]
Using the definition of $\delta_{i,\alpha}^{\mathrm{opt}}$, we obtain
\[
\begin{aligned}
R_{i,K_\alpha}^\alpha(\widehat\vartheta_i^\alpha)
&\le
\widehat R_{i,H,K_\alpha}^\alpha(\widehat\vartheta_i^\alpha)
+\varepsilon_{i,\alpha}^{\mathrm{gen}}=
\inf_{\vartheta_i^\alpha\in(\mathbb R^d)^M}
\widehat R_{i,H,K_\alpha}^\alpha(\vartheta_i^\alpha)
+\delta_{i,\alpha}^{\mathrm{opt}}
+\varepsilon_{i,\alpha}^{\mathrm{gen}}\\
&\le
\inf_{\vartheta_i^\alpha\in(\mathbb R^d)^M}
R_{i,K_\alpha}^\alpha(\vartheta_i^\alpha)
+2\varepsilon_{i,\alpha}^{\mathrm{gen}}
+\delta_{i,\alpha}^{\mathrm{opt}}.
\end{aligned}
\]
By the reference-risk identity, the excess risk above is equal to that based on $R_i^\alpha$. Combining this equality with the finite-node benchmark, we obtain
\[
\begin{aligned}
r_{i,\alpha}
&=R_i^\alpha(\widehat\vartheta_i^\alpha)\le
\inf_{\vartheta_i^\alpha\in(\mathbb R^d)^M}
R_i^\alpha(\vartheta_i^\alpha)
+2\varepsilon_{i,\alpha}^{\mathrm{gen}}
+\delta_{i,\alpha}^{\mathrm{opt}}\le
\frac{V_i^\alpha}{M}
+2\varepsilon_{i,\alpha}^{\mathrm{gen}}
+\delta_{i,\alpha}^{\mathrm{opt}},
\end{aligned}
\]
which gives the local quadrature error bound in the main text.

The reference-risk identity is a pointwise identity for every fixed candidate configuration and can therefore be evaluated at the trained configuration. It does not require the trained configuration to be independent of the training data, since its dependence on the training data is accounted for by the uniform generalization discrepancy and the optimization residual.

\subsection{Continuation-risk comparison}
\label{app:continuation-risk-comparison}

\noindent\textit{Proof of Proposition~\ref{prop:continuation_error_decomposition}.}
Apply the standard empirical-risk comparison \cite[Chapter~4]{ShalevShwartzBenDavid2014} to the fixed target $\widehat y_i^{\mathrm{BQL}}$ and the continuation function class $\mathcal H_i^{\mathrm{con}}$. By the definitions of the generalization and optimization residuals,
\begin{equation}
\begin{aligned}
r_{i,\mathrm{con}}
&=R_i^{\mathrm{con}}(\widehat y_i^{\mathrm{con}})\\
&\le\widehat R_{i,N_{\mathrm{con}}}^{\mathrm{con}}(\widehat y_i^{\mathrm{con}})+\varepsilon_{i,\mathrm{con}}^{\mathrm{gen}}\\
&=\inf_{v\in\mathcal H_i^{\mathrm{con}}}\widehat R_{i,N_{\mathrm{con}}}^{\mathrm{con}}(v)+\delta_{i,\mathrm{con}}^{\mathrm{opt}}+\varepsilon_{i,\mathrm{con}}^{\mathrm{gen}}\\
&\le\inf_{v\in\mathcal H_i^{\mathrm{con}}}R_i^{\mathrm{con}}(v)+2\varepsilon_{i,\mathrm{con}}^{\mathrm{gen}}+\delta_{i,\mathrm{con}}^{\mathrm{opt}}\\
&=(\eta_{i,\mathrm{con}}^{\mathrm{app}})^2+2\varepsilon_{i,\mathrm{con}}^{\mathrm{gen}}+\delta_{i,\mathrm{con}}^{\mathrm{opt}}.
\end{aligned}
\label{eq:continuation_risk_comparison_appendix}
\end{equation}
The two generalization contributions arise from comparing the empirical and population risks at the trained surrogate and at their respective infima. The final equality follows from the definition of the best-approximation error. No population-risk minimizer is required, since the argument uses infima throughout.
\hfill$\square$

\section{Architecture-Dependent Continuation Approximation}
\label{app:architecture}
This appendix establishes architecture-dependent approximation error bounds for the continuation surrogate under the stated uniform regularity and moment assumptions. These assumptions are imposed on the continuation target functions and the associated forward states, and are not asserted to follow automatically from the network architecture or the training procedure. The resulting estimates characterize the best-approximation capabilities of the admissible function classes under these conditions, rather than providing unconditional guarantees on the regularity of the trained networks or their training errors.

\paragraph{Common Notation and Standing Assumptions.}
We consider a sequence of uniform time partitions $\{\pi_N\}_{N\geq1}$ of $[0,T]$, where $\pi_N=\{t_i=iT/N\}_{i=0}^N$ and $\Delta t=T/N$. Throughout this appendix, the notation introduced below is understood as the specialization of the corresponding notation in the preceding sections to the partition $\pi_N$. In particular, for each fixed $N$, $\pi$ is identified with $\pi_N$, and we write $X_i^{P,\pi_N}$ and $Q_{i,N}=(P,X_i^{P,\pi_N})$ in place of $X_i^{P,\pi}$ and $Q_i$, respectively. Similarly, $\widehat y_{i,N}^{\mathrm{BQL}}$ denotes the continuation target $\widehat y_i^{\mathrm{BQL}}$ associated with the partition $\pi_N$.
Let $d_q=d_p+d$. For each $N\geq1$, $i=0,\ldots,N-1$, and $R\geq1$, define $D_R=\operatorname{int}(\mathcal P)\times(-R,R)^d$ and $\Omega_R=\overline{D_R}$.

Assume that the continuation targets are uniformly Lipschitz continuous with respect to the joint parameter--state variable: there exists a constant $L_y>0$, independent of $N$ and $i$, such that
\begin{equation}
\sup_{N\geq1}\sup_{0\leq i<N}\operatorname{Lip}\big(\widehat y_{i,N}^{\mathrm{BQL}};\mathcal P\times\mathbb R^d\big)\leq L_y.
\label{eq:uniform_lipschitz_continuation}
\end{equation}

Assume also that there exist $\delta>0$ and constants $C_{\mathrm{mom}},C_u>0$, independent of $N$ and $i$, such that
\begin{equation}
\sup_{N\geq1}\max_{0\leq i\leq N}\mathbb E[(1+|X_i^{P,\pi_N}|)^{2+\delta}]\leq C_{\mathrm{mom}},\qquad |\widehat y_{i,N}^{\mathrm{BQL}}(p,x)|\leq C_u(1+|x|),
\label{eq:FNN_uniform_tail_assumptions}
\end{equation}
for all $p\in\mathcal P$, $x\in\mathbb R^d$, $N\geq1$, and $0\leq i<N$.

\begin{lemma}[FNN approximation error bound]
\label{lem:FNN_approximation_bound}
Define $G_{i,N}(R):=\|\widehat y_{i,N}^{\mathrm{BQL}}\|_{W^{1,\infty}(D_R)}$,
where $\|\phi\|_{W^{1,\infty}(D_R)}:=\max_{|\nu|\le1}\operatorname*{ess\,sup}_{q\in D_R}|D^\nu\phi(q)|$. Here, $\nu=(\nu_1,\ldots,\nu_{d_q})\in\mathbb N_0^{d_q}$, $|\nu|:=\sum_{j=1}^{d_q}\nu_j$, and $D^\nu f$ denotes the corresponding weak derivative, with $D^0f:=f$. Under the common assumptions, there exists a constant $C_G>0$, independent of $N$, $i$, and $R$, such that
\[
G_{i,N}(R)\le C_G(1+R),\qquad R\ge1.
\]

Let $\eta_{i,N,\mathrm{FNN}}^{\mathrm{app}}$ denote the best-approximation error of $\widehat y_{i,N}^{\mathrm{BQL}}$ in the admissible FNN class $\mathcal H_{i,N}^{\mathrm{FNN}}(S_{i,N})$ with network size bound $S_{i,N}$, which bounds both the total number of weights and the number of computation units. Then, for sufficiently large $S_{i,N}$,
\begin{equation}
(\eta_{i,N,\mathrm{FNN}}^{\mathrm{app}})^2
\le C(1+R)^4
\left(\frac{S_{i,N}}{\log(1+S_{i,N})}\right)^{-2/d_q}
+C_\tau R^{-\delta},\qquad R\ge1,
\label{eq:FNN_approximation_bound}
\end{equation}
where $C$ is independent of $N$, $i$, $R$, and $S_{i,N}$, and $C_\tau$ is independent of $N$, $i$, and $R$.
\end{lemma}

\begin{proof}
By the uniform Lipschitz assumption in the common notation, there exists a constant $L_y>0$, independent of $N$ and $i$, such that
\[
|\widehat y_{i,N}^{\mathrm{BQL}}(q)-\widehat y_{i,N}^{\mathrm{BQL}}(q')|\le L_y|q-q'|,\qquad q,q'\in\mathcal P\times\mathbb R^d.
\]
By Rademacher's theorem~\cite{EvansGariepy2015}, the weak first derivatives exist almost everywhere and satisfy
\[
\operatorname*{ess\,sup}_{q\in D_R}|D^\nu\widehat y_{i,N}^{\mathrm{BQL}}(q)|\le L_y,\qquad |\nu|=1.
\]
Moreover, the uniform linear-growth condition and $|x|\le\sqrt d R$ on $\Omega_R$ imply
\[
\|\widehat y_{i,N}^{\mathrm{BQL}}\|_{L^\infty(D_R)}\le C_u(1+\sqrt d R)\le C(1+R).
\]
Consequently, by the definition of the Sobolev norm,
\[
G_{i,N}(R):=\|\widehat y_{i,N}^{\mathrm{BQL}}\|_{W^{1,\infty}(D_R)}\le C_G(1+R),
\]
where $C_G$ is independent of $N$, $i$, and $R$. In particular, $G_{i,N}(R)<\infty$.

We first control the approximation error outside $\Omega_R$. Define $B_{i,N}(R):=1+\|\widehat y_{i,N}^{\mathrm{BQL}}\|_{L^\infty(\Omega_R)}$.
Let $v_{i,N}^{\mathrm{FNN}}$ be a candidate approximating network and denote its output clipped to $[-B_{i,N}(R),B_{i,N}(R)]$ by $\widetilde v_{i,N}^{\mathrm{FNN}}$. Since $|\widehat y_{i,N}^{\mathrm{BQL}}|\le B_{i,N}(R)$ on $\Omega_R$, clipping does not increase the approximation error there. Moreover, $|\widetilde v_{i,N}^{\mathrm{FNN}}|\le B_{i,N}(R)$ everywhere, and hence
\begin{equation}
\begin{aligned}
&\mathbb E\!\left[
|\widehat y_{i,N}^{\mathrm{BQL}}(Q_{i,N})-\widetilde v_{i,N}^{\mathrm{FNN}}(Q_{i,N})|^2
\mathbf 1_{\{Q_{i,N}\notin\Omega_R\}}
\right]\\
&\qquad\le
2\mathbb E\!\left[
|\widehat y_{i,N}^{\mathrm{BQL}}(Q_{i,N})|^2
\mathbf 1_{\{Q_{i,N}\notin\Omega_R\}}
\right]
+2B_{i,N}(R)^2\mathbb P(Q_{i,N}\notin\Omega_R).
\end{aligned}
\label{eq:FNN_tail_error_bound}
\end{equation}
Accordingly, define
\[
\tau_{i,N}(R):=
2\mathbb E\!\left[
|\widehat y_{i,N}^{\mathrm{BQL}}(Q_{i,N})|^2
\mathbf 1_{\{Q_{i,N}\notin\Omega_R\}}
\right]
+2B_{i,N}(R)^2\mathbb P(Q_{i,N}\notin\Omega_R).
\]

The linear-growth condition gives
\[
B_{i,N}(R)\le C_B(1+R),
\]
where $C_B$ is independent of $N$, $i$, and $R$. Moreover, since $P\in\mathcal P$ almost surely, Markov's inequality and the uniform $(2+\delta)$-moment bound yield
\[
\begin{aligned}
\mathbb P(Q_{i,N}\notin\Omega_R)
&\le \mathbb P(|X_i^{P,\pi_N}|>R)\\
&\le \frac{\mathbb E[(1+|X_i^{P,\pi_N}|)^{2+\delta}]}{(1+R)^{2+\delta}}\\
&\le C_{\mathrm{mom}}(1+R)^{-(2+\delta)}.
\end{aligned}
\]
The linear-growth condition further gives
\[
\begin{aligned}
&\mathbb E\!\left[
|\widehat y_{i,N}^{\mathrm{BQL}}(Q_{i,N})|^2
\mathbf 1_{\{Q_{i,N}\notin\Omega_R\}}
\right]\\
&\qquad\le C_u^2\mathbb E\!\left[
(1+|X_i^{P,\pi_N}|)^2
\mathbf 1_{\{|X_i^{P,\pi_N}|>R\}}
\right]\\
&\qquad\le C_u^2 C_X(1+R)^{-\delta}.
\end{aligned}
\]
Consequently,
\[
\tau_{i,N}(R)\le 2C_{\mathrm{mom}}(C_u^2+C_B^2)(1+R)^{-\delta}\le C_\tau R^{-\delta},\qquad R\ge1,
\]
where $C_\tau:=2C_X(C_u^2+C_B^2)$ is independent of $N$, $i$, and $R$.

Next, we control the approximation error inside $\Omega_R$. If $G_{i,N}(R)=0$, the target vanishes almost everywhere on $\Omega_R$, and the desired bound follows by choosing the zero network and applying the truncation estimate. Otherwise, let $T_R:\Omega_R\to[0,1]^{d_q}$ be the affine bijection defined by
\[
(T_R(p,x))_j=\frac{p_j-a_j}{b_j-a_j},\quad 1\le j\le d_p,\qquad
(T_R(p,x))_{d_p+k}=\frac{x_k+R}{2R},\quad 1\le k\le d,
\]
where $\mathcal P=\prod_{j=1}^{d_p}[a_j,b_j]$. Define the normalized target $\overline y_{i,N,R}(z):=G_{i,N}(R)^{-1}\widehat y_{i,N}^{\mathrm{BQL}}(T_R^{-1}(z))$.
By the definition of $G_{i,N}(R)$, the uniform Lipschitz assumption, and the fixed geometry of $\mathcal P$, the affine transformation gives
\[
\|\overline y_{i,N,R}\|_{L^\infty((0,1)^{d_q})}\le1,
\qquad
\|\overline y_{i,N,R}\|_{W^{1,\infty}((0,1)^{d_q})}\le C_T(1+R),
\]
where $C_T$ is independent of $N$, $i$, and $R$. In particular, setting $A_R:=C_T(1+R)$
with $C_T$ sufficiently large, the further normalized function $h_{i,N,R}(z):=A_R^{-1}\overline y_{i,N,R}(z)$
satisfies $\|h_{i,N,R}\|_{W^{1,\infty}((0,1)^{d_q})}\le1$.

By Theorem~1 of Yarotsky~\cite{Yarotsky2017} with smoothness order $n=1$, for every sufficiently small $\varepsilon>0$, there exists a ReLU network approximating $h_{i,N,R}$ with uniform error at most $\varepsilon$ and at most $C_{d_q}\varepsilon^{-d_q}(1+\log(1/\varepsilon))$
trainable parameters, where $C_{d_q}$ is independent of $N$, $i$, and $R$. Rescaling the network output by $G_{i,N}(R)A_R$ and composing with $T_R$ gives an approximation to $\widehat y_{i,N}^{\mathrm{BQL}}$ on $\Omega_R$ with uniform error at most $C G_{i,N}(R)(1+R)\varepsilon$.
Inverting the network complexity estimate and accounting for the output-rescaling and clipping overhead, for all sufficiently large $S_{i,N}$ there exists an admissible clipped network $\widetilde v_{i,N}^{\mathrm{FNN}}$ with at most $S_{i,N}$ trainable parameters such that
\begin{equation}
\|\widehat y_{i,N}^{\mathrm{BQL}}-\widetilde v_{i,N}^{\mathrm{FNN}}\|_{L^\infty(\Omega_R)}
\le C G_{i,N}(R)(1+R)
\left(\frac{S_{i,N}}{\log(1+S_{i,N})}\right)^{-1/d_q},
\label{eq:FNN_uniform_approximation}
\end{equation}
where $C$ is independent of $N$, $i$, $R$, and $S_{i,N}$.

Using $\widetilde v_{i,N}^{\mathrm{FNN}}$ as a candidate in the definition of the best-approximation error and splitting the expectation over $\Omega_R$ and its complement, we obtain
\begin{equation}
\begin{aligned}
(\eta_{i,N,\mathrm{FNN}}^{\mathrm{app}})^2
&\le \mathbb E\!\left[
|\widehat y_{i,N}^{\mathrm{BQL}}(Q_{i,N})-\widetilde v_{i,N}^{\mathrm{FNN}}(Q_{i,N})|^2
\right]\\
&\le
\|\widehat y_{i,N}^{\mathrm{BQL}}-\widetilde v_{i,N}^{\mathrm{FNN}}\|_{L^\infty(\Omega_R)}^2
+\tau_{i,N}(R)\\
&\le
C G_{i,N}(R)^2(1+R)^2
\left(\frac{S_{i,N}}{\log(1+S_{i,N})}\right)^{-2/d_q}
+\tau_{i,N}(R).
\end{aligned}
\label{eq:FNN_approximation_with_tail}
\end{equation}
Finally, using $G_{i,N}(R)\le C_G(1+R)$ and $\tau_{i,N}(R)\le C_\tau R^{-\delta}$, we obtain
\[
(\eta_{i,N,\mathrm{FNN}}^{\mathrm{app}})^2
\le C(1+R)^4
\left(\frac{S_{i,N}}{\log(1+S_{i,N})}\right)^{-2/d_q}
+C_\tau R^{-\delta},
\]
which establishes the desired approximation bound.
\end{proof}

\begin{lemma}[Hybrid RBF--FNN approximation error bound]
\label{lem:hybrid_approximation_bound}
Let the hybrid continuation function class be defined by
\[
\mathcal H_{i,N}^{\mathrm{hyb}}(S_{i,N},N_c):=
\left\{
v^{\mathrm{RBF}}+v^{\mathrm{FNN}}:
v^{\mathrm{RBF}}\in\mathcal H_{i,N}^{\mathrm{RBF}}(N_c),\
v^{\mathrm{FNN}}\in\mathcal H_{i,N}^{\mathrm{FNN}}(S_{i,N})
\right\}.
\]
Here, $S_{i,N}$ denotes the network size bound for the FNN branch, which bounds both the total number of weights and the number of computation units, while $N_c$ denotes the number of centers in the RBF branch. By setting the output weights and bias of the RBF component to zero, the hybrid architecture contains the FNN class as a special case. Consequently, under the assumptions of Lemma~\ref{lem:FNN_approximation_bound}, for sufficiently large $S_{i,N}$,
\begin{equation}
(\eta_{i,N,\mathrm{con}}^{\mathrm{app}})^2
\leq
C(1+R)^4
\left(\frac{S_{i,N}}{\log(1+S_{i,N})}\right)^{-2/d_q}
+C_\tau R^{-\delta},
\qquad R\geq1,
\label{Feq:hybrid_approximation_bound}
\end{equation}
where $C$ and $C_\tau$ are independent of $N$, $i$, $S_{i,N}$, and $R$.
\end{lemma}

\begin{proof}
By setting the output weights and bias of the RBF component to zero, any function in $\mathcal H_{i,N}^{\mathrm{FNN}}(S_{i,N})$ can be represented as a member of the hybrid class. Hence,
\[
\mathcal H_{i,N}^{\mathrm{FNN}}(S_{i,N})
\subseteq
\mathcal H_{i,N}^{\mathrm{hyb}}(S_{i,N},N_c).
\]
It follows from the definition of the best-approximation errors that
\[
\eta_{i,N,\mathrm{con}}^{\mathrm{app}}
\leq
\eta_{i,N,\mathrm{FNN}}^{\mathrm{app}}.
\]
Squaring both sides and applying Lemma~\ref{lem:FNN_approximation_bound} yield
\[
(\eta_{i,N,\mathrm{con}}^{\mathrm{app}})^2
\leq
C(1+R)^4
\left(\frac{S_{i,N}}{\log(1+S_{i,N})}\right)^{-2/d_q}
+C_\tau R^{-\delta}.
\]
This completes the proof.
\end{proof}

\begin{remark}
The above estimate is obtained by considering a special case of the hybrid architecture in which the RBF component is set to zero, so that the approximation is provided entirely by the FNN component. It therefore establishes an upper bound for the best-approximation error of the hybrid class without explicitly exploiting the approximation capability of the RBF component. A sharper approximation analysis that accounts for the contribution of the RBF component and its interaction with the FNN component requires additional theoretical developments, which are left for future work.
\end{remark}

\section{Quadrature Accuracy Evaluation}
\label{app:quadrature_evaluation}

Let $N_{\mathrm{map}}$ denote the number of test trajectories used to evaluate the quadrature errors for each parameter--time configuration $(t_i,p_j)$. For the $r$-th test trajectory associated with parameter configuration $p_j$, let $X_i^{(r,j)}$ denote its state at time $t_i$. Consistent with the joint parameter--state variable $q=(p,x)$ introduced previously, we define the corresponding test input by $q_{i,j}^{(r)}=(p_j,X_i^{(r,j)}), r=1,\ldots,N_{\mathrm{map}}$.
Using the integrand $F_i^Y(q,\vartheta)$ defined previously, the conditional expectation to be approximated is
\[
I_{i,j}^{(r)}
=
\mathbb E_{\xi\sim\mathcal N(0,I_d)}
\left[
F_i^Y\left(q_{i,j}^{(r)},\sqrt{\Delta t_i}\xi\right)
\right].
\]

We consider three numerical approximations of this conditional expectation. First, the reference value is estimated using $K_{\mathrm{ref}}$ independent Gaussian samples:
\[
\widehat I_{i,j}^{\mathrm{ref},(r)}
=
\frac{1}{K_{\mathrm{ref}}}
\sum_{k=1}^{K_{\mathrm{ref}}}
F_i^Y\left(q_{i,j}^{(r)},\sqrt{\Delta t_i}\zeta_{i,k}\right),
\qquad
\zeta_{i,k}\sim\mathcal N(0,I_d).
\]
Second, the BQL approximation is computed using the $M$ learned Y-quadrature nodes:
\[
\widehat I_{i,j}^{\mathrm{BQL},(r)}
=
\frac{1}{M}
\sum_{m=1}^{M}
F_i^Y\left(q_{i,j}^{(r)},\vartheta_{i,m}^{Y,\mathrm{BQL}}\right).
\]
Third, the untrained Gaussian MC approximation is obtained using $M$ independent Gaussian nodes in each realization:
\[
\widehat I_{i,j}^{\mathrm{G},(s,r)}
=
\frac{1}{M}
\sum_{m=1}^{M}
F_i^Y\left(q_{i,j}^{(r)},\sqrt{\Delta t_i}\xi_{i,m}^{(s)}\right),
\qquad
\xi_{i,m}^{(s)}\sim\mathcal N(0,I_d),
\]
where $s=1,\ldots,N_{\mathrm{rep}}$ indexes independent Gaussian MC realizations. Thus, BQL and untrained Gaussian MC use the same number $M$ of quadrature nodes in each realization.

We report the mean absolute quadrature errors (MAEs) and the corresponding improvement factor (IF), defined as follows:
\begin{equation}
\begin{aligned}
\mathrm{MAE}_{i,j}^{\mathrm{BQL}}
=
\frac{1}{N_{\mathrm{map}}}
\sum_{r=1}^{N_{\mathrm{map}}}
\left|
\widehat I_{i,j}^{\mathrm{BQL},(r)}
-
\widehat I_{i,j}^{\mathrm{ref},(r)}
\right|, \quad
\mathrm{MAE}_{i,j}^{\mathrm{G}}
=
\frac{1}{N_{\mathrm{rep}}N_{\mathrm{map}}}
\sum_{s=1}^{N_{\mathrm{rep}}}
\sum_{r=1}^{N_{\mathrm{map}}}
\left|
\widehat I_{i,j}^{\mathrm{G},(s,r)}
-
\widehat I_{i,j}^{\mathrm{ref},(r)}
\right|.
\end{aligned}
\label{eq:quadrature_mae}
\end{equation}
The improvement factor is defined by
\begin{equation}
\mathrm{IF}_{i,j}
=
\frac{\mathrm{MAE}_{i,j}^{\mathrm{G}}}
{\mathrm{MAE}_{i,j}^{\mathrm{BQL}}+10^{-15}}.
\label{eq:quadrature_improvement_factor}
\end{equation}
Here, $N_{\mathrm{rep}}$ denotes the number of independent Gaussian MC realizations, which is set to $16$ in the default configuration, while $K_{\mathrm{ref}}$ is set to $50000$. The improvement factor measures the error reduction achieved by BQL relative to untrained Gaussian MC quadrature, with values greater than one indicating higher accuracy of BQL. The reported quadrature errors are evaluated relative to the high-sample MC reference estimates rather than the exact conditional expectations.

For the time-dependent error curves, we further average the MAEs over all tested parameter configurations. Specifically, letting $N_p$ denote the number of tested parameter configurations, the averaged errors at time $t_i$ are given by
\begin{equation}
\overline{\mathrm{MAE}}_{i}^{A}
=
\frac{1}{N_p}
\sum_{j=1}^{N_p}
\mathrm{MAE}_{i,j}^{A},
\qquad A\in\{\mathrm{BQL},\mathrm{G}\}.
\label{eq:averaged_quadrature_mae}
\end{equation}

\end{document}